\documentclass[a4paper, 9pt]{article}
\usepackage{latexsym}
\usepackage{amsmath}
\usepackage{graphics}
\usepackage{amsthm}
\usepackage{bm}
\usepackage{amssymb}
\usepackage{makeidx}
\usepackage{multicol}
\usepackage[all]{xy}
\newfont{\Fr}{eufm10}
\newfont{\Sc}{eusm10}
\newfont{\Bb}{msbm10}
\newfont{\Am}{msam10}
\newfont{\am}{msam7}
\numberwithin{equation}{section}
\newtheorem{theorem}{Theorem}[section]
\newtheorem{proposition}[theorem]{Proposition}
\newtheorem{lemma}[theorem]{Lemma}
\newtheorem{corollary}[theorem]{Corollary}
\newtheorem{claim}{Claim}{\bf}{\it}

\newtheorem{ftheorem}{Theorem}{\bf}{\it}
{\bf}{\it}
\theoremstyle{definition}
\newtheorem{definition}[theorem]{Definition}

\newtheorem{convention}[theorem]{Convention}

\newtheorem{fdefinition}[ftheorem]{Definition}{\bf}{\rm}

\theoremstyle{remark}

\newtheorem{remark}[theorem]{Remark}

\newtheorem{definition and corollary}[theorem]{Definition and Corollary}

\newtheorem{fexample}[ftheorem]{Example}{\it}{\rm}

\newcommand{\MID}{\! \! \mid}

\newcommand{\h}{\mathfrac{\h}}

\title{An exotic Deligne-Langlands correspondence for symplectic groups}
\author{Syu \textsc{Kato}
          \footnote{Graduate School of Mathematical Sciences, University of Tokyo, 3-8-1 Meguro Komaba 153-8914, Japan.} \footnote{Current address: Research Institute for Mathematical Sciences, Kyoto University, Oiwake Kita-Shirakawa Sakyo Kyoto 606-8502, Japan. \tt{E-mail:kato@kurims.kyoto-u.ac.jp}} \footnote{The author was partially supported by JSPS Research Fellowship for Young Scientists (PD) 15-10371 and JSPS Grant-in-Aid for Young Scientists (B) 20-740011 during this research.}}

\begin{document}
\maketitle

\begin{abstract}
Let $G = \mathop{Sp} ( 2n, \mathbb C )$ be a complex symplectic group. We introduce a $G \times ( \mathbb C ^{\times} ) ^{\ell + 1}$-variety $\mathfrak N _{\ell}$, which we call the $\ell$-exotic nilpotent cone. Then, we realize the Hecke algebra $\mathbb H$ of type $C _{n} ^{(1)}$ with three parameters via equivariant algebraic $K$-theory in terms of the geometry of $\mathfrak N _2$. This enables us to establish a Deligne-Langlands type classification of simple $\mathbb H$-modules under a mild assumption on parameters. As applications, we present a character formula and multiplicity formulas of $\mathbb H$-modules.
\end{abstract}


\begin{center}
{\sf Table of Contents}
\end{center}
\vskip -3mm
{\tt{\small

\hskip 15mm \S 1 \hskip 2mm Preparatory materials \\
\vskip -4.5mm \hskip 15mm \S 2 \hskip 2mm Hecke algebras and exotic nilpotent cones\\
\vskip -4.5mm \hskip 15mm \S 3 \hskip 2mm Clan decomposition\\
\vskip -4.5mm \hskip 15mm \S 4 \hskip 2mm On stabilizers of exotic nilpotent orbits\\
\vskip -4.5mm \hskip 15mm \S 5 \hskip 2mm Semisimple elements attached to $G \backslash \mathfrak N _1$\\
\vskip -4.5mm \hskip 15mm \S 6 \hskip 2mm A vanishing theorem\\
\vskip -4.5mm \hskip 15mm \S 7 \hskip 2mm Standard modules and an induction theorem\\
\vskip -4.5mm \hskip 15mm \S 8 \hskip 2mm Exotic Springer correspondence\\
\vskip -4.5mm \hskip 15mm \S 9 \hskip 2mm A deformation argument on parameters\\
\vskip -4.5mm \hskip 15mm \S 10 \hskip 0.4mm Main Theorems\\
\vskip -4.5mm \hskip 15mm \S 11 \hskip 0.4mm Consequences\\}
}

\section*{Introduction}
In their celebrated paper \cite{KL2}, Kazhdan and Lusztig gave a classification of simple modules of an affine Hecke algebra $\mathbb H$ with one-parameter in terms of the geometry of nilpotent cones. (It is also done by Ginzburg, c.f. \cite{CG}.) Since some of the affine Hecke algebras admit two or three parameters, it is natural to extend their result to multi-parameter cases. (It is called the unequal parameter case.) Lusztig realized the ``graded version" of $\mathbb H$ (with unequal parameters) via several geometric means \cite{L1, L2, L3} (c.f. \cite{L4}) and classified their representations in certain cases. Unfortunately, his geometries admit essentially only one parameter. As a result, his classification is restricted to the case where all of the parameters are certain integral power of a single parameter. It is enough for his main interest, the study of representations of $p$-adic groups (c.f. \cite{L5}). However, there are many areas of mathematics which wait for the full-representation theory of Hecke algebras with unequal parameters (see e.g. Macdonald's book \cite{M} and its featured review in MathSciNet).

In this paper, we give a realization of all simple modules of the Hecke algebra of type $C _n ^{(1)}$ with three parameters by introducing a variety which we call {\it the $\ell$-exotic nilpotent cone} (c.f. \S \ref{exgeom}). Our framework works for all parameters and realizes the whole Hecke algebra (Theorem \ref{iHecke}) and its specialization to each central character. Unfortunately, the study of our geometry becomes harder for some parameters and the result becomes less explicit in such cases. Even so, our result gives a definitive classification of simple modules of affine Hecke algebras of type $B _n ^{(1)}$ and $C _n ^{(1)}$ for almost all parameters including so-called real central character case. (See the argument after Theorem \ref{iDL}.)

Let $G$ be the complex symplectic group $\mathop{Sp} ( 2n, \mathbb C )$. We fix its Borel subgroup $B$ and a maximal torus $T \subset B$. Let $R$ be the root system of $( G, T )$. We embed $R$ into a $n$-dimensional Euclid space $\oplus _i \mathbb C \epsilon _i$ as $R = \{ \pm \epsilon _i \pm \epsilon _j \} \cup \{ \pm 2 \epsilon _i \}$. We define $V _1 := \mathbb C ^{2n}$ and $V _2 := ( \wedge ^2 V _1 ) / \mathbb C$. For each non-negative integer $\ell$, we put $\mathbb V _{\ell} := V _1 ^{\oplus \ell} \oplus V _2$ and call it {\it the $\ell$-exotic representation}. Let $\mathbb V ^+ _{\ell}$ be the positive part of $\mathbb V _{\ell}$ (for precise definition, see \S \ref{prem}). We define
$$F _{\ell} := G \times ^B \mathbb V ^+ _{\ell} \subset G \times ^B \mathbb V _{\ell} \cong G / B \times \mathbb V _{\ell}.$$
Composing with the second projection, we have a map
$$\mu _{\ell} : F _{\ell} \longrightarrow \mathbb V _{\ell}.$$
We denote the image of $\mu _{\ell}$ by $\mathfrak N _{\ell}$. This is the $G$-variety which we refer as {\it the $\ell$-exotic nilpotent cone}. We put $Z _{\ell} := F _{\ell} \times _{\mathfrak N _{\ell}} F _{\ell}$. Let $G _{\ell}:= G \times ( \mathbb C ^{\times} ) ^{\ell + 1}$. We have a natural $G _{\ell}$-action on $F _{\ell}$ (and $Z _{\ell}$). (In fact, the variety $F _{\ell}$ admits an action of $G \times \mathop{GL} ( \ell, \mathbb C) \times \mathbb C ^{\times}$. We use only a restricted action in this paper.)

Assume that $\mathbb H$ is the Hecke algebra with unequal parameters of type $C _n ^{( 1 )}$ (c.f. Definition \ref{Checke}). This algebra has three parameters $q _0, q _1, q _2$. All affine Hecke algebras of classical type with two parameters are obtained from $\mathbb H$ by suitable specializations of parameters (c.f. Remark \ref{remHecke}).

\begin{ftheorem}[= Theorem \ref{ehecke}]\label{iHecke}
We have an isomorphism
$$\mathbb H \stackrel{\cong}{\longrightarrow} \mathbb C \otimes _{\mathbb Z} K ^{G _2} ( Z _2 )$$
as algebras.
\end{ftheorem}

The Ginzburg theory suggests a classification of simple $\mathbb H$-modules by the $G$-conjugacy classes of the following Langlands parameters:

\begin{fdefinition}[Langlands parameters]
\item ${\bf 1)}$ A triple $\vec{q} := ( q _0, q _1, q _2 ) \in ( \mathbb C ^{\times} ) ^3$ is said to be admissible if $q _0 \neq q _1$, $q _2$ is not a root of unity of order $\le 2 n$, $q _0 q _1 ^{\pm 1} \neq q _2 ^m$ for $| m | < n$;
\item ${\bf 2)}$ A pair $( a, X ) = ( s, \vec{q}, X _0 \oplus X _1 \oplus X _2 ) \in G _2 \times \mathfrak N _2$ is called an admissible parameter iff $s$ is semisimple, $\vec{q}$ is admissible, and $s X _i = q _i X _i$ for $i = 0, 1, 2$.
\end{fdefinition}

For $a = ( s, \vec{q} ) \in G _{2}$, we put $G ( s ) := Z _G ( s )$ and $G _2 ( a ) := Z _{G _2} ( a )$.

Notice that our Langlands parameters do not have additional data as in the usual Deligne-Langlands-Lusztig correspondence. This is because the (equivariant) fundamental groups of orbits are always trivial (c.f. Theorem \ref{pi1 desc}). Instead, we have the following kind of difficulty:

\begin{fexample}[Non-regular parameters]\label{fcrit}
Let $G = \mathop{Sp} ( 4, \mathbb C )$ and let $a = ( \exp ( r \epsilon _1 + ( r + \pi \sqrt{- 1} ) \epsilon _2 ), e ^r, - e ^r, - e ^{2 r} ) \in T \times ( \mathbb C ^{\times} ) ^3$ ($r \in \mathbb C \backslash \pi \sqrt{- 1} \mathbb Q$). Then, the number of $G _2 ( a )$-orbits in $\mathfrak N _2 ^{a}$ is eight, while the number of corresponding representations of $\mathbb H$ is six. (c.f. Enomoto \cite{E}) These orbits contain weight vectors of $\epsilon _1 + \epsilon _2$ or ``$\epsilon _1$ \& $\epsilon _2$".
\end{fexample}

Now we state the main theorem of this paper:

\begin{ftheorem}[= Theorem \ref{EDL}]\label{iDL}
The set of $G$-conjugacy classes of admissible parameters is in one-to-one correspondence with the set of isomorphism classes of simple $\mathbb H$-modules if $q _2$ is not a root of unity of order $\le 2n$, and $q _0 q _1 ^{\pm 1} \neq q _2 ^{\pm m}$ holds for every $0 \le m < n$.
\end{ftheorem}

We treat a slightly more general case in Theorem \ref{DLmain} including Example \ref{fcrit}. Since the general condition is rather technical, we state only a part of it here.

By imposing an additional relation $q _0 + q _1 = 0$, the algebra $\mathbb H$ specializes to an extended Hecke algebra $\mathbb H _B$ of type $B _n ^{(1)}$ with two-parameters. (c.f. Remark \ref{remHecke}.) Therefore, Theorem \ref{iDL} also gives a definitive classification of simple $\mathbb H _B$-modules except for $- q _0 ^2 = q _2 ^{m}$ ($\left| m \right| < n$) or $q _2$ is a root of unity of order $\le 2 n$.

Let us illustrate an example which (partly) explains the title ``exotic":

\begin{fexample}[Equal parameter case]
Let $G = \mathop{Sp} ( 4, \mathbb C )$. Let $s = \exp ( r \epsilon _1 + r \epsilon _2 ) \in T$ ($r \in \mathbb C \backslash \pi \sqrt{- 1} \mathbb Q$). Fix $a _{0} = ( s, e ^{2r} ) \in G \times \mathbb C ^{\times}$ and $a = ( s, e ^r, - e ^r, e ^{2r} ) \in G _2$. Let $\mathcal N$ be the nilpotent cone of $G$. Then, the sets of $G ( s )$-orbits of $\mathcal N ^{a _{0}}$ and $\mathfrak N ^{a} _2$ are responsible for the usual and our exotic Deligne-Langlands correspondences. The number of $G ( s )$-orbits in $\mathcal N ^{a _{0}}$ is three. (Corresponding to root vectors of $\emptyset$, $2 \epsilon _1$, and ``$2 \epsilon _1$ \& $2 \epsilon _2$") The number of $G ( s )$-orbits in $\mathfrak N ^{a} _2$ is four. (Corresponding to weight vectors of $\emptyset$, $\epsilon _1$, $\epsilon _1 + \epsilon _2$, and ``$\epsilon _1$ \& $\epsilon _1 + \epsilon _2$") On the other hand, the actual number of simple modules arising in this way is four (c.f. Ram \cite{R} and \cite{E}).
\end{fexample}

The organization of this paper is as follows:

In \S 1, we fix notation and introduce exotic nilcones and related varieties. In particular, we present geometric structures involved in our varieties as much as we need in the later sections. In \S 2, we prove Theorem \ref{iHecke}, which connects our varieties with an affine Hecke algebra $\mathbb H$ of type $C _n ^{(1)}$. In order to simplify the study of representation theory of $\mathbb H$, we divide our varieties into a product of primitive ones in \S 3. In \S 4, we prove that the stabilizers of exotic nilpotent orbits are connected, which implies that ``the Lusztig part" of the Deligne-Langlands-Lusztig parameter should be always trivial in our situation. Unfortunately, we have no nice parabolic subgroup as Kazhdan-Lusztig employed in \cite{KL2}. We construct some explicit semisimple element out of each orbit in \S 5 for the sake of compensation. We introduce the notion of exotic Springer fibers and prove its odd-term vanishing result in \S 6, under the assumption that the parameters are sufficiently nice (including admissible case). Its proof essentially relies on the argument of \S 5. We define our standard modules as the total homology group of exotic Springer fibers in \S 7. At the same time, we present an induction theorem, which claims that they behave well under inductions. In \S 8, we present an analogue of the Springer correspondence for exotic nilcones. In order to prove Theorem \ref{iDL}, we still need two additional structural results. One is that our geometric structure is preserved by replacing the central character by a suitable real positive one. The other is that we can embed the corresponding finite Weyl group into the graded version of $\mathbb H$. Our proofs of both results essentially use admissibility of parameters. These results occupy \S 9. With the knowledge of all of the previous sections except for \S7, we prove Theorem \ref{iDL} in \S 10. The last section \S 11 concerns with applications, which are straight-forward consequences of Ginzburg theory assuming the results presented in earlier sections.

{\small
{\bf Acknowledgment:} The present from of this paper\footnote{{\bf Note:} After the original version of this paper is circulated (in 2006, with different argument and weaker conclusion in Theorem \ref{iDL}, and consequently give a classification of $\mathbb H$-modules only with a help of Lusztig's results \cite{L1,L2,L3}), there appeared two kinds of related works. One is the study of geometry which is connected to our nilcone by Achar-Henderson \cite{AH}, Enomoto \cite{E2}, Finkelberg-Ginzburg-Travkin \cite{FGT}, Springer \cite{Sp}, Travkin \cite{Tr}, and the other is the classification of tempered dual by Opdam and Solleveld \cite{OS, OS2, So}. For the former, I have included explanations about the situation as much as I could in order to avoid potential problems. For the latter, we are preparing another paper \cite{CK} in this direction.} is heavily benefited from the comments from Pramod Achar, Susumu Ariki, Michel Brion, Masaki Kashiwara, Michael Finkelberg, Anthony Henderson, George Lusztig, Hiraku Nakajima, Eric Opdam, Midori Shiota, Toshiyuki Tanisaki, Masahiko Yoshinaga, and discussion with Naoya Enomoto, Hisayosi Matumoto, Eric Vasserot, Tonny A. Springer. The author wants to express his deep gratitude to all of them for their kindness, warmness, and tolerance. In particular, Professor Ariki kindly arranged him an opportunity to talk at a seminar at RIMS. The author wishes to express his gratitude to him and all the participants of the seminar. Last, but not least, the author would like to thank the referee for his kindness and patience.}

\section{Preparatory materials}\label{prem}
Let $G := \mathop{Sp} ( 2n, \mathbb C )$. Let $B$ be a Borel subgroup of $G$. Let $T$ be a maximal torus of $B$. Let $X ^* ( T )$ be the character group of $T$. Let $R$ be the root system of $(G, T)$ and let $R ^+$ be its positive part defined by $B$. We embed $R$ and $R ^+$ into a $n$-dimensional Euclid space $\mathbb E = \oplus _i \mathbb C \epsilon _i$ with standard inner product as:
$$R ^+ = \{ \epsilon _i \pm \epsilon _j \} _{i < j} \cup \{ 2 \epsilon _i \} \subset \{ \pm \epsilon _i \pm \epsilon _j \} \cup \{ \pm 2 \epsilon _i \} = R \subset \mathbb E.$$
By the inner product, we identify $\epsilon _i$ with its dual basis. We put $\epsilon _{i} := - \epsilon _{- i}$ when $-n \le i < 0$. We put $\alpha _i := \epsilon _i - \epsilon _{i + 1}$ ($i = 1, \ldots, n - 1$) or $2 \epsilon _n$ ($i = n$). Let $W$ be the Weyl group of $( G, T )$. For each $\alpha _i$, we denote the reflection of $\mathbb E$ corresponding to $\alpha _i$ by $s _i$. Let $\ell : W \rightarrow \mathbb Z _{\ge 0}$ be the length function with respect to $( B, T )$. We denote by $\dot{w} \in N _G ( T )$ a lift of $w \in W$. For a subgroup $H \subset G$ containing $T$, we put ${} ^w H := \dot{w} H \dot{w} ^{-1}$. For a group $H$ and its element $s$, we put $H ( s ) := Z _H ( s )$. For a subset $S \subset H$, we put $H ( S ) := \cap _{s \in S} H ( s )$. We denote the identity component of $H$ by $H ^{\circ}$. We denote by $R ( H )$ and $R ( H ) _s$ the representation ring of $H$ and its localization along the evaluation at $s \in H$, respectively. For each $\alpha \in R$, we denote the corresponding one-parameter unipotent subgroup of $G$ (with respect to $T$) by $U _{\alpha}$. We define $\mathfrak g, \mathfrak t, \mathfrak g ( s ),$ etc$\ldots$ to be the Lie algebras of $G, T, G ( s ),$ etc$\ldots$, respectively.

For a $T$-module $V$, we define its weight $\lambda$-part (with respect to $T$) as $V [\lambda]$. We define the positive part $V ^+$ and negative part $V ^-$ of $V$ as
$$V ^+ := \bigoplus _{\lambda \in \mathbb Q _{\ge 0} R ^+ - \{ 0 \} } V [ \lambda ], \text{ and } V ^- := \bigoplus _{\lambda \in \mathbb Q _{\le 0} R ^+ - \{ 0 \} } V [ \lambda ],$$
respectively.
We denote the set of $T$-weights of $V$ by $\Psi ( V )$.

In this paper, a segment is a set of integers $I$ written as $I = [ i _1, i _2 ] \cap \mathbb Z$ for some integers $i _1 \le i _2$. By abuse of notation, we may denote $I$ by $[ i _1, i _2 ]$. For a segment $I$, we set $I ^* := I$ (if $0 \not\in I$) or $I - \{ 0 \}$ (if $0 \in I$). We denote the absolute value function by $| \bullet | : \mathbb C \rightarrow \mathbb R _{\ge 0}$. We set $\Gamma _0 := 2 \pi \sqrt{-1} \mathbb Z \subset \mathbb C$ and set $\exp : \mathbb E \to T$ to be the exponential map. We normalize the map $\exp$ so that $\ker \exp \cong \sum _{i = 1} ^n \Gamma _0 \epsilon _i$.

A variety in this paper is a quasi-projective reduced scheme of finite type over $\mathbb C$. Its points are closed points. If an algebraic group $H$ acts on a variety $\mathcal X$, then we denote the stabilizer of the $H$-action at $x \in \mathcal X$ by $\mathsf{Stab} _H x$. For each $h \in H$, we denote by $\mathcal X ^h$ the $h$-fixed point set of $\mathcal X$. For a variety $\mathcal X$, we denote by $H _{\bullet} ( \mathcal X )$ the Borel-Moore homology groups with coefficients $\mathbb C$.

\subsection{Exotic nilpotent cones}\label{exgeom}
Let $\ell = 0, 1,$ or $2$. We define $V _1 := \mathbb C ^{2n}$ (vector representation) and $V _2 := ( \wedge ^2 V _1 ) / \mathbb C$. These representations have $B$-highest weights $\epsilon _1$ and $\epsilon _1 + \epsilon _2$, respectively. We put $\mathbb V _{\ell} := V _1 ^{\oplus \ell} \oplus V _2$ and call it {\it the $\ell$-exotic representation of $\mathop{Sp} ( 2 n )$}. For $\ell \ge 1$, the set of non-zero weights of $\mathbb V _{\ell}$ is in one-to-one correspondence with $R$ as
\begin{eqnarray}
R \ni \begin{cases} \pm 2 \epsilon _i \leftrightarrow \pm \epsilon _i & \in \Psi ( V _1 )\\
\pm \epsilon _i \pm \epsilon _j \leftrightarrow \pm \epsilon _i \pm \epsilon _j & \in \Psi ( V _2 )
\end{cases}. \label{one-to-one}
\end{eqnarray}
We define
$$F _{\ell} := G \times ^B \mathbb V ^+ _{\ell} \subset G \times ^B \mathbb V _{\ell} \cong G / B \times \mathbb V _{\ell}.$$
Composing with the second projection, we have a map
$$\mu _{\ell} : F _{\ell} \longrightarrow \mathbb V _{\ell}.$$
We denote the image of $\mu _{\ell}$ by $\mathfrak N _{\ell}$. We call this variety {\it the $\ell$-exotic nilpotent cone}. By abuse of notation, we may denote the map $F _{\ell} \rightarrow \mathfrak N _{\ell}$ also by $\mu _{\ell}$.

\begin{convention}
For the sake of simplicity, we define objects $F$, $\mathfrak N$, $\mathbb V$, $\mu$, etc... to be the objects $F _{\ell}$, $\mathfrak N _{\ell}$, $\mathbb V _{\ell}$, $\mu _{\ell}$ etc... with $\ell = 1$.
\end{convention}

We summarize some basic geometric properties of $\mathfrak N _{\ell}$:

\begin{theorem}[Geometric properties of $\mathfrak N _{\ell}$]\label{fgeom}
We have the following:
\begin{enumerate}
\item The defining ideal of $\mathfrak N _{\ell}$ is generated by $G$-invariant polynomials of $\mathbb C [ \mathbb V _{\ell} ]$ without constant terms;
\item The variety $\mathfrak N _{\ell}$ is normal;
\item For $\ell = 1, 2$, the map $\mu _{\ell}$ is a birational projective morphism onto $\mathfrak N _{\ell}$;
\item Every fiber of the map $\mu _{\ell}$ is connected;
\begin{center}
In the below, we present properties which are valid only for the $\ell = 1$ case.
\end{center}
\item The set of $G$-orbits in $\mathfrak N _1$ is finite;
\item The map $\mu _1$ is stratified semi-small with respect to the stratification of $\mathfrak N _1$ given by $G$-orbits.
\end{enumerate}
\end{theorem}

\begin{proof}
The proof is given after Lemma \ref{irrc} since we need extra notation.
\end{proof}

\begin{lemma}\label{lci}
We have a natural identification
$$F _{\ell} \cong \{ ( g B, X ) \in G / B \times \mathbb V _{\ell} ; X \in g \mathbb V ^+ _{\ell} \}.$$
\end{lemma}
\begin{proof}
Straightforward.
\end{proof}

Let $G _{\ell} := G \times ( \mathbb C ^{\times} ) ^{\ell + 1}$. We define a $G _{\ell}$-action on $\mathfrak N _{\ell}$ as
$$G _{\ell} \times \mathfrak N _{\ell} \ni ( g, q _{2 - \ell}, \ldots, q _2 ) \times ( X _{2 - \ell} \oplus \cdots \oplus X _2 ) \mapsto ( q _{2 - \ell} ^{- 1} g X _{2 - \ell} \oplus \cdots \oplus q _2 ^{- 1} g X _2 ) \in \mathfrak N _{\ell}.$$
(Here we always regard $X _{2 - \ell}, \ldots, X _1 \in V _1$ and $X _2 \in V _2$.) Similarly, we have a natural $G _{\ell}$-action on $F _{\ell}$ which makes $\mu _{\ell}$ a $G _{\ell}$-equivariant map. We define $Z _{\ell} := F _{\ell} \times _{\mathfrak N _{\ell}} F _{\ell}$. By Lemma \ref{lci}, we have
$$Z _{\ell} := \{ ( g _1 B, g _2 B, X ) \in ( G / B ) ^2 \times \mathbb V _{\ell} ; X \in g _1 \mathbb V _{\ell} ^+ \cap g _2 \mathbb V _{\ell} ^+ \}.$$
We put
$$Z _{\ell} ^{123} := \{ ( g _1 B, g _2 B, g _3 B, X ) \in ( G / B ) ^3 \times \mathbb V _{\ell} ; X \in g _1 \mathbb V _{\ell} ^+ \cap g _2 \mathbb V _{\ell} ^+  \cap g _3 \mathbb V _{\ell} ^+ \}.$$
We define $p _i : Z _{\ell} \ni ( g _1 B, g _2 B, X ) \mapsto ( g _i B, X ) \in F _{\ell}$ and $p _{ij} : Z ^{123} _{\ell} \ni ( g _1 B, g _2 B, g _3 B, X ) \mapsto ( g _i B, g _j B, X ) \in Z _{\ell}$ ($i, j \in \{ 1, 2, 3 \}$). We also put $\tilde{p} _i : F _{\ell} \times F _{\ell} \rightarrow F _{\ell}$ as the first and second projections ($i = 1, 2$). (Notice that the meaning of $p _{i}, \tilde{p} _{i}, p _{ij}$ depends on $\ell$. The author hopes that there occurs no confusion on it.)

\begin{lemma}
The maps $p _i$ and $p _{ij}$ $(1 \le i < j \le 3)$ are projective.
\end{lemma}

\begin{proof}
The fibers of the above maps are given as the subsets of $G / B$ defined by incidence relations. It is automatically closed, and we obtain the result.
\end{proof}

We have a projection
$$\pi _{\ell} : Z _{\ell} \ni ( g _1 B, g _2 B, X ) \mapsto ( g _1 B, g _2 B ) \in G / B \times G / B.$$

For each $w \in W$, we define a point ${\mathsf p} _w := B \times \dot{w} B \in G / B \times G / B$. This point is independent of the choice of $\dot{w}$. We put $\mathsf O _w := G {\mathsf p} _w \subset G / B \times G / B$. By the Bruhat decomposition, we have
\begin{eqnarray}
G / B \times G / B = \bigsqcup _{w \in W} \mathsf O _w. \label{BD}
\end{eqnarray}

\begin{lemma}\label{irrc}
The variety $Z _{\ell}$ $(\ell = 1,2)$ consists of $\left| W \right|$-irreducible components. Moreover, the dimensions of all of the irreducible components of $Z$ are equal to $\dim F$.
\end{lemma}

\begin{proof}
We first prove the assertion for $Z = Z _1$. By $(\ref{BD})$, the structure of $Z$ is determined by the fibers over ${\mathsf p} _w$. We have
$$\pi ^{-1} ( {\mathsf p} _w ) = \mathbb V ^+ \cap \dot{w} \mathbb V ^+.$$
By the dimension counting using (\ref{one-to-one}), we deduce
\begin{align*}
\dim \mathbb V ^+ \cap \dot{w} \mathbb V ^+ = & \dim V ^+ _1 \cap \dot{w} V ^+ _1 + \dim V ^+ _2 \cap \dot{w} V ^+ _2 \\ = & \# ( R _l ^+ \cap w R _l ^+ ) + \# ( R _s ^+ \cap w R _s ^+ ) = N - \ell ( w ),
\end{align*}
where $N := \dim \mathbb V ^+ = \dim G / B$ and $R _l ^+, R ^+ _s$ are the sets of long and short positive roots, respectively. As a consequence, we deduce
$$\dim \pi ^{- 1} ( \mathsf O _w ) = N + \ell ( w ) + N - \ell ( w ) = 2 N.$$
Thus, each $\overline{\pi ^{- 1} ( \mathsf O _w )}$ is an irreducible component of $Z$. Moreover, we have $\pi ^{- 1} ( \mathsf O _1 ) \cong F$, which implies that the dimensions of irreducible components of $Z$ are equal to $\dim F$.\\
Next, we prove the assertion for $Z _2$. By forgetting the first $V _1$-factor, we have a surjective map $\eta: Z _2 \rightarrow Z$. We have a surjective map $\eta ^{\prime} : Z \rightarrow Z _0$ given by forgetting the $V _1$-factor. The fiber of $( \eta ^{\prime} \circ \eta )$ at $x \in Z$ is isomorphic to the two-fold product of the fiber of $\eta ^{\prime}$ at $\eta ^{\prime} ( x )$. The latter fiber is isomorphic to the vector space $V _1 ^+ \cap g V _1 ^+$ when $\pi ( x ) = ( 1, g ) {\mathsf p} _1$. Therefore, the preimage of each irreducible component of $Z$ gives an irreducible component of $Z _2$. These irreducible components are distinct since their images under $\eta$ must be distinct. Hence, the number of irreducible components of $Z _2$ is equal to the number of irreducible components of $Z$ as desired.
\end{proof}

\begin{proof}[Proof of Theorem \ref{fgeom}]
The weight distribution of $\mathbb V^+$ and the Hesselink theory (c.f. \cite{P} Theorem 1) claims that $\mu _{\ell}$ gives a birational projective morphism onto an irreducible component of the Hilbert nilcone of $\mathbb V _{\ell}$. Here the Hilbert nilcone of $\mathbb V _{\ell}$ is an irreducible normal variety by Dadok-Kac \cite{DK} or Schwarz \cite{Sc}. In particular, our variety $\mathfrak N _{\ell} \subset \mathbb V _{\ell}$ is the Hilbert nilcone itself. Therefore, we obtain 1--3). 4) is an immediate consequence of 2), 3), and the Zariski main theorem (c.f. \cite{CG} 3.3.26). 5) is proved as a part of Proposition \ref{wN1}. We show 6). Let $\widehat{\mathbb O}$ be the inverse image of a $G$-orbit $G. X = \mathbb O \subset \mathfrak N$ under the map $\mu \circ p _2$. Then, we have
$$\dim \mathbb O + 2 \dim \mu ^{-1} ( X ) \le \dim \widehat{\mathbb O}.$$
The dimension of the RHS is less than or equal to $\dim F$, which is the (constant) dimension of irreducible components of $Z$. In particular, we have
$$\dim \mathbb O + 2 \dim \mu ^{-1} ( X ) \le \dim \mathfrak N = \dim F,$$
which implies that $\mu$ is semi-small.
\end{proof}

By a general result of \cite{Gi} p135 (c.f. \cite{CG} 2.7), the $G _{\ell}$-equivariant $K$-group of $Z _{\ell}$ becomes an associative algebra via the map
$$\star : K ^{G _{\ell}} ( Z _{\ell} ) \times K ^{G _{\ell}} ( Z _{\ell} ) \ni ([ \mathcal E], [\mathcal F]) \mapsto \sum _{i \ge 0} ( - 1 ) ^i [ \mathbb R ^{i} ( p _{13} ) _* ( p _{12} ^* \mathcal E \otimes ^{\mathbb L} p _{23} ^* \mathcal F )]  \in K ^{G _{\ell}} ( Z _{\ell} ).$$
Moreover, the $G _{\ell}$-equivariant $K$-group of $F _{\ell}$ becomes a representation of $K ^{G _{\ell}} ( Z _{\ell} )$ as
$$\circ : K ^{G _{\ell}} ( Z _{\ell} ) \times K ^{G _{\ell}} ( F _{\ell} ) \ni ([ \mathcal E], [\mathcal K]) \mapsto \sum _{i \ge 0} ( - 1 ) ^i [ \mathbb R ^{i} ( p _{1} ) _* ( \mathcal E \otimes ^{\mathbb L} \tilde{p} _{2} ^* \mathcal K )] \in K ^{G _{\ell}} ( F _{\ell} ).$$
Here we regard $\mathcal E$ as a sheaf over $F _{\ell} \times F _{\ell}$ via the natural embedding $Z _{\ell} \subset F _{\ell} \times F _{\ell}$.

\subsection{Definition of parameters}\label{sJNF}
In this subsection, we present the definitions of parameters which we need in the sequel. First, we put $a _0 := ( 1, 1, -1, 1 ) \in G _2$. (The value $a _0$ is special in the sense it naturally gives the Weyl group of type $C$ in our framework. C.f. \S \ref{esc})

\begin{definition}[Configuration of semisimple elements]
\item {\bf 1)} An element $a = ( s, q _0, q _1, q _2 ) \in G _2$ is called pre-admissible iff $s$ is semisimple, $q _0 \neq q _1$, $q _2$ is not a root of unity of order $\le 2 n$.
\item {\bf 2)} An element $a \in G _2$ is called finite if $\mathfrak N _2 ^{a}$ has only finitely many $G _2 ( a )$-orbit.
\item {\bf 3)} A pre-admissible element $a = ( s, q _0, q _1, q _2 )$ is called admissible if $q _0 q _1 ^{\pm 1} \neq q _2 ^{\pm m}$ holds for every $0 \le m < n$.
\end{definition}

For a pre-admissible element $a = ( s, q _0, q _1, q _2 )$, we put
$$\mathbb V ^{a} _2 = V _1 ^{(s, q _0)} \oplus V _1 ^{(s, q _1)} \oplus V _2 ^{(s, q _2)} \subset V _1 \oplus V _1 \oplus V _2 = \mathbb V _2.$$

In the below, we may denote $( q _0, q _1, q _2 ) \in ( \mathbb C ^{\times} ) ^3$ by $\vec{q}$ for the sake of simplicity.

Let $a = ( s, \vec{q} ) \in G _2$ be a pre-admissible element such that $s \in T$. We sometimes denote it as
$$s = \exp \left( \sum _{i = 1} ^n \log _i ( s ) \epsilon _i \right) \in \exp ( \mathbb E ) \cong T,$$
where $\log _i ( s ) \in \mathbb C$. 

\begin{remark}
The values of $\log _i ( s )$ are determined modulo $\Gamma _0$. Here we understand that $\log _i ( s )$ is a fixed choice of a representative in $\log _i ( s ) + \Gamma _0$.
\end{remark}

\begin{definition}[Admissible parameters]
\item {\bf 1)} A pre-admissible parameter is a pair
$$\nu = ( a, X ) = ( s, \vec{q}, X _1 \oplus X _2 ) \in G _2 \times \mathfrak N _1$$
such that $a$ is pre-admissible, $( s - q _0 ) ( s - q _1 ) X _1 = 0$, and $s X _2 = q _2 X _2$;
\item For a pre-admissible $a \in G _2$, we denote by $\Lambda _{a}$ the set of $G ( s )$-conjugacy classes of pre-admissible parameters of the form $( a, Y )$, where $Y \in \mathbb V$.
\item {\bf 2)} A pre-admissible parameter $\nu = ( a, X )$ is called admissible if $a$ is admissible.
\end{definition}

\subsection{Orbit structures arising from $\mathfrak N _{\ell}$}\label{os}

In the below, we fix vectors in $V _1$ and $V _2$ as follows:
\begin{itemize}
\item For each $i \in [ -n, n ] ^*$, we define $0 \neq x _i \in V _1$ as a non-zero vector of weight $\epsilon _i$;
\item For each distinct $i, j \in [ - n, n ] ^*$, we define $y _{ij} \in V _2$ to be a non-zero vector of weight $\epsilon _{i} - \epsilon _{j}$.
\end{itemize}

The following is a slight enhancement of the good basis of Ohta \cite{Oh} (1.3).

\begin{definition}[Signed partitions]
Let $\mathbf J := \{ J _1, J _2, \ldots \}$ be a collection of sequence of elements of $[ - n, n ] ^*$. (I.e. each $J _k \in \mathbf J$ is a sequence $( J _k ^1, J _k ^2, \ldots )$ in $[-n,n] ^*$.) We put $J ^{+} _k = ( |J _k ^1|, |J _k ^2 |, \ldots )$ for each $k = 1,2,\ldots$. We call $\mathbf J$ a signed partition of $n$ if and only if $\{ J _1 ^+, J _2 ^+, \ldots \}$ gives a subdivision of $[1,n]$. I.e. we have
$$[ 1, n ] = \bigcup _{k \ge 1} J _k ^+ = \bigcup _{k \ge 1} \{ |j| ; j \in J _k \} \text{ and } J _k ^+ \cap J _{k ^{\prime}} ^+ = \emptyset \text{ for } k \neq k ^{\prime}.$$
For each member $J$ of a signed partition $\mathbf J$, we define a subtorus
$$T _J := \exp \sum _{i \in J} \mathbb C \epsilon _{i} \subset T.$$
Let $\lambda := ( \lambda _1 \ge \lambda _2 \ge \cdots )$ be a partition of $n$. Then, we regard it as a signed partition by setting
$$J ^i _j := i + \sum _{k = 1} ^{j - 1} \lambda _k \text{ if } \lambda _j \neq 0 \text{ and } 1 \le i \le \lambda _{j}.$$
\end{definition}

\begin{definition}[Foot functions]\label{ff}
Let $\ell = 0, 1,$ or $2$. A collection of $\ell$-tuple of functions $\delta _k : [ -n , n ] ^* \rightarrow \{ 0, 1 \}$ for $1 \le k \le \ell$ is called a $\ell$-foot function of $n$. We denote a $\ell$-foot function $\{ \delta _k \} _{k = 1} ^{\ell}$ by $\vec{\delta}$.
\end{definition}

Notice that Definition \ref{ff} claims that $\vec{\delta} = \emptyset$ when $\ell = 0$.

\begin{definition}[Marked partitions, blocks, and normal forms]
Let $\ell$ be as in Definition \ref{ff}. We refer a pair $\sigma = ( \mathbf J, \vec{\delta} )$ consisting of a signed partition and a $\ell$-foot function of $n$ as a $\ell$-marked partition if the following condition holds:
\begin{itemize}
\item For each $J \in \mathbf J$ and $m = 1, \ldots, \ell$, we have
$$\# \{ j \in J ; \delta _m ( j ) = 1 \} + \# \{ j \in J ; \delta _m ( - j ) = 1 \} \le 1.$$
\end{itemize}
For each $J \in \mathbf J$, we define the $\ell$-block $\mathbf v _{\sigma} ^J = \mathbf v _{\sigma, 1} ^J + \mathbf v _{\sigma, 2} ^J \in \mathbb V$ associated to $( \sigma, J ) = ( \sigma, \{ J ^1, J ^2, \ldots \})$ as:
\begin{align*}
\mathbf v _{\sigma, 1} ^J := & \sum _{j \in J} \sum _{k = 1} ^{\ell} ( \delta _k ( j ) x _{j} + \delta _k ( - j ) x _{ - j} ) \in V _1\\
\mathbf v _{\sigma, 2} ^J := & \sum _{j \ge 1} y _{J ^{j}, J ^{j + 1}} \in V _2,
\end{align*}
where we regard $y _{J _k ^{j}, J _k ^{j + 1}} \equiv 0$ whenever $J _k$ nor $J _k ^{j + 1}$ is non-existent.\\
A $\ell$-normal form $\mathbf v _{\sigma} = \mathbf v _{\sigma, 1} + \mathbf v _{\sigma, 2} \in \mathbb V$ associated to $\sigma$ is defined as:
\begin{eqnarray*}
\mathbf v _{\sigma, 1} := \sum _{J \in \mathbf J} \mathbf v _{\sigma, 1} ^J \in V _1, \text{ and }
\mathbf v _{\sigma, 2} := \sum _{J \in \mathbf J} \mathbf v _{\sigma, 2} ^J \in V _2.
\end{eqnarray*}
\end{definition}

\begin{remark}
We regard that $\ell$-normal forms are elements of $\mathbb V = \mathbb V _1$, regardless the value of $\ell$. 
\end{remark}

\begin{definition}[Strict normal forms]\label{snf}
A $\ell$-marked partition $\sigma = ( \mathbf J, \vec{\delta} )$ is called strict if and only if the following four conditions hold:
\begin{enumerate}
\item $\mathbf J$ is obtained from a partition $\lambda$ of $n$;
\item We have $\delta _2 \equiv 0$ and $\delta _1 ( j ) = 0$ for every $j \in [ -n,  -1]$;
\end{enumerate}
Before stating the rest of the conditions, we introduce extra notation. Assume the above two conditions. If we have $\delta _1 ( j ) = 1$ for $j \in J$, then we set $\underline{\#} J := \# \{ j ^{\prime} \in J ; j ^{\prime} \le j \}$ and $\overline{\#} J := \# \{ j ^{\prime} \in J ; j ^{\prime} > j \}$.
\begin{itemize}
\item[3.] Let $k < m$ be two integers and let $\mathbf J = \{ J _1, J _2, \ldots \}$. Then, we have $\delta _1 \MID _{J _m} \equiv 0$ if $\# J _k = \# J _m$;
\item[4.] Let $J, J ^{\prime} \in \mathbf J$ be a pair such that $\delta _1 ( j ) = 1 = \delta _1 ( j ^{\prime} )$ for some $j \in J$ and $j ^{\prime} \in J ^{\prime}$. If $\# J > \# J ^{\prime}$, then we have
$$\underline{\#} J > \underline{\#} J ^{\prime} \text{ and } \overline{\#} J > \overline{\#} J ^{\prime}.$$
\end{itemize}
Conditions are not applicable when $\delta _2$ or $\delta _1$ are non-existent. Notice that only the first condition survives when $\ell = 0$. A normal form attached to a strict $\ell$-marked partition is called a $\ell$-strict normal form.
\end{definition}

In the below, we refer foot functions, blocks, normal forms..., to be the $1$-foot functions, $1$-blocks, $1$-normal forms..., respectively. Moreover, we naturally identify strict $1$-normal forms and strict $2$-normal forms since $\delta _2 \equiv 0$ for $2$-strict marked partitions.

Let $\mathsf{Irrep} W$ be the set of isomorphism classes of irreducible $W$-modules.

\begin{theorem}[Orbit description of $\mathfrak N _1$]\label{N1}
We have:
\begin{enumerate}
\item The set of strict $1$-normal forms is in one-to-one correspondence with the the set of $G$-orbits of $\mathfrak N _1$;
\item We have $\# ( G \backslash \mathfrak N _1 ) = \# \mathsf{Irrep} W$;
\item For each $X \in \mathfrak N _1$, the group $\mathsf{Stab} _G X$ is connected.
\end{enumerate}
\end{theorem}

\begin{remark}
The original form of the proof of Theorem \ref{N1} (in \cite{K}) employs explicit calculation using basis. In the meantime, Springer \cite{Sp} gives a base-free proof (with stronger consequences). The proof given here is somewhat the mixture of the both, which the author gives it for the sake of completeness. Note that the closure relation of the orbits of $\mathfrak N _1$ is calculated by Achar-Henderson \cite{AH}.
\end{remark}

The proof of Theorem \ref{N1} is obtained as a combination of Proposition \ref{2con} and Theorem \ref{eS} by using the knowledge of the following:

\begin{proposition}[Weak version of Theorem \ref{N1}]\label{wN1}
We have:
\begin{enumerate}
\item Each $G$-orbit of $\mathfrak N _1$ contains a strict normal form;
\item The number of elements of the set of strict marked partitions is less than or equal to $\# \mathsf{Irrep} W$.
\end{enumerate}
\end{proposition}

\begin{proof}
By a result of Ohta-Sekiguchi \cite{Se, Oh}, the set of strict $0$-marked partitions are in one-to-one correspondence with the set of $G$-orbits of $\mathfrak N _0$ via the assignment $\sigma \mapsto G \mathbf v _{\sigma}$. We have
$$\mathbb C [ \mathbb V ] ^G \cap \mathbb C [ \mathbb V _{0} ] = \mathbb C [ \mathbb V _{0} ] ^G,$$
which gives the natural projection map
$$\mathfrak N _{1} \longrightarrow \mathfrak N _{0}$$
obtained from the natural projection $\mathbb V _{1} \to \mathbb V _{0}$. (In fact we have $\mathbb C [ \mathbb V ] ^G = \mathbb C [ \mathbb V _{0} ] ^G$. But this fact is not used here.) It follows that each orbit of $\mathfrak N _1$ contains a vector of type $v = v _1 \oplus \mathbf v _{2, \lambda}$, where $\lambda$ is a partition of $n$ regarded as a strict $0$-marked partition $( \mathbf J, \emptyset )$ in a natural way.\\
Consider the action of
$$G ^{\prime} := \mathop{Sp} ( 2 \lambda _1 ) \times \mathop{Sp} ( 2 \lambda _2 ) \times \cdots \subset \mathop{Sp} ( 2 n ),$$
which are embedded so that $T \subset G ^{\prime}$ and $V _1$ restricted to $G ^{\prime}$ has the form
$$\mathsf{Res} ^G _{G ^{\prime}} V _1 = \bigoplus _{k \ge 1} V _1 ^{k}$$
such that $V _1 ^{k}$ is a vector representation of $\mathop{Sp} ( 2 \lambda _k )$ with $T$-weights $\pm \epsilon _i$ for $i = 1 + \sum _{j = 1} ^{k - 1} \lambda _j, \ldots, \lambda _k + \sum _{j = 1} ^{k-1} \lambda _j$.\\
Let $\omega$ be the symplectic form on $V _1$ which is preserved by $G$. For each $k$, we put $\omega _k := \omega \MID _{V _1 ^{k}}$. We have $v = \sum _{k \ge 1} v _k$, where $v _k = v _{1, k} \oplus \mathbf v _{2, \lambda} ^{J _k} \in V _1 ^k \oplus \wedge ^2 V _1 ^k$. We consider an identification of $y _{ij}$ ($i,j \in J_k$) with a matrix such that $y _{ij} x _k = c _{ij} x _i$ ($k = j$), $- c _{ij} x _{j}$ ($k = - i$), $0$ (otherwise), for some $c _{ij} \in \mathbb C$. We arrange $\{ c _{ij} \} _{i,j}$ so that $\wedge ^2 V _1 ^{k}$ is $\mathop{Sp} ( 2 \lambda _k )$-equivariantly identified with the subset of $\mathrm{End} V _1 ^{k}$ such that
$$\omega _k ( y _{ij} x, x^{\prime} ) = \omega _k ( x, y _{ij} x^{\prime} ) \text{ for each } x, x ^{\prime} \in V _1 ^k \text{ and } i, j \in ( J _k \cup - J _k ).$$
The Ohta-Sekiguchi result asserts that this gives an identification of $\mathop{Sp} ( 2 \lambda _k ) \mathbf v _{2, \lambda} ^{J _k}$ and the set of linear nilpotent endomorphisms on $V _1 ^{k}$ of maximal rank ($= \dim V _1 ^k - 2$) which preserve $\omega _k$. Since $v _{1, k}$ can be complemented to a suitable choice of a standard basis of $V _1 ^k$ (as a symplectic vector space), we deduce that a suitable change of symplectic basis makes $v _{1, k}$ into one of $x _i$ ($i >0$). This implies that $v$ can be transformed into a $1$-normal form associated to $( \mathbf J, \delta _1 ) = ( \lambda, \delta _1 )$ which satisfies \ref{snf} 1) and 2).\\
Now for each $k < k ^{\prime}$, we examine the $\mathop{Sp} ( 2 \lambda _k + 2 \lambda _{k ^{\prime}} )$-orbit which contains $v _{k} + v _{k ^{\prime}} \in V _1 ^{k} \oplus V _1 ^{k ^\prime}$. We have $\lambda _k \ge \lambda _{k ^{\prime}}$ by \ref{snf} 1). We put $\xi := \mathbf v _{2, \sigma} ^{J _k} + \mathbf v _{2, \sigma} ^{J _{k ^{\prime}}}$. The $\mathop{Sp} ( 2 \lambda _k + 2 \lambda _{k ^{\prime}} )$-conjugacy class of $\xi$ is the set of nilpotent endomorphisms of $V _1 ^{k} \oplus V _1 ^{k ^\prime}$ which preserve $\omega$ and have $( \lambda _k, \lambda _k, \lambda _{k ^{\prime}}, \lambda _{k ^{\prime}} )$ as its (nilpotent) Jordan form. If $v _{1, k} = 0$ or $v _{1, k ^{\prime}} =0$ hold, then \ref{snf} 3) and 4) are satisfied for the pair $( J _k, J _{k ^{\prime}} )$. Hence, we assume $v _{1, k} \neq 0 \neq v _{1, k ^{\prime}}$ in the below. We have
\begin{align*}
& \xi ^{\underline{\#} J _k} v _{1, k} = 0, \xi ^{\underline{\#} J _k - 1} v _{1,k} \neq 0, \text{ and }\xi ^{\underline{\#} J _{k ^{\prime}}} v _{1, k ^{\prime}} = 0, \xi ^{\underline{\#} J _{k ^{\prime}} - 1} v _{1,k ^{\prime}} \neq 0\\
& v _{1, k} \in \mathrm{Im} \xi ^{\overline{\#} J _k}, v _{1,k} \not\in \mathrm{Im} \xi ^{\overline{\#} J _k + 1}, \text{ and } v _{1, k ^{\prime}} \in \mathrm{Im} \xi ^{\overline{\#} J _{k ^{\prime}}}, v _{1, k ^{\prime}} \not\in \mathrm{Im} \xi ^{\overline{\#} J _{k ^{\prime}} + 1}.
\end{align*}
If $\underline{\#} J _k \le \underline{\#} J _{k ^{\prime}}$ or $\overline{\#} J _k \le \overline{\#} J _{k ^{\prime}}$ holds, then we can regard $v _{1, k} + v _{1, k ^{\prime}}$ as a part of a standard basis of a $( \mathbf v _{2,\sigma} )$-stable symplectic subspace of $( V _1 ^{k ^{\prime}} \oplus V _1 ^{k} )$ isomorphic to $V _1 ^{k ^{\prime}}$ or $V _1 ^{k}$, respectively. When $\lambda _k = \lambda _{k ^{\prime}}$, we use this to change our normal form so that the corresponding marked partition satisfies \ref{snf} 3). When $\lambda _k > \lambda _{k ^{\prime}}$, we use this to transform our marked partition into another marked partition which satisfies \ref{snf} 4) for the pair $( J _k, J _{k ^{\prime}} )$. If we make changes to our marked partition in one of the above two procedures, then $\mathbf J$ is unchanged, $\delta _1$ on one of $\{ J _k, J _{k ^{\prime}} \}$ is unchanged, but $\delta _1$ on the other becomes $0$. By repeating these procedures for every possible pair $k < k ^{\prime}$, we complete the proof of the first assertion.\\
For the second assertion, recall that $\mathsf{Irrep} W$ is parametrized by the set of ordered pair of partitions $( \lambda ^1, \lambda ^2 )$ which sum up to $n$. We define two-partitions out of a strict marked partition $\sigma$ as
$$\lambda ^1 _k + \lambda ^2 _k = \lambda _k, \text{ and } \lambda ^2 _k := \begin{cases} \underline{\#} J _k & (\mathbf v _{1, \sigma} ^{J _k} \neq 0)\\ \max \{0, \underline{\#} J _{k ^{\prime}}, \lambda _k - \overline{\#} J _{k ^{\prime\prime}} ; k ^{\prime} > k > k ^{\prime\prime} \} & (otherwise)\end{cases}$$
for each $k$, where the set we choose its maximal is formed only from these $J _{k ^{\prime}}$ and $J _{k ^{\prime\prime}}$ for which $\overline{\#}$ and $\underline{\#}$ are defined. It is clear that two sequences $\lambda ^1, \lambda ^2$ sum up to $n$. By \ref{snf} 4), we deduce that
$$\lambda _k - \overline{\#} J _{k ^{\prime\prime}} < \underline{\#} J _{k} \text{ (this is equivalent to }\overline{\#} J _{k ^{\prime\prime}} > \overline{\#} J _{k})$$
holds for $k ^{\prime\prime} < k$ (such that $\overline{\#}$ and $\underline{\#}$ are defined for both $J _k$ and $J _{k ^{\prime\prime}}$). It follows that $\lambda ^2$ is a partition. (I.e. $\{ \lambda _i ^2 \} _i$ is a decreasing sequence.) By the symmetry of $\overline{\#}$ and $\underline{\#}$ in \ref{snf}, we conclude that $\lambda ^1$ is also a partition.\\
Therefore, it suffices to prove that the pairs of partitions formed by strict marked partitions are equal only if the marked partitions are equal. (Since this gives the injectivity of the above assignment.) For this, we assume that two strict marked partitions $\sigma = ( \mathbf J, \delta _1 )$ and $\sigma ^{\prime} = ( \mathbf J ^{\prime}, \delta _1 ^{\prime} )$ gives the same pair $( \lambda ^1, \lambda ^2 )$ to deduce contradiction. We can assume that $\mathbf J = \mathbf J ^{\prime}$ since $\lambda = \lambda ^1 + \lambda ^2$. Hence, their difference is concentrated in their foot function. By \ref{snf} 3) and 4), we deduce that the foot functions are non-trivial on $J _k = J _{k}  ^{\prime}$ if and only if
$$\lambda ^2 _k \neq \max \{ \lambda ^2 _j, \lambda _k - \lambda ^1 _i ;  \lambda _j \neq \lambda _k \neq \lambda _i, i < k < j \}$$
and $\lambda _k \neq \lambda _{k - 1}$. Moreover, the value of the foot functions on $J _k$ are determined by the value of $\lambda ^2 _k$ if they are non-trivial. Since this system has a unique solution, we deduce $\sigma = \sigma ^{\prime}$, which is contradiction. Thus, the pair of partitions recovers a strict marked partition uniquely, which completes the proof of the second assertion.
\end{proof}

\begin{theorem}[Orbit structure of $\mathfrak N _2 ^a$]\label{EJNF}
Let $\nu = (a, X) = ( s, \vec{q}, X )$ be an admissible parameter. Then, there exists $g \in G$ such that:
$$g s g ^{-1} \in T \text{ and } g X \text{ is a normal form.}$$
\end{theorem}

\begin{proof}
Postponed to \S \ref{sten}.
\end{proof}

We have a natural $W$-action $\cdot$ on $[-n, n] ^*$ by setting
$$s _i \cdot j := \begin{cases} \pm ( i + 1 ) & (j = \pm i) \\ \pm i & (j = \pm ( i + 1 )) \\ j & (otherwise) \end{cases} \text{ for } i = 1, \ldots, n - 1, \text{ and } s _n \cdot j = \begin{cases} - j & (j = \pm n) \\ j & (otherwise) \end{cases}.$$
Using this, we define the $W$-action $\cdot$ on the set of $\ell$-marked partitions as:\\
For $w \in W$ and $\sigma = ( \mathbf J, \vec{\delta} ) = ( \{ J _1, J _2, \ldots, \}, \{ \delta _1, \ldots, \delta _{\ell} \})$, we set
$$w \cdot \sigma := ( \{ w \cdot J _1, w \cdot J _2, \ldots, \}, \{ w \cdot \delta _1, \ldots, w \cdot \delta _{\ell} \}),$$
where we set
$$w \cdot ( J _1 ^1, J _1 ^2 ,\ldots ) = ( ( w \cdot J ) _1 ^1, ( w \cdot J ) _1 ^2 ,\ldots ) := ( w \cdot J _1 ^1, w \cdot J _1 ^1, \ldots )$$
and $w \cdot \delta _k ( j ) := \delta _k ( w ^{-1} \cdot j )$. Notice that we have $w \cdot J _k = ( w \cdot J ) _k$ and $w \cdot J _k ^j = ( w \cdot J ) _k ^j$ for every $k, j$ in this action.

\begin{lemma}\label{stabTJ}
Let $\sigma = ( \mathbf J, \vec{\delta} )$ be a marked partition which is a $W$-translation of a strict marked partition. Then, we have
$$\mathbb C ^{\times} \mathbf v _{1, \sigma} \oplus \mathbb C ^{\times} \mathbf v _{2, \sigma} \subset T \mathbf v _{\sigma}, \text{ and } \mathbb C ^{\times} \mathbf v _{1, \sigma} ^J \oplus \mathbb C ^{\times} \mathbf v _{2, \sigma} ^J \subset T _J \mathbf v _{\sigma} ^J \text{ for each }J \in \mathbf J.$$
\end{lemma}

\begin{proof}
Since $T _J \cap T _{J ^{\prime}} = \{ 1 \}$, it suffices to prove the second assertion. Let $\mathbf v _{\sigma} ^{J} = \sum _{\xi \in \Xi} v _{\xi}$ be the $T$-eigen-decomposition of $\mathbf v _{\sigma} ^{J}$. Then, we have $\# \Xi = \# J$ and $\dim T _J = \# J$. Moreover, the weights appearing in $\Xi$ are linearly independent. Hence, we have the scalar multiplications of each $v _{\xi}$, which implies the result.
\end{proof}

\begin{corollary}[of the proof of Lemma \ref{stabTJ}]
Let $\sigma$ be a strict marked partition. Let $w \in W$. Then, we have  $\mathbf v _{w \cdot \sigma} \in G \mathbf v _{\sigma}$. \hfill $\Box$
\end{corollary}

\subsection{Structure of simple modules}\label{gen conv}
We put $T _{\ell} := T \times ( \mathbb C ^{\times} ) ^{\ell + 1}$. Let $a \in T _{\ell}$. Let $\mu ^a : F _{\ell} ^a \rightarrow \mathfrak N _{\ell} ^a$ denote the restriction of $\mu _{\ell}$ to $a$-fixed points.\\
We review the convolution realization of simple modules in our situation. The detailed constructions are found in \cite{CG} 5.11, 8.4 or \cite{Gi} \S 5. For its variant, see \cite{J1}.\\
The properties we used to apply the Ginzburg theory are: 1) $Z _{\ell} = F _{\ell} \times _{\mathfrak N _{\ell}} F _{\ell}$; 2) $F _{\ell}$ is smooth; 3) $\mu _{\ell}$ is projective; 4) $R ( G _{\ell} ) \subset K ^{G _{\ell}} ( Z _{\ell} )$ is central; and 5) $H _{\bullet} ( Z _{\ell} )$ is spanned by algebraic cycles.\\
Let $\mathbb C _{a}$ be the unique residual field of $\mathbb C \otimes _{\mathbb Z} R ( G _{\ell} ) _a$ or $\mathbb C \otimes _{\mathbb Z} R ( T _{\ell} ) _a$. The Thomason localization theorem yields ring isomorphisms
$$\mathbb C _{a} \otimes _{R ( G _{\ell} )} K ^{G _{\ell}} ( Z _{\ell} ) \stackrel{\cong}{\longrightarrow} \mathbb C _{a} \otimes _{R ( G _{\ell} ( a ) )} K ^{G _{\ell} (a )} ( Z _{\ell} ^{a} ) \stackrel{\cong}{\longrightarrow} \mathbb C _{a} \otimes _{R ( T _{\ell} )} K ^{T _{\ell}} ( Z _{\ell} ^{a} ).$$
Moreover, we have the Riemann-Roch isomorphism
$$\mathbb C _{a} \otimes _{R ( T _{\ell} )} K ^{T _{\ell}} ( Z _{\ell} ^{a} ) \cong K ( Z _{\ell} ^{a} ) \stackrel{RR}{\longrightarrow} H _{\bullet} ( Z _{\ell} ^{a} )\cong \mathrm{Ext} ^{\bullet} ( \mu ^{a} _* \mathbb C _{F _{\ell} ^{a}}, \mu ^{a} _* \mathbb C _{F _{\ell} ^{a}} ).$$
By the equivariant Beilinson-Bernstein-Deligne (-Gabber) decomposition theorem (c.f. Saito \cite{S} 5.4.8.2), we have
$$\mu ^{a} _* \mathbb C _{F _{\ell} ^{a}} \cong \bigoplus _{\mathbb O \subset \mathfrak N _{\ell} ^{a}, \chi, d} L _{\mathbb O, \chi, d} \boxtimes IC ( \mathbb O, \chi ) [ d ],$$
where $\mathbb O \subset \mathfrak N _{\ell} ^{a}$ is a $G ( s )$-stable subset such that $\mu ^{a}$ is locally trivial along $\mathbb O$, $\chi$ is an irreducible local system on $\mathbb O$, $d$ is an integer, $L _{\mathbb O, \chi, d}$ is a finite dimensional vector space, and $IC ( \mathbb O, \chi )$ is the minimal extension of $\chi$. Moreover, the set of $\mathbb O$'s such that $L _{\mathbb O, \chi, d} \neq 0$ (for some $\chi$ and $d$) forms a subset of an algebraic stratification in the sense of \cite{CG} 3.2.23. It follows that:

\begin{theorem}[Ginzburg \cite{Gi} Theorem 5.2]\label{Gd}
The set of simple modules of $K ^{G _{\ell}} ( Z _{\ell} )$ for which $R ( G _{\ell} )$ acts as the evaluation at $a$ is in one-to-one correspondence with the set of isomorphism classes of irreducible $G _{\ell} ( a )$-equivariant perverse sheaves appearing in $\mu ^{a} _* \mathbb C _{F _{\ell} ^{a}}$ $($up to degree shift$)$. \hfill $\Box$
\end{theorem}

\section{Hecke algebras and exotic nilpotent cones}\label{fHecke}
We retain the setting of the previous section. We put ${\bm G} = G _2$, ${\bm T} := T _2$, ${\bm G} := F _2$, ${\bm \mu} := \mu _2$, ${\bm Z} := Z _2$, and ${\bm \pi} := \pi _2$. Most of the arguments in this section are exactly the same as \cite{CG} 7.6 if we replace ${\bm G}$ by $G \times \mathbb C ^{\times}$, $\mathfrak N _2$ by the usual nilpotent cone, ${\bm \mu}$ by the moment map, ${\bm F}$ by the cotangent bundle of the flag variety, and ${\bm Z}$ by the Steinberg variety. Therefore, we frequently omit the detail and make pointers to \cite{CG} 7.6 in which the reader can obtain a correct proof merely replacing the meaning of symbols as mentioned above.

We put $\mathcal A _{\mathbb Z} := \mathbb Z [ \mathbf q _0 ^{\pm 1}, \mathbf q _1 ^{\pm 1}, \mathbf q _2 ^{\pm 1} ]$ and $\mathcal A := \mathbb C \otimes _{\mathbb Z} \mathcal A _{\mathbb Z} = \mathbb C [ \mathbf q _0 ^{\pm 1}, \mathbf q _1 ^{\pm 1}, \mathbf q _2 ^{\pm 1} ]$.

\begin{definition}[Hecke algebras of type $C _n ^{( 1 )}$]\label{Checke}
A Hecke algebra of type $C _n ^{( 1 )}$ with three parameters is an associative algebra $\mathbb H$ over $\mathcal A$ generated by $\{ T _i \} _{i = 1} ^n$ and $\{ e ^{\lambda} \} _{\lambda \in X ^* ( T ) }$ subject to the following relations:
\item {\bf (Toric relations)} For each $\lambda, \mu \in X ^* ( T )$, we have $e ^{\lambda} \cdot e ^{\mu} = e ^{\lambda + \mu}$ (and $e ^0 = 1$);
\item {\bf (The Hecke relations)} We have
$$( T _i + 1 ) ( T _i - \mathbf q _2 ) = 0 \text{ ($1 \le i < n$)  and  } ( T _n + 1 ) ( T _n + \mathbf q _0 \mathbf q _1 ) = 0;$$
\item {\bf (The braid relations)} We have
\begin{center}
$T _i T _j = T _j T _i$ (if $\left| i - j \right| > 1$), $( T _n T _{n - 1} ) ^2 = ( T _{n - 1} T _n ) ^2$,\\
$T _i T _{i + 1} T _i = T _{i + 1} T _i T _{i + 1}$ (if $1 \le i < n - 1$);
\end{center}
\item {\bf (The Bernstein-Lusztig relations)} For each $\lambda \in X ^{*} ( T )$, we have
$$T _i e ^{\lambda} - e ^{s _i \lambda} T _i = \begin{cases} ( 1 - \mathbf q _2 ) \frac{e ^{\lambda} - e ^{s _i \lambda}}{e ^{\alpha _i} - 1} & \text{ ($i \neq n$)}\\ \frac{( 1 + \mathbf q _0 \mathbf q _1 ) - ( \mathbf q _0 + \mathbf q _1 ) e ^{\epsilon _n}}{e ^{\alpha _n} - 1} ( e ^{\lambda} - e ^{s _n \lambda} ) & \text{ ($i = n$)} \end{cases}.$$
\end{definition}

\begin{remark}\label{remHecke}
{\bf 1)} The standard choice of parameters $( t _0, t _1, t _n )$ is: $t _1 ^2 = \mathbf q _2$, $t _n ^2 = - \mathbf q _0 \mathbf q _1$, and $t _n ( t _0 - t _0 ^{- 1} ) = ( \mathbf q _0 + \mathbf q _1 )$. This yields
$$T _n e ^{\lambda} - e ^{s _n \lambda} T _n =  \frac{1 - t _n ^2 - t _n ( t _0 - t _0 ^{- 1} ) e ^{\epsilon _n}}{e ^{2 \epsilon _n} - 1} ( e ^{\lambda} - e ^{s _n \lambda} );$$
{\bf 2)} If $n = 1$, then we have $T _1 = T _n$ in Definition \ref{Checke}. In this case, we have $\mathbb H \cong \mathbb C [ \mathbf q _2 ^{\pm 1} ] \otimes _{\mathbb C} \mathbb H _0$, where $\mathbb H _0$ is the Hecke algebra of type $A _1 ^{(1)}$ with two-parameters $( \mathbf q _0, \mathbf q _1 )$;\\
{\bf 3)} An extended Hecke algebra of type $B _n ^{( 1 )}$ with two-parameters considered in \cite{E} is obtained by requiring $\mathbf q _0 + \mathbf q _1 = 0$. An equal parameter extended Hecke algebra of type $B _n ^{(1)}$ is obtained by requiring $\mathbf q _0 + \mathbf q _1 = 0$ and $\mathbf q _1 ^2 = \mathbf q _2$. An equal parameter Hecke algebra of type $C _n ^{(1)}$ is obtained by requiring $\mathbf q _2 = - \mathbf q _0 \mathbf q _1$ and $( 1 + \mathbf q _0 ) ( 1 + \mathbf q _1 ) = 0$.
\end{remark}

For each $w \in W$, we define two closed subvarieties of ${\bm Z}$ as
$$Z _{\le w} := {\bm \pi} ^{-1} ( \overline{\mathsf O _w} ) \text{ and } Z _{< w} := Z _{\le w} \backslash {\bm \pi} ^{-1} ( \mathsf O _w ).$$

Let $\lambda \in X ^* ( T )$. Let $\mathcal L _{\lambda}$ be the pullback of the line bundle $G \times ^{B} \lambda ^{- 1}$ over $G / B$ to ${\bm F}$. Clearly $\mathcal L _{\lambda}$ admits a $\mathbf G$-action by letting $( \mathbb C ^{\times} ) ^3$ act on $\mathcal L _{\lambda}$ trivially. We denote the operator $[\tilde{p} _1 ^* \mathcal L _{\lambda} \otimes ^{\mathbb L}\bullet ]$ by ${\mathtt e} ^{\lambda}$. By abuse of notation, we may denote ${\mathtt e} ^{\lambda} ( 1 )$ by ${\mathtt e} ^{\lambda}$ (in $K ^{\mathbf G} ( {\bm Z} )$). Let $\mathtt q _0 \in R ( \{ 1 \} \times \mathbb C ^{\times} \times \{ 1 \} \times \{ 1 \}  ) \subset R ( {\bm G} )$, $\mathtt q _1 \in R ( \{ 1 \} \times \{ 1 \} \times \mathbb C ^{\times} \times \{ 1 \}  ) \subset R ( {\bm G} )$, and $\mathtt q _2 \in R ( \{ 1 \} \times \{ 1 \} \times \{ 1 \} \times \mathbb C ^{\times} ) \subset R ( {\bm G} )$ be the inverse of degree-one characters. (I.e. $\mathtt q _2$ corresponds to the inverse of the scalar multiplication on $V _2$.) By the operation ${\mathtt e} ^{\lambda}$ and the multiplication by $\mathtt q _i$, each of $K ^{{\bm G}} ( Z _{\le w} )$ admits a structure of $R ({\bm T} )$-modules.

Each $Z _{\le w} \backslash Z _{< w}$ is a ${\bm G}$-equivariant vector bundle over an affine fibration over $G / B$ via the composition of $ {\bm \pi}$ and the second projection. Therefore, the cellular fibration Lemma (or the successive application of localization sequence) yields:

\begin{theorem}[c.f. \cite{CG} 7.6.11]\label{cell}
We have
$$K ^{\bm G} ( Z _{\le w} ) = \bigoplus _{v \in W ; \mathsf O _v \subset \overline{\mathsf O} _w} R ( {\bm T} ) [ \mathcal O _{Z _{\le v}} ].$$
\end{theorem}

For each $i = 1, 2, \ldots, n$, we put $\mathbb O _{i} := \overline{{\bm \pi} ^{-1} ( \mathsf O _{s _i} )}$. We define $\tilde{T} _i := [ \mathcal O _{\mathbb O _{i}} ]$ for each $i = 1, \ldots, n$. 

\begin{theorem}[c.f. Proof of \cite{CG} 7.6.12]\label{genH}
The set $\{ [ \mathcal O _{Z _{\le 1}} ], \tilde{T} _i, {\mathtt e} ^{\lambda} ; 1 \le i \le n, \lambda \in X ^* ( T ) \}$ is a generator set of $K ^{\bm G} ( {\bm Z} )$ as $\mathcal A _{\mathbb Z}$-algebras.
\end{theorem}
\begin{proof}
The tensor product of structure sheaves corresponding to vector subspaces of a vector space is the structure sheaf of their intersection. Taking account into that, the proof of the assertion is exactly the same as \cite{CG} 7.6.12.
\end{proof}

By the Thom isomorphism, we have an identification
\begin{eqnarray}
K ^{\bm G} ( {\bm F} ) \cong K ^{\bm G} ( G / B ) \cong R ( {\bm T} ) = \mathcal A _{\mathbb Z} [ T ].\label{Thom}
\end{eqnarray}
We normalize the images of $[ \mathcal L _{\lambda} ]$ and $\mathtt q _i$  ($i = 0, 1, 2$) under (\ref{Thom}) as $e ^{\lambda}$ and $\mathbf q _i$, respectively.
\begin{theorem}[c.f. \cite{CG} Claim 7.6.7]\label{inclbasic}
The homomorphism
$$\circ : K ^{\bm G} ( {\bm Z} ) \longrightarrow \mathrm{End} _{R ( {\bm G} )} K ^{\bm G} ( {\bm F} )$$
is injective. \hfill $\Box$
\end{theorem}

\begin{proposition}\label{simple}
We have
\begin{enumerate}
\item $[ \mathcal O _{Z _{\le 1}} ] = 1 \in \mathrm{End} _{R ( {\bm G} )} K ^{\bm G} ( {\bm F} )$;
\item $\tilde{T} _i \circ e ^{\lambda} = ( 1 - \mathbf q _2 e ^{\alpha _i} ) \frac{e ^{\lambda} - e ^{s _i \lambda - \alpha _i}}{1 - e ^{- \alpha _i}}$ for every $\lambda \in X ^* ( T )$ and every $1 \le i < n$;
\item $\tilde{T} _n \circ e ^{\lambda} = ( 1 - \mathbf q _0 e ^{\frac{1}{2} \alpha _n} ) ( 1 - \mathbf q _1 e ^{\frac{1}{2} \alpha _n} ) \frac{e ^{\lambda} - e ^{s _n \lambda - \alpha _n}}{1 - e ^{- \alpha _n}}$ for every $\lambda \in X ^* ( T )$.
\end{enumerate}
\end{proposition}

\begin{proof}
The component $Z _{\le 1}$ is equal to the diagonal embedding of ${\bm F}$. In particular, both of the first and the second projections give isomorphisms between $Z _{\le 1}$ and ${\bm F}$. It follows that 
\begin{eqnarray*}
[ \mathcal O _{Z _{\le 1}} ] \circ [ \mathcal L _{\lambda} ] = \sum _{i \ge 0} ( - 1 ) ^i [ \mathbb R ^i ( p _1 ) _* \left( \mathcal O _{Z _{\le 1}} \otimes ^{\mathbb L} \tilde{p} _2 ^* \mathcal L _{\lambda} \right) ] \\
= [ \mathbb R ^0 ( p _1 ) _* \left( \mathcal O _{Z _{\le 1}} \otimes \tilde{p} _2 ^* \mathcal L _{\lambda} \right) ] = [ \mathcal L _{\lambda} ],
\end{eqnarray*}
which proves 1). For each $i = 1, \ldots, n$, we define $\mathbb V ^+ ( i ) := \mathbb V ^+ _2 \cap \dot{s} _i \mathbb V ^+ _2$. Let $P _i := B \dot{s} _i B \sqcup B$ be a parabolic subgroup of $G$ corresponding to $s _i$. Each $\mathbb V ^+ ( i )$ is $B$-stable. Hence, it is $P _i$-stable. We have
$${\bm \pi} ( \mathbb O _i ) = \overline{\mathsf O} _{s _i} = G ( 1 \times P _i ) \mathsf p _1 \subset G / B \times G / B.$$
The product $( 1 \times P _i ) \mathsf p _1 \times \mathbb V ^+ ( i )$ is a $B$-equivariant vector bundle. Here we have $G \cap ( B \times P _i ) = B$. Hence, we can induce it up to a $G$-equivariant vector bundle $\tilde{\mathbb V} ( i )$ on ${\bm \pi} ( \mathbb O _i )$. By means of the natural embedding of $G$-equivariant vector bundles
$${\bm F} = G \times ^B \mathbb V _2 ^+ \hookrightarrow G \times ^B \mathbb V _2 \cong G \times \mathbb V _2,$$
we can naturally identify ${\bm \pi} ^{-1} ( \mathsf{p} _{s _i} )$ with $\mathbb V ^+ ( i )$. Since $\mathbb V ^+ ( i )$ is $P _i$-stable, we conclude ${\bm \pi} ^{-1} ( \mathsf{p} _{s _i} ) \cong \mathbb V ^+ ( i )$ as $P _i$-modules. As a consequence, we conclude $\tilde{\mathbb V} ( i ) \cong \mathbb O _i$. Let $\breve{F} ( i ) := G \times ^B ( \mathbb V ^+ _2 / \mathbb V ^+ ( i ) )$. It is a ${\bm G}$-equivariant quotient bundle of ${\bm F}$. The rank of $\breve{F} ( i )$ is one ($1 \le i < n$) or two ($i = n$). Let $\breve{Z} _{\le s _i}$ be the image of $Z _{\le s _i}$ under the quotient map ${\bm F} \times {\bm F} \rightarrow \breve{F} ( i ) \times \breve{F} ( i )$. We obtain the following commutative diagram:
\begin{center}
\xymatrix{
{\bm F} \ar[d] & & Z _{\le s _i} \ar[ll] \ar[rr] \ar[d] & & {\bm F} \ar[d] \\
\breve{F} ( i ) & & \breve{Z} _{\le s _i} \ar[ll] \ar[rr] & & \breve{F} ( i )}
\end{center}
Here the above objects are smooth $\mathbb V ^+ ( i )$-fibrations over the bottom objects. Therefore, it suffices to compute the convolution operation of the bottom line. We have $\breve{Z} _{\le s _i} = \overline{\mathsf O} _{s _i} \cup \triangle ( \breve{F} ( i ) )$, where $\triangle : \breve{F} ( i ) \hookrightarrow \breve{F} ( i ) ^2$ is the diagonal embedding. Let $\breve{p} _j : \overline{\mathsf O} _{s _i} \rightarrow G / B$ ($j = 1, 2$) be projections induced by the natural projections of $G / B \times G / B$. By construction, each $\breve{p} _j$ is a $\bm G$-equivariant $\mathbb P ^1$-fibration. Let $\breve{\mathcal L} _{\lambda}$ be the pullback of $G \times ^{B} \lambda ^{- 1}$ to $\breve{F}  ( i )$. We deduce
\begin{eqnarray*}
\tilde{T} _i \circ [ \breve{\mathcal L} _{\lambda} ] = \sum _{i \ge 0} ( - 1 ) ^i [ \mathbb R ^i ( \breve{p} _1 ) _* ( \mathcal O _{\overline{\mathsf O} _{s _i}} \otimes ^{\mathbb L} ( \mathcal O _{\breve{F}  ( i )} \boxtimes \breve{\mathcal L} _{\lambda} ) ]\\
= \sum _{i \ge 0} ( - 1 ) ^i [ \mathbb R ^{i} ( \breve{p} _1 ) _* \breve{p} _2 ^* \imath _* ( G \times ^B \lambda ^{- 1} ) ] = \left[ G \times ^B [ \frac{e ^{\lambda} - e ^{s _i \lambda - \alpha _i}}{1 - e ^{- \alpha _i}} ] \right],
\end{eqnarray*}
where $\imath : G / B \hookrightarrow \breve{F} ( i )$ is the zero section, and $[ \frac{e ^{\lambda} - e ^{s _i \lambda - \alpha _i}}{1 - e ^{- \alpha _i}} ] \in R ( T ) \cong R ( B )$ is a virtual $B$-module. Here the ideal sheaf associated to $G / B \subset \breve{F} ( i )$ represents $\mathbf q _2 [ \breve{\mathcal L} _{\alpha _i} ]$ in $K ^{\bm G} ( \breve{F} ( i ) )$ ($1 \le i < n$) or corresponds to $\mathbf q _0 \breve{\mathcal L} _{\epsilon _n} + \mathbf q _1 \breve{\mathcal L} _{\epsilon _n} \subset \mathcal O _{\breve{F} ( i )}$ ($i = n$). In the latter case, divisors corresponding to $\mathbf q _0 \breve{\mathcal L} _{\epsilon _n}$ and $\mathbf q _1 \breve{\mathcal L} _{\epsilon _n}$ are normal crossing. Thus, we have $[ \mathbf q _0 \breve{\mathcal L} _{\epsilon _n} \cap \mathbf q _1 \breve{\mathcal L} _{\epsilon _n} ] = \mathbf q _0 \mathbf q _1 [ \breve{\mathcal L} _{2 \epsilon _n} ]$. In particular, we deduce
$$[ \mathbf q _0 \breve{\mathcal L} _{\epsilon _n} + \mathbf q _1 \breve{\mathcal L} _{\epsilon _n} ] = \mathbf q _0 [ \breve{\mathcal L} _{\epsilon _n}]  + \mathbf q _1 [ \breve{\mathcal L} _{\epsilon _n} ] - \mathbf q _0 \mathbf q _1 [ \breve{\mathcal L} _{2 \epsilon _n} ] \in K ^{\bm G} ( \breve{F} ( n ) ).$$ Therefore, we conclude
$$\tilde{T} _i \circ e ^{\lambda} = \begin{cases} 
( 1 - \mathbf q _2 e ^{\alpha _i} ) \frac{e ^{\lambda} - e ^{s _i \lambda - \alpha _i}}{1 - e ^{- \alpha _i}} & \text{($1 \le i < n$)}\\
( 1 - \mathbf q _0 e ^{\frac{\alpha _n}{2}} ) ( 1 - \mathbf q _1 e ^{\frac{\alpha _n}{2}} ) \frac{e ^{\lambda} - e ^{s _n \lambda - \alpha _n}}{1 - e ^{- \alpha _n}} & \text{($i = n$)}
\end{cases}$$
as desired.
\end{proof}

The following representation of $\mathbb H$ is usually called the basic representation or the anti-spherical representation:

\begin{theorem}[Basic representation c.f. \cite{M} 4.3.10]\label{basic}
There is an injective $\mathcal A$-algebra homomorphism
$$\varepsilon : \mathbb H \rightarrow \mathrm{End} _{\mathcal A} \mathcal A [ T ],$$
defined as $\varepsilon ( e ^{\lambda} ) := e ^{\lambda} \cdot$ $(\lambda \in X ^* ( T ))$ and
$$
\varepsilon ( T _i ) e ^{\lambda} := \begin{cases} \frac{e ^{\lambda} - e ^{s _i \lambda}}{e ^{\alpha _i} - 1} - \mathbf q _2 \frac{e ^{\lambda} - e ^{s _i \lambda + \alpha _i}}{e ^{\alpha _i} - 1} & (\text{$1 \le i < n$}) \\ \frac{e ^{\lambda} - e ^{s _n \lambda}}{e ^{\alpha _n} - 1} + \mathbf q _0 \mathbf q _1 \frac{e ^{\lambda} - e ^{s _n \lambda + \alpha _n}}{e ^{\alpha _n} - 1} - ( \mathbf q _0 + \mathbf q _1 ) e ^{\epsilon _n} \frac{e ^{\lambda} - e ^{s _n \lambda}}{e ^{\alpha _n} - 1} & (\text{$i = n$}) \end{cases}.
$$
\end{theorem}

\begin{theorem}[Exotic geometric realization of Hecke algebras]\label{ehecke}
We have an isomorphism
$$\mathbb H \stackrel{\cong}{\longrightarrow} \mathbb C \otimes _{\mathbb Z} K ^{\bm G} ( {\bm Z} ),$$
as algebras.
\end{theorem}

\begin{proof}
Consider an assignment $\vartheta$
$$\underline{e} ^{\lambda} \mapsto \mathsf e ^{\lambda}, \underline{T _i} \mapsto \begin{cases}
\tilde{T} _i - ( 1 - \mathtt q _2 ( {\mathsf e} ^{\alpha _i} + 1 )) & (1 \le i < n)\\
 \tilde{T} _i + ( \mathtt q _0 + \mathtt q _1 ) \mathsf e ^{\epsilon _n} - ( 1 + \mathtt q _0 \mathtt q _1 ( {\mathsf e} ^{\alpha _n} + 1 )) & (i = n) \end{cases}.$$
By means of the Thom isomorphism, the above assignment gives an action of an element of the set $\{ \underline{e} ^{\lambda} \} \cup \{ \underline{T _i} \} _{i = 1} ^n$ on $\mathcal A [ T ]$. We have
\begin{align*}
\vartheta ( \underline{e} ^{\lambda} ) e ^{\mu} = & e ^{\lambda + \mu}\\
\vartheta ( \underline{T _i} ) e ^{\lambda} = & \left( \tilde{T} _i - ( 1 - \mathbf q _2 ( e ^{\alpha _i} + 1 )) \right) e ^{\lambda} = ( 1 - \mathbf q _2 e ^{\alpha _i} ) \frac{e ^{\lambda} - e ^{s _i \lambda - \alpha _i}}{1 - e ^{- \alpha _i}} - e ^{\lambda} + \mathbf q _2 ( e ^{\alpha _i} + 1 ) e ^{\lambda}\\
= & ( \frac{e ^{\lambda} - e ^{s _i \lambda - \alpha _i}}{1 - e ^{- \alpha _i}} - \frac{e ^{\lambda} - e ^{\lambda - \alpha _i}}{1 - e ^{- \alpha _i}}) - \mathbf q _2 e ^{\alpha _i} ( \frac{e ^{\lambda} - e ^{s _i \lambda - \alpha _i}}{1 - e ^{- \alpha _i}} - \frac{e ^{\lambda} - e ^{\lambda - 2 \alpha _i}}{1 - e ^{- \alpha _i}} ) = \varepsilon ( T _i ) e ^{\lambda}\\
\vartheta ( \underline{T _n} ) e ^{\lambda} = & \left( \tilde{T} _n + ( \mathbf q _0 + \mathbf q _1 ) e ^{\epsilon _n} - ( 1 + \mathbf q _0 \mathbf q _1 ( e ^{\alpha _n} + 1 )) \right) e ^{\lambda}\\
= & ( 1 - \mathbf q _0 e ^{\epsilon _n} ) ( 1 - \mathbf q _1 e ^{\epsilon _n} )  \frac{e ^{\lambda} - e ^{s _n \lambda - \alpha _n}}{1 - e ^{- \alpha _n}} - e ^{\lambda} + ( \mathbf q _0 + \mathbf q _1 ) e ^{\lambda + \epsilon _n} - \mathbf q _0 \mathbf q _1 ( e ^{\alpha _n} + 1 ) e ^{\lambda}\\
= & ( \frac{e ^{\lambda} - e ^{s _n \lambda - \alpha _n}}{1 - e ^{- \alpha _n}} - \frac{e ^{\lambda} - e ^{\lambda - \alpha _n}}{1 - e ^{- \alpha _n}}) + \mathbf q _0 \mathbf q _1  e ^{\alpha _n} ( \frac{e ^{\lambda} - e ^{s _n \lambda - \alpha _n}}{1 - e ^{- \alpha _n}} - \frac{e ^{\lambda} - e ^{\lambda - 2 \alpha _n}}{1 - e ^{- \alpha _n}} )\\
 & - ( \mathbf q _0 + \mathbf q _1 ) ( \frac{e ^{\lambda + \epsilon _n} - e ^{s _n \lambda - \epsilon _n}}{1 - e ^{- \alpha _n}} - \frac{e ^{\lambda + \epsilon _n} - e ^{\lambda - \epsilon _n}}{1 - e ^{- \alpha _n}} ) = \varepsilon ( T _n ) e ^{\lambda}.
\end{align*}
This identifies $\mathbb C \otimes _{\mathbb Z} K ^{\bm G} ( {\bm F} )$ with the basic representation of $\mathbb H$ via the correspondence $e ^{\lambda} \mapsto \underline{e} ^{\lambda}$ and $T _i \mapsto \underline{T _i}$. In particular, it gives an inclusion $\mathbb H \subset \mathbb C \otimes _{\mathbb Z} K ^{\bm G} ( {\bm Z} )$. Here we have $T _i \in \tilde{T} _i + R ( {\bm T} )$ for $1 \le i \le n$. It follows that $\mathbb C \otimes _{\mathbb Z} K ^{\bm G} ( {\bm Z} ) \subset \mathbb H$, which yields the result.
\end{proof}

\begin{theorem}[Bernstein c.f. \cite{CG} 7.1.14 and \cite{M} 4.2.10]\label{bcenter}
The center $Z ( \mathbb H )$ of $\mathbb H$ is naturally isomorphic to $\mathbb C \otimes _{\mathbb Z} R ( \mathbf G )$. \hfill $\Box$
\end{theorem}

\begin{corollary}
The center of $K ^{\bm G} ( {\bm Z} )$ is $R ( {\bm G} )$. \hfill $\Box$
\end{corollary}

For a semisimple element $a \in \bm G$, we define
$$\mathbb H _{a} := \mathbb C _{a} \otimes _{Z ( \mathbb H )} \mathbb H \quad \text{ (c.f. \S \ref{gen conv})}$$
and call it the specialized Hecke algebra.

\begin{theorem}\label{HKmain}
Let $a \in \bm G$ be a semisimple element. We have an isomorphism
$$\mathbb H _{a} \cong \mathbb C \otimes _{\mathbb Z} K ( {\bm Z} ^{a} )$$
as algebras.
\end{theorem}

\begin{proof}
This is a combination of \cite{CG} 6.2.3 and 5.10.11. (See also \cite{CG} 8.1.6.)
\end{proof}

\begin{convention}\label{ssnoneq}
Let $a = ( s, \vec{q} ) \in \bm G$ be a pre-admissible element. We define $Z _+ ^{a}$ to be the image of ${\bm Z} ^{a}$ under the natural projection defined by
$${\bm Z} \ni ( g _1 B, g _2 B, X _0, X _1, X _2 ) \mapsto ( g _1 B, g _2 B, X _0 + X _1, X _2 ) \in Z.$$
Let $F _+ ^{a}$ be the image of $Z _+ ^{a}$ via the first (or the second) projection. Let $\mu _+ ^{a}$ be the restriction of $\mu$ to $F _+ ^{a}$. We denote its image by $\mathfrak N ^{a} _+$.
By the assumption $q _0 \neq q _1$, we have $F _+ ^{a} \cong {\bm F} ^{a}$, $Z _+ ^{a} \cong {\bm Z} ^{a}$, and $\mathfrak N _+ ^{a} \cong \mathfrak N _2 ^{a}$.
\end{convention}

\begin{corollary}\label{SH}
Keep the setting of Convention \ref{ssnoneq}. We have an isomorphism
$$\mathbb H _{a} \cong \mathbb C \otimes _{\mathbb Z} K ( Z ^{a} _+ )$$
as algebras. \hfill $\Box$
\end{corollary}

For the later use, let us introduce our last class of parameter here.

\begin{definition}[Regular parameters]\label{lan par}
A pre-admissible parameter $( a, X )$ is called regular iff there exists a direct factor $A [ d ] \subset ( \mu ^{a} _+ ) _* \mathbb C _{F ^{a} _+}$, where $A$ is a simple ${\bm G} ( a )$-equivariant perverse sheaf on $\mathfrak N ^{a} _+$ such that $\mathrm{supp} A = \overline{{\bm G} ( a ) X}$ and $d$ is an integer.\\
We denote by $\mathfrak R _a$ the set of ${\bm G} ( a )$-conjugacy classes of regular pre-admissible parameters of the form $( a, X )$ ($X \in \mathfrak N ^a _+$).
\end{definition}

\section{Clan decomposition}\label{clan}
We work under the same setting as in \S \ref{fHecke}.

\begin{definition}[Clans]\label{LC}
Let $a = ( s, \vec{q} )$ be a pre-admissible element such that $s \in T$. Let $q _2 = e ^{r _2}$. We put $\Gamma := r _2 \mathbb Z + \Gamma _0$. A clan associated to $a$ is a maximal subset ${\mathbf c} \subset [ 1, n ]$ with the following property: For each two elements $i, j \in {\mathbf c}$, there exists a sequence $i = i _0, i _1, \ldots, i _m = j$ (in ${\mathbf c}$) such that
$$\{ \log _{i _k} ( s ) \pm \log _{i _{k + 1}} ( s ) \} \cap \{ \pm r _2 + \Gamma _0, \Gamma _0 \} \neq \emptyset \quad \text{for each } 0 \le k < m.$$
We have a disjoint decomposition
$$[1, n] = \bigsqcup _{{\mathbf c} \in \mathcal C _a} {\mathbf c},$$
where each ${\mathbf c}$ is a clan associated to $a$ and $\mathcal C _a$ is the set of clans associated to $a$. For a clan ${\mathbf c}$, we put $n ^{\mathbf c} := \# {\mathbf c}$.
\end{definition}

We assume the setting of Definition \ref{LC} in the rest of this section unless stated otherwise. At the level of Lie algebras, we have a decomposition
$$\mathfrak g ( s ):= \mathfrak t \oplus \bigoplus _{\scriptsize \begin{matrix}i < j, \sigma _1, \sigma _2 \in \{ \pm 1 \},\\ \sigma _1 \log _i ( s ) + \sigma _2 \log _j ( s ) \equiv 0\end{matrix}} \mathfrak g ( s ) [ \sigma _1 \epsilon _i + \sigma _2 \epsilon _j] \oplus \bigoplus _{\scriptsize \begin{matrix}i \in [1, n], \sigma \in \{ \pm 1 \},\\ 2 \log _i ( s ) \equiv 0\end{matrix}} \mathfrak g ( s ) [ \sigma 2 \epsilon _i ],$$
where $\equiv$ means modulo $\Gamma _0$. For each ${\mathbf c} \in \mathcal C _a$, we define a Lie algebra $\mathfrak g ( s ) _{\mathbf c}$ as the Lie subalgebra of $\mathfrak g ( s )$ defined as
$$\bigoplus _{i \in {\mathbf c}} \mathbb C \epsilon _i \oplus \bigoplus _{\scriptsize \begin{matrix}i < j \in {\mathbf c}, \sigma _1, \sigma _2 \in \{ \pm 1 \},\\ \sigma _1 \log _i ( s ) + \sigma _2 \log _j ( s ) \equiv 0\end{matrix}} \mathfrak g ( s ) [ \sigma _1 \epsilon _i + \sigma _2 \epsilon _j] \oplus \bigoplus _{\scriptsize \begin{matrix}i \in {\mathbf c}, \sigma \in \{ \pm 1 \},\\ 2 \log _i ( s ) \equiv 0\end{matrix}} \mathfrak g ( s ) [ \sigma 2 \epsilon _i ],$$
where $\equiv$ means modulo $\Gamma _0$. Moreover, we have
\begin{eqnarray}
\mathfrak g ( s ) = \bigoplus _{{\mathbf c} \in \mathcal C _a} \mathfrak g ( s ) _{\mathbf c}. \label{wdecg}
\end{eqnarray}
In particular, we have $[ \mathfrak g ( s ) _{\mathbf c}, \mathfrak g ( s ) _{{\mathbf c} ^{\prime}} ] = 0$ unless ${\mathbf c} = {\mathbf c} ^{\prime}$. Let $G ( s ) _{\mathbf c}$ be the connected subgroup of $G ( s )$ which has $\mathfrak g ( s ) _{\mathbf c}$ as its Lie algebra.

The following theorem is a consequence of Steinberg's centralizer theorem and the Borel-de Siebenthal theorem, for which we present a proof for the reference purpose.

\begin{theorem}[Centralizer theorem for symplectic groups]\label{Steinberg}
Let $A \subset T$ be an algebraic subgroup. Then, the group $G ( A )$ is connected.
\end{theorem}

\begin{proof}
By a Lie algebra calculation, the group $G ( A ) ^{\circ}$ is generated by $T$ and unipotent one-parameter subgroups $U _{\alpha}$ ($\alpha \in R$) such that $\alpha ( A ) = \{ 1 \}$. By a repeated use of the Borel-de Siebenthal theorem \cite{BS}, the root system of $G ( A ) ^{\circ}$ is the product of the standard presentations of the root systems of
\begin{eqnarray}
\mathop{GL} ( m, \mathbb C ), \mathop{SL} ( 2, \mathbb C ),\text{ or } \mathop{Sp} (2m, \mathbb C ).\label{clBS}
\end{eqnarray}
In particular, the derived group of $G ( A ) ^{\circ}$ must be simply connected.
Now we prove the theorem by induction on the cardinality $k$ of a generating set of $A$. (Notice that the word "generating" means the Zariski closure of the group generated by a given subset of $T$ is $A$. Hence, we can assume the finiteness of the cardinality of a such set.) The case $k = 1$ is an immediate consequence of Steinberg's centralizer theorem (c.f. \cite{C85} 3.5.6). If the assertion is true for smaller $k$, then it suffices to consider the centralizer of a semi-simple element in a group listed at (\ref{clBS}). This is again connected by Steinberg's centralizer theorem. Therefore, the induction proceeds and we obtain the result.
\end{proof}

\begin{lemma}
We have $G ( s ) = \prod _{{\mathbf c} \in \mathcal C _a} G ( s ) _{\mathbf c}$.
\end{lemma}
\begin{proof}
By ($\ref{wdecg}$), it is clear that $\prod _{{\mathbf c} \in \mathcal C _a} G ( s ) _{\mathbf c}$ is equal to the identity component of $G ( s )$. Since $G$ is a simply connected semi-simple group, it follows that $G ( s )$ is connected by Theorem \ref{Steinberg}. In particular, we have $G ( s ) \subset \prod _{{\mathbf c} \in \mathcal C _a} G ( s ) _{\mathbf c}$ as desired.
\end{proof}
We denote $B \cap G ( s ) _{\mathbf c}$ and ${} ^w B \cap G ( s ) _{\mathbf c}$ by $B ( s ) _{\mathbf c}$ and ${} ^w B ( s ) _{\mathbf c}$, respectively.

\begin{convention}
We denote by $\mathbb V ^{a}$ the image of $\mathbb V ^{a} _2$ to $\mathbb V$ via the map
$$\mathbb V _2 \ni ( X _0 \oplus X _1 \oplus X _2 ) \mapsto ( ( X _0 + X _1 ) \oplus X _2 ) \in \mathbb V.$$
Since $q _0 \neq q _1$, we have $\mathbb V ^{a} \cong \mathbb V _2 ^{a}$.
\end{convention}

For each ${\mathbf c} \in \mathcal C _a$, we define
$$\mathbb V ^{a} _{{\mathbf c}} := \sum _{i, j \in {\mathbf c}, \sigma _1, \sigma _2, \sigma _3 \in \{ \pm 1 \}} \mathbb V ^{a} [ \sigma _1 \epsilon _i + \sigma _2 \epsilon _j] \oplus \mathbb V ^{a} [ \sigma _3 \epsilon _i ].$$
It is clear that $\mathbb V ^{a} = \bigoplus _{{\mathbf c} \in \mathcal C _a} \mathbb V ^{a} _{{\mathbf c}}$. By the comparison of weights, the $\mathfrak g ( s ) _{\mathbf c}$-action on $\mathbb V ^{a} _{{\mathbf c} ^{\prime}}$ is trivial unless ${\mathbf c} = {\mathbf c} ^{\prime}$.

\begin{remark}
Since ${\mathbf c}$ is not an integer and we do not use $\mathbb V _{\ell}$ in the rest of this paper (except for \S 7), we use the notation $\mathbb V _{\mathbf c} ^{a}$. The author hopes the reader not to confuse $\mathbb V _{\mathbf c} ^{a}$ with $( \mathbb V _{\ell} )^{a}$.
\end{remark}

\begin{lemma}\label{odecomp}
Let $\mathbb O \subset \mathfrak N ^{a} _+$ be a ${\bm G} ( a )$-orbit. Let $\mathbb O _{{\mathbf c}}$ denote the image of $\mathbb O$ under the natural projection $\mathbb V ^{a} \rightarrow \mathbb V ^{a} _{{\mathbf c}}$. Then, we have a product decomposition $\mathbb O = \oplus _{{\mathbf c} \in \mathcal C _a} \mathbb O _{{\mathbf c}}$.
\end{lemma}

\begin{proof}
Let $X \in \mathbb V ^{a}$. There exists a family $\{ X _{{\mathbf c}} \} _{{\mathbf c} \in \mathcal C _a}$ ($X _{{\mathbf c}} \in \mathbb V ^{a} _{{\mathbf c}}$) such that $X = \sum _{{\mathbf c} \in \mathcal C _a} X _{{\mathbf c}}$. We have $G ( s ) X = \bigoplus _{\mathbf c \in \mathcal C _a} G ( s ) _{\mathbf c} X _{\mathbf c}$. For each of $i = 0, 1$, the clan $\mathbf c \in \mathcal C _a$ such that $( V _1 ^{(s, q _i)} \cap \mathbb V _{\mathbf c} ) \neq \{ 0 \}$ is at most one since clans are determined by the $s$-eigenvalues of $V _1$. Let $\mathbf{c} ^i$ ($i = 1, 2$) be the unique clan such that $( V _1 ^{(s, q _i)} \cap \mathbb V _{\mathbf c ^i} ) \neq \{ 0 \}$. Let $\mathbb G _{\mathbf c}$ be the product of scalar multiplications of $V _1 ^{(s, q _i)}$ such that $V _1 ^{(s, q _i)} \cap \mathbb V ^{a} _{\mathbf c} \neq \{ 0 \}$. Since the set of $a$-fixed points of a conic variety in $\mathbb V$ is conic, we have $( G ( s ) _{\mathbf{c}} \times ( \mathbb C ^{\times} ) ^3 ) X _{\mathbf c} = ( G ( s ) _{\mathbf c} \times \mathbb G _{\mathbf c} ) X _{\mathbf c}$. We have $\prod _{\mathbf c \in \mathcal C _a} ( G ( s ) _{\mathbf c} \times \mathbb G _{\mathbf c} ) \subset {\bm G} ( a )$. It follows that
$${\bm G} ( a ) X = \bigoplus _{{\mathbf c} \in \mathcal C _a} {\bm G} ( a ) X _{{\mathbf c}} = \bigoplus _{{\mathbf c} \in \mathcal C _a} ( G ( s ) _{\mathbf{c}} \times \mathbb G _{\mathbf c} ) X _{{\mathbf c}} = \bigoplus _{{\mathbf c} \in \mathcal C _a} \mathbb O _{{\mathbf c}}$$
as desired.
\end{proof}

For each $w \in W$, we define
$$F ^{a} _+ ( w ) := G ( s ) \times ^{{} ^w B ( s )} ( \dot{w} \mathbb V ^+ \cap \mathbb V ^{a} ).$$
Similarly, we define
$$F ^{a} _+ ( w, {\mathbf c} ) := G ( s ) _{\mathbf c} \times ^{{} ^w B ( s ) _{\mathbf c}} ( \dot{w} \mathbb V ^+ \cap \mathbb V ^{a} _{{\mathbf c}} )$$
for each ${\mathbf c} \in \mathcal C _a$.

\begin{lemma}
We have $F ^{a} _+ = \cup _{w \in W} F ^{a} _+ ( w )$.
\end{lemma}

\begin{proof}
The set of $a$-fixed points of $G / B$ is a disjoint union of flag varieties of $G ( s )$. It follows that each point of $F ^{a} _+$ is $G ( s )$-conjugate to a point in the fiber over a $T$-fixed point of $G / B$.
\end{proof}

The local structures of these connected components are as follows.

\begin{lemma}\label{pdecomp}
For each $w \in W$, we have
$$F ^{a} _+ ( w ) \cong \prod _{{\mathbf c} \in \mathcal C _a} F ^{a} _+ ( w, {\mathbf c} ).$$
\end{lemma}

\begin{proof}
The set $\mathbb V ^{a} _{{\mathbf c}}$ is $T$-stable for each ${\mathbf c} \in \mathcal C _a$. Hence, we have
$$F ^{a} _+ ( w ) = G ( s ) \times ^{{} ^w B ( s )} ( \dot{w} \mathbb V ^+ \cap \mathbb V ^{a}) \cong G ( s ) \times ^{{} ^w B ( s )} ( \bigoplus _{\mathbf c \in \mathcal C _a} ( \dot{w} \mathbb V ^+ \cap \mathbb V ^{a} _{{\mathbf c}} ) ).$$
Since we have $G ( s ) / B ( s ) \cong \prod _{\mathbf c \in \mathcal C _a} G ( s ) _{\mathbf c} / B ( s ) _{\mathbf c}$, we deduce
$$G ( s ) \times ^{{} ^w B ( s )} ( \dot{w} \mathbb V ^+ \cap \mathbb V ^{a} _{\mathbf c} ) \cong \prod _{{\mathbf c} ^{\prime} \in \mathcal C _a} G ( s ) _{{\mathbf c} ^{\prime}} \times ^{{} ^w B ( s ) _{{\mathbf c} ^{\prime}}} ( \dot{w} \mathbb V ^+ \cap \mathbb V ^{a} _{{\mathbf c}} \cap \mathbb V ^{a} _{{\mathbf c} ^{\prime}} ).$$
Here the RHS is isomorphic to
$$F ^{a} _+ ( w, {\mathbf c} ) \times \prod _{{\mathbf c} \neq {\mathbf c} ^{\prime}} G ( s ) _{{\mathbf c} ^{\prime}} / {} ^w B ( s ) _{{\mathbf c} ^{\prime}}.$$
Gathering these information yields the result.
\end{proof}

We define a map ${} ^w \mu _{{\mathbf c}} ^{a}$ by
$${} ^w \mu _{{\mathbf c}} ^{a} : F ^{a} _+ ( w, {\mathbf c} ) = G ( s ) _{{\mathbf c}} \times ^{{} ^w B ( s ) _{{\mathbf c}}} ( \dot{w} \mathbb V ^+ \cap \mathbb V ^{a} _{{\mathbf c}} ) \longrightarrow \mathbb V ^{a} _{{\mathbf c}}.$$

We put $G _{\mathbf c} := \mathop{Sp} ( 2 n ^{\mathbf c} )$ and $s _{\mathbf c} := \exp ( \sum _{i \in \mathbf c} ( \log _i ( s ) ) \epsilon _i ) \in T$. We have embeddings
$$s = \prod _{{\mathbf c} \in \mathcal C _a} s _{\mathbf c} \in \prod _{{\mathbf c} \in \mathcal C _a} \mathop{Sp} ( 2 n ^{\mathbf c} ) \subset \mathop{Sp} ( 2 n ),$$
induced by the following identifications:
\begin{eqnarray}
\mathfrak g ( s ) _{\mathbf c} = \mathfrak g _{\mathbf c} ( s _{\mathbf c} ) \subset \left( \bigoplus _{i \in \mathbf c} \mathbb C \epsilon _i \right) \oplus \bigoplus _{\scriptsize \begin{matrix}\alpha = \sigma _1 \epsilon _i + \sigma _2 \epsilon _j \neq 0\\ \sigma _1, \sigma _2 \in \{ \pm 1 \}, i, j \in \mathbf c\end{matrix}} \mathfrak g [ \alpha ] = \mathfrak g _{\mathbf c}.\label{identg}
\end{eqnarray}

Note that we have $G ( s ) _{\mathbf c} = G _{\mathbf c} ( s _{\mathbf c} ) \subsetneq G _{\mathbf c}$ in general.

Let $\mathbb V ( {\mathbf c} )$ be the $1$-exotic representation of $G _{\mathbf c}$. We have a natural embedding $\mathbb V _{\mathbf c} ^{a} \subset \mathbb V ( {\mathbf c} )$ of $G ( s )_{\mathbf c}$-modules. (The $G ( s )_{\mathbf c}$-module structure on the RHS is given by the restriction of the $G _{\mathbf c}$-action.)

Let $\nu = ( a, X )$ be a pre-admissible parameter. We have a family of pre-admissible parameters $\nu _{\mathbf c} := ( s _{\mathbf c}, \vec{q}, X _{\mathbf c} )$ of $G _{\mathbf c}$'s such that $s = \prod _{\mathbf c} s _{\mathbf c}$, $X = \oplus _{\mathbf c} X _{\mathbf c}$. We denote
$$\nu = \prod _{\mathbf c \in \mathcal C _a} \nu _{\mathbf c}$$
and call it the clan decomposition of $\nu$. Let $W _{a} := \prod _{\mathbf{c} \in \mathcal C _a} N _{G _{\mathbf{c}}} ( T ) / T$. By Lemma \ref{pdecomp}, we conclude that
\begin{eqnarray}
\bigcup _{w \in W _{a}} F ^{a} _+ ( w ) \subset F ^{a} _+\label{exdiag}
\end{eqnarray}
is the product of the $F ^{a} _+$'s obtained by replacing the pair $( G , \nu )$ by $( G _{\mathbf c}, \nu _{\mathbf c} )$ for all $\mathbf c \in \mathcal C _a$.

\begin{proposition}[Clan decomposition of $\mu ^a$]\label{cd}
For each $w \in W$, we have
$$\mu ^{a} _+ \MID _{F ^{a} _+ ( w )} \cong \prod _{{\mathbf c} \in  \mathcal C _a} {} ^w \mu _{{\mathbf c}} ^{a}.$$
In particular, every irreducible direct summand $A$ of $( \mu _+ ^{a} ) _* \mathbb C _{F ^{a} _+}$ is written as an external product of $G ( s )$-equivariant sheaves appearing in $( {} ^w \mu _{{\mathbf c}} ^{a} )_* \mathbb C _{F ^{a} _+ ( w, {\mathbf c} )}$ $($up to degree shift$)$.
\end{proposition}

\begin{proof}
The first assertion follows from the combination of Lemma \ref{odecomp}, Lemma \ref{pdecomp}, and the definition of ${} ^w \mu _{\mathbf c} ^{a}$. We have $\mathbb C _{F ^{a} _+} = \bigoplus _{F ^{a} _+ ( w ) \subset F ^{a} _+} \mathbb C _{F ^{a} _+ ( w )}$. A direct summand of $( \mu _+ ^{a} ) _* \mathbb C _{F ^{a} _+}$ is a direct summand of $( \mu ^{a} _+ ) _* \mathbb C _{F ^{a} _+ ( w )}$ for some $w \in W$. Since
$$( \mu ^{a} _+ ) _* \mathbb C _{F ^{a} _+ ( w )} \cong \boxtimes _{{\mathbf c}} ( {} ^w \mu _{\mathbf c} ^{a} ) _* \mathbb C _{F ^{a} _+ ( w, {\mathbf c} )},$$
the second assertion follows.
\end{proof}

\begin{corollary}\label{cp}
Let $\nu = ( a, X )$ be a pre-admissible parameter. Then, it is regular if and only if $\nu _{\mathbf c}$ is a regular pre-admissible parameter of $G _{\mathbf c}$ for every ${\mathbf c} \in \mathcal C _a$.
\end{corollary}

\begin{proof}
Let $W _{0} := N _{G ( s )} ( T ) / T \subset W$. We have a natural inclusion $W _{0} \subset W _{a}$. Here we have
$$\mu ^{a} _+ = \bigsqcup _{w \in W / W _{0}} \mu ^{a} _+ \MID _{F ^{a} _+ ( w )},$$
where we regard $W / W _{0} \subset W$ by taking some representative. For each $w \in W$, there exists $v \in W _{a}$ such that ${} ^w \mathbb V ^+ \cap \mathbb V ^{a} = {} ^{v} \mathbb V ^+ \cap \mathbb V ^{a} \subset \mathbb V ^{a}$. Moreover, we can choose $v$ so that ${} ^w B ( s ) _{\mathbf c} = {} ^v B ( s ) _{\mathbf c}$ holds for each ${\mathbf c} \in \mathcal C _a$. As a consequence, all $F ^{a} _+ ( w )$ are isomorphic to one of $F ^{a} _+ ( w )$ ($w \in W _{a}$) as ${\bm G} ( a )$-varieties, together with maps $\mu ^{a} _+ \MID _{F ^{a} _+ ( w )}$ to $\mathbb V ^{a}$. Therefore, $\nu$ is regular if and only if an intersection cohomology complex with its support $\overline{{\bm G} ( a ) \nu}$ (with degree shift) appears in $( \mu ^{a} _+ ) _* \mathbb C _{F ^{a} _+ ( w )}$ for some $w \in W _{a}$. Hence, Proposition \ref{cd} implies the result.
\end{proof}

Corollary \ref{cp} reduces the analysis of the decomposition pattern of $( \mu _+ ^{a} ) _* \mathbb C _{F ^{a} _+}$ into the case that $\nu$ has a unique clan.

\section{On stabilizers of exotic nilpotent orbits}\label{sten}

We retain the setting of \S \ref{fHecke}.

\begin{lemma}\label{red2}
Let $H$ be a connected linear algebraic group and let $\mathcal X$ be a variety with $H$-action. Let $H = H _r H _u$ be a Levi decomposition of $H$ with $H _r$ its reductive part. If $\mathsf{Stab} _{H _r} x$ is connected for $x \in \mathcal X$, then so is $\mathsf{Stab} _{H} x$.
\end{lemma}

\begin{proof}
Assume to the contrary to deduce contradiction. Let $h \in \mathsf{Stab} _{H} x$ be an element which is not in the identity component. Let $h = h _r h _u \in H _r H _u$ be its Levi decomposition. For some $k > 1$, we have $h ^k \in ( \mathsf{Stab} _{H} x ) ^{\circ}$. This implies the existence of $g \in ( \mathsf{Stab} _{H} x ) ^{\circ}$ which satisfies $h ^k = g ^k$. Let $g = g _r g _u$ be the Levi decomposition. We have $H _u \triangleleft H$, which claims $h _r ^k = g _r ^k$. Replacing $h$ by $g ^{-1} h$, we further assume $h ^k _r = 1$. Here we have $h ^k \in ( \mathsf{Stab} _{H _u} x ) ^{\circ} = \mathsf{Stab} _{H _u} x$. Put $u := h ^k \in \mathsf{Stab} _{H _u} x$. Let $U$ denote the group given as the Zariski closure of the group generated by $h$. We have its connected component decomposition $U = U _0 \sqcup U _1 \sqcup U _2 \sqcup \cdots$, where $1, u \in U _0$ and $h \in U _1$. Since $h _r$ is of finite order, $U _0$ is unipotent and each of $U _i$ is a homogeneous $U _0$-space. Let $U _0 ^{(m)}$ be the $m$-th lower central subgroup of $U _0$. For each $m$, the adjoint $h$-action preserves $U _0 ^{(m)}$. It follows that if $u \in U _{0} ^{(m)}$ for some $m$, then we have $( h u _m ) ^k \in U _{0} ^{(m)}$ for each $u _m \in U _0 ^{(m)}$. Moreover, we have $h u h ^{-1} \equiv u \!\! \mod U _0 ^{(m + 1)}$ by $u = h ^k$. Since $U _0 ^{(m)} / U _0 ^{(m+1)}$ is abelian, we deduce
$$\{ ( h u _m ) ^k ; u _m \in U _0 ^{(m)} \} /  U _0 ^{(m+1)} = \{ \bar{u}_m \in U _0 ^{(m)} / U _0 ^{(m+1)} ; h \bar{u}_m h ^{-1} = \bar{u}_m \} \subset U _0 ^{(m)} / U _0 ^{(m+1)}$$
for each $m$. The second term contains $1 \!\! \mod U _0 ^{(m+1)}$. We have $U _0 ^{(m)} = \{ 1 \}$ for $m \gg 0$. Hence, we can change $h$ if necessary to assume $h ^k = 1$, which implies that $h$ is semisimple. Therefore, $h$ belongs to $\mathsf{Stab}_{H_r}x$. An element of finite order is always semisimple, hence its unipotent part is trivial. Thus, we have $h_u=1$ if $h^r=1$. Therefore, we have contradiction and the result follows.
\end{proof}

\begin{theorem}[Igusa \cite{Ig} Lemma 8, Springer \cite{Sp}]\label{Ig}
Let $\lambda = ( \lambda _1 \ge \lambda _2 \ge \cdots )$ be a partition of $n$. We regard it as a $0$-marked partition. Then, the reductive part of $\mathsf{Stab} _G \mathbf v _{\lambda}$ is
$$L _{\lambda} := \mathop{Sp} ( 2n _1, \mathbb C ) \times \mathop{Sp} ( 2n _2, \mathbb C ) \times \cdots,$$
where the sequence $( n _1, n _2, \ldots )$ are the number of $\lambda _i$'s which share the same value. Moreover, we have
\begin{eqnarray}
\mathsf{Res} ^G _{L _{\lambda}} V _1 = \bigoplus _{i \ge 1} V ( i ) ^{\oplus \lambda _i},\label{resV}
\end{eqnarray}
where $V ( i )$ is the vector representation of $\mathop{Sp} ( 2 n _i )$ with trivial actions of $\mathop{Sp} ( 2 n _j )$ $(j \neq i)$. \hfill $\Box$
\end{theorem}

\begin{remark}
Igusa's result is not as precise as Theorem \ref{Ig}. But we can deduce from its proof without difficulty. Springer \cite{Sp} contains more precise statement.
\end{remark}

\begin{corollary}
Keep the setting of Theorem \ref{Ig}. Then, we can choose a maximal torus of $L _{\lambda}$ inside $T$.
\end{corollary}

\begin{proof}
Let $\sigma = ( \mathbf J, \emptyset )$ be the $0$-marked partition corresponding to $\lambda$. By Lemma \ref{stabTJ}, we have $\mathbb C ^{\times} \subset \mathsf{Stab} _{T _J} \mathbf v ^J _{\sigma}$ for each $J \in \mathbf J$. It follows that $L _{\lambda} \cap T$ contains a torus of dimension $( \sum _{i \ge 0} n _i )$, which implies the result.
\end{proof}

\begin{proposition}\label{2con}
Let $X \in \mathfrak N _2$. Then, $\mathsf{Stab} _G X$ is connected.
\end{proposition}

\begin{proof}
Let $X = ( X _0 \oplus X _1 ) \oplus \mathbf v _{2, \lambda}$, where $\lambda$ is a partition of $n$ regarded as a $0$-marked partition. It suffices to show that the action of $\mathsf{Stab} _G \mathbf v _{2, \lambda}$ on $( X _0 \oplus X _1 )$ has connected stabilizer. Let $L _{\lambda} U _{\lambda}$ be the Levi decomposition of $\mathsf{Stab} _G \mathbf v _{2, \lambda}$. By Lemma \ref{red2}, it is sufficient to show that the stabilizer of $L _{\lambda}$ on $( X _0 \oplus X _1 )$ is connected. By Theorem \ref{Ig}, it suffices to prove that the $G$-stabilizer of finite set of elements in $V _1$ is connected. By a repeated use of Lemma \ref{red2}, it suffices to prove that the $G$-stabilizer of one element in $V _1$ has $\mathop{Sp} ( 2n - 2 )$ as its (reductive) Levi factor. We denote the element $v \in V _1$ and fix a symplectic form on $V _1$ which is preserved by $G$. Then, it is easy to see that $\mathsf{Stab} _G v$ preserves $\mathbb C v$ and the compliment space $v ^{\perp}$ of $V _1$ with respect to the symplectic form. Thus, its Levi component is given as a subgroup of
$$\mathbb C ^{\times} \times \mathop{Sp} ( 2n - 2 ) = \left( \mathbb C ^{\times} \times \mathop{GL} ( 2n - 2, \mathbb C ) \times \mathbb C ^{\times} \right) \cap \mathop{Sp} ( 2n ) \subset \mathop{GL} ( V _1 ),$$
which fixes $v$. (Here the middle group is the Levi component of $\mathop{GL} ( V _1 )$ which preserves a partial flag $\{ 0 \} \subset \mathbb C v \subset v ^{\perp} \subset V _1$.) Therefore, it is $\mathop{Sp} ( 2n - 2 )$ as desired.
\end{proof}

\begin{remark}
Springer \cite{Sp} contains an explicit description of the $G$-stabilizer of each strict normal form. As is seen easily from the proof of Proposition \ref{2con}, it is not hard to write down the $G$-stabilizer of a point of $\mathfrak N _2$ assuming \cite{Sp}.
\end{remark}

\begin{corollary}[of the proof of Proposition \ref{2con}]\label{rsymp}
For each $X \in \mathfrak N _2$, the reductive part of $\mathsf{Stab} _G X$ is a product of symplectic groups. \hfill $\Box$
\end{corollary}

\begin{theorem}[Refined form of Theorem \ref{EJNF}]\label{rEJNF}
Let $\nu = ( a, X ) = ( s, \vec{q}, X _1 \oplus X _2 ) \in {\bm T} \times \mathfrak N$ be a pre-admissible parameter which is admissible or $a = a _0$. Then, we have a clan decomposition
$$\nu = \prod _{\mathbf c \in \mathcal C _a} \nu ^{\mathbf c} = \prod _{\mathbf c \in \mathcal C _a} ( s _{\mathbf c}, \vec{q}, X _{\mathbf c} )$$
with the following properties:
\begin{itemize}
\item Each $\nu _{\mathbf c}$ is an admissible parameter;
\item There exists $g \in G$ such that:
$$g s g ^{-1} \in T \text{ and each } g X _{\mathbf c} \text{ is a strict normal form.}$$
\end{itemize}
\end{theorem}

\begin{proof}
If $a = a _0$, then we have $\mathfrak N _2 ^a \cong \mathfrak N$. Hence, the result reduces to Proposition \ref{wN1} 1).

Thus, we assume that $a$ is admissible. Since the admissibility condition depends only on the configuration of $\vec{q}$, the clan decomposition preserves admissibility. Hence, it suffices to prove the case that $\mathbf c = [ 1, n ]$ is the unique clan of $\mathcal C _a$. Then, two distinct eigenvalues $t _1, t _2$ of $s$ on $V _1$ satisfies
$$t _1 t _2 \text{ or } t _1 /t _2 = q _2 ^m \text{, where } |m| < n.$$
It follows that at least one of $q _0$ or $q _1$ does not appear as a $s$-eigenvalue of $V _1$ by the admissibility condition. Therefore, we can assume $( s - q _1 ) X _1 = 0$ by swapping the roles of $q _0$ and $q _1$ if necessary.

Let us take $G$-conjugate to assume that $X = \mathbf v _{\sigma}$ for a strict marked partition $\sigma = ( \mathbf J, \vec{\delta} )$. By the description of the $G$-stabilizer of $\mathbf v _{\sigma}$, we deduce that we can choose a maximal torus of $\mathsf{Stab} _G \mathbf v _{\sigma}$ inside of $T$. By Lemma \ref{stabTJ} and the fact $\mathbf v _{\sigma}$ is a strict normal form, we deduce that a (possibly disconnected) maximal torus of $\mathsf{Stab} _{\bm G} \mathbf v _{\sigma}$ is taken inside ${\bm T}$. Therefore, we conclude that $( a, \mathbf v _{\sigma} )$ is a strict normal form after taking conjugate of $a$ by the $\mathsf{Stab} _{\bm G} \mathbf v _{\sigma}$-action (or the $\mathsf{Stab} _{G} \mathbf v _{\sigma}$-action).
\end{proof}

\begin{corollary}\label{V10}
Let $a = ( s, q _0, q _1, q _2 ) \in {\bm T}$ be an admissible element. If $\mathcal C _a$ consists of a unique clan $[1, n]$, then we have either $V _1 ^{(s, q _0)} = \{ 0 \}$ or $V _1 ^{(s, q _1)} = \{ 0 \}$.
\end{corollary}

\begin{proof}
See the second paragraph of the proof of Proposition \ref{rEJNF}.
\end{proof}

\begin{theorem}\label{pi1 desc}
Let $( a, X ) = ( s, \vec{q}, X )$ be a pre-admissible parameter. Then, $\mathsf{Stab} _{G ( s )}  X$ is connected.
\end{theorem}

\begin{proof}
The group $\mathsf{Stab} _{G ( s )}  X$ consists of elements of $\mathsf{Stab} _{G}  X$ which commute with $s$ inside $G$. Moreover, this is equal to $( G \cap ( \mathsf{Stab} _{\mathbf G}  X ) ( a ) )$. Let $\mathbf L$ be the Levi part of $\mathsf{Stab} _{\mathbf G}  X$. By Lemma \ref{red2}, the desired component group is the same as that of $( G \cap \mathbf L ( a ) )$. By Corollary \ref{rsymp}, we deduce that $\mathbf L$ is a product of symplectic groups and a torus $T ^{\perp}$ which injects into $( \mathbb C ^{\times} ) ^3$ via the second projection $\mathbf G = G \times ( \mathbb C ^ {\times} ) ^ 3 \to ( \mathbb C ^ {\times} ) ^ 3$.

We take $G$-conjugation if necessary to assume that $X_2$ is a strict $0$-normal form corresponding to a partition $\lambda$ of $n$ and $s \in T$ by Theorem \ref{rEJNF}. Then, the semi-simple groups contributing $\mathbf L$ are direct factors of the subgroup of the group $L = L _{\lambda}$ borrowed from Theorem \ref{Ig} which fix the both of $X _0$ and $X _1$. The $T$-action on $V _1$ is compatible with the restriction of (\ref{resV}). It follows that we have a sequence of semi-simple elements $s _1, s _2, \ldots,$ in $L$ such that
$$Z _L ( s ) = \{ g \in L ; g s = s g \} = \bigcap _{j \ge 1} L ( s _j ).$$
Let $A$ be the Zariski closure of the group generated by $s _1, \ldots$ in $L$. The condition that an element of $L$ fixes $X _0$ or $X _1$ can be translated into a condition that a collection of vectors $\{ X _0 ^1, X _0 ^2, \ldots,\}$ or $\{ X _1 ^1, X _1 ^2, \ldots,\}$ of $\oplus _i V ( i )$ obtained from $X _0$ or $X _1$ by (\ref{resV}) is fixed, respectively. We put $S \subset \oplus _i V (i)$ to be the $A$-span of all the vectors in $\{X _0 ^1,\ldots\} \cup \{X _1 ^1,\ldots\}$. The condition that an element of $L ( A )$ fixes $X _0$ and $X _1$ is the same as fixing each element of $S$. Here the subgroup $L ^{\prime}$ of $L$ which fixes $S$ is a product of (probably smaller) symplectic groups as in the proof of Proposition \ref{2con}. Moreover, the subgroup of $L ( A )$ which fixes $S$ is isomorphic $L ^{\prime} ( A ^{\prime} )$ for some torus $A ^{\prime} \subset L ^{\prime}$ obtained as the Zariski closure of elements of $L ^{\prime}$ which acts as the same as $s _1, \ldots$ to $S ^{\perp} / S$. (Here $S ^{\perp}$ is the orthogonal complement of $S$ with respect to the $G$-invariant symplectic form on $V _1$.)

Therefore, we deduce that $\mathsf{Stab} _{L ( A )} ( X _0 \oplus X _1 )$ is written as a product of the centralizer of some subgroups of maximal torus in symplectic groups. Since each of such groups are connected by Theorem \ref{Steinberg}, we conclude the result.
\end{proof}

\section{Semisimple elements attached to $G \backslash \mathfrak N _1$}\label{sss}
We keep the setting of the previous section.

Let $\sigma := ( \mathbf J, \vec{\delta} )$ be a strict marked partition. Let $\lambda = ( \lambda _1 \ge \lambda _2 \ge \cdots )$ be the partition of $n$ corresponding to $\mathbf J = \{ J _1, J _2, \ldots \}$.

We fix a sequence of positive real numbers $\gamma _0, \gamma _1, \ldots, \gamma _n > ( n + 1 ) \gamma$ such that
\begin{eqnarray}
\{ 2 \gamma _i, 2 \gamma _j, \gamma _i + \gamma _j, \gamma _i - \gamma _j \} \cap ( \Gamma + \mathbb Z \gamma ) = \emptyset \label{cnfg}
\end{eqnarray}
holds for every pair of distinct numbers $(i,j)$ in $[ 0, n ]$.

\begin{remark}
Our choice of $\{ \gamma _k \} _k$ and $\gamma$ are possible since $\mathbb C$ is an extension of the field $\mathbb Q ( q _2, \sqrt{-1}, \pi )$ with infinite transcendental degree.
\end{remark}

We define a semi-simple element $s _{\sigma} \in T$ as follows:
\begin{itemize}
\item If $\delta _1 \MID _{J _k} \equiv 0$, then we set $\log _j s _{\sigma} = \gamma _k - j \gamma$ for each $j \in J _k$; 
\item If $\delta _1 ( j _0 ) = 1$ for $j _0 \in J _k$, then we set $\log _j s _{\sigma} = \gamma _0 - ( j - j _0 ) \gamma$ for each $j \in J _k$.
\end{itemize}

By the definition of strict marked partitions, the choice of $j _0$ is unique for each $J \in \mathbf J$. Hence, $s _{\sigma}$ is uniquely determined. We put $a _{\sigma} := ( s _{\sigma}, e ^{\gamma _0}, 1, e ^{\gamma} ) \in {\bm T}$.

\begin{lemma}\label{asv}
In the above setting, we have $a _{\sigma} \mathbf v _{\sigma} = \mathbf v _{\sigma}$.
\end{lemma}

\begin{proof}
It suffices to prove $( s _{\sigma}, e ^{\gamma _0}, 1, e ^{\gamma} ) \mathbf v _{\sigma} ^J = \mathbf v _{\sigma} ^J$ for each $J \in \mathbf J$. Let $J = J _k$. Then, $\mathbf v _{2, \sigma} ^{J _k}$ is a sum of $y _{i, i + 1}$ for $i, i + 1 \in J _k$, which has $s _{\sigma}$-eigenvalue
$$e ^{( \gamma _k - i \gamma - ( \gamma _k - ( i + 1 ) \gamma ) )} = e ^{\gamma} \text{ or }e ^{( \gamma _0 - ( i - j _0 ) \gamma - ( \gamma _0 - ( i + 1 - j _0 ) \gamma ) )} = e ^{\gamma},$$
where the latter case occurs only if $\delta ( j _0 ) = 1$ for some $j _0 \in J _k$. Hence, we have $s _{\sigma} \mathbf v _{2, \sigma} ^{J _k} = e ^{\gamma} \mathbf v _{2, \sigma} ^{J _k}$. Moreover, we have $s _{\sigma} x _i = e ^{\gamma _0} x _i$ if $\delta _1 ( i ) = 1$. In particular, we have $s _{\sigma} \mathbf v _{1, \sigma} ^{J _k} = e ^{\gamma _0} \mathbf v _{1, \sigma} ^{J _k}$. These calculations imply the desired result.
\end{proof}

Fix a real number ${\mathsf r} > 0$. We define $D _{\sigma} \in T$ to be
$$\log _i D _{\sigma} = \begin{cases} 0 & (\log _i s _{\sigma} \not\in \gamma _0 + \Gamma )\\ - {\mathsf r} ( \underline{\#} J _k ) & (i \in J _k \ni \exists j _0, \delta _1 ( j _0 ) = 1 ) \end{cases}$$

Consider a parabolic subgroup $P _{\sigma}$ of $G ( s _{\sigma} )$:
$$P _{\sigma} := \{ g \in G ( s _{\sigma} ) ; \lim _{N \to \infty} \mathrm{Ad} ( D _{\sigma} ^N ) g \in G ( s _{\sigma} ) \}.$$
It is well-known that $P _{\sigma}$ is a parabolic subgroup of $G ( s _{\sigma} )$. Let $w _{\sigma}$ be the shortest element of $W$ such that
$$\left< w _{\sigma} R ^+, D _{\sigma} \right> \le 1.$$
It is straight-forward to see
$${} ^{w _{\sigma}} B \cap G ( s _{\sigma} ) \subset P _{\sigma}.$$

\begin{lemma}
For a strict marked partition $\sigma$, we have $\mathbf v _{\sigma} \in \mathbb V ^{a _{\sigma}} \cap {} ^{w _{\sigma}} \mathbb V ^+$.
\end{lemma}

\begin{proof}
By Lemma \ref{asv}, it suffices to prove $\mathbf v _{\sigma} \in {} ^{w _{\sigma}} \mathbb V ^+$. The definition of $w _{\sigma}$ implies that
\begin{enumerate}
\item $x _i \in {} ^{w _{\sigma}} \mathbb V ^+$ if and only if {\bf a)} $D _{\sigma} ( \epsilon _{i} ) < 1$ or {\bf b)} $D _{\sigma} ( \epsilon _{i} ) = 1$ and $i > 0$;
\item $y _{ij} \in {} ^{w _{\sigma}} \mathbb V ^+$ if and only if {\bf a)} $D _{\sigma} ( \epsilon _{i} - \epsilon _{j} ) < 1$ or {\bf b)} $D _{\sigma} ( \epsilon _{i} - \epsilon _{j} ) = 1$ and $\epsilon _i - \epsilon _j \in R ^+$.
\end{enumerate}
Since $\mathbf v _{1, \sigma}$ is sum of $x _i$ with $D _{\sigma}$-eigenvalue $< 1$, we have $\mathbf v _{1, \sigma} \in {} ^{w _{\sigma}} \mathbb V ^+$. The vector $\mathbf v _{2, \sigma}$ have $D _{\sigma}$-eigenvalue $1$. By construction, a strict normal form is contained in $\mathbb V ^+$. Therefore, we conclude $\mathbf v _{2, \sigma} \in {} ^{w _{\sigma}} \mathbb V ^+$, which completes the proof.
\end{proof}

\begin{proposition}\label{prehom}
Let $\sigma = ( \mathbf J, \vec{\delta} )$ be a strict marked partition. Then, we have an inclusion
$$P _{\sigma} \mathbf v _{\sigma} \subset \mathbb V ^{a _{\sigma}} \cap {} ^{w _{\sigma}} \mathbb V ^+,$$
which is dense open.
\end{proposition}

Before giving the proof of Proposition \ref{prehom}, we count the set of weights we concern in its proof:

\begin{lemma}\label{prehomwt}
Keep the setting of Proposition \ref{prehom}. Then, the set $\Psi ( \mathbb V ^{a _{\sigma}} \cap {} ^{w _{\sigma}} \mathbb V ^+ )$ is given by the following list:
\begin{enumerate}
\item $\epsilon _i - \epsilon _{i + 1}$ for each $i, i + 1 \in J \in \mathbf J$;
\item $\epsilon _{i + j _0} - \epsilon _{i + j _1 + 1}$ if the following conditions hold:
\begin{itemize}
\item $i + j _0, j _0 \in J _k$, and $i + j _1 + 1, j _1 \in J _{k^{\prime}}$ for some $k, k^{\prime}$;
\item $\delta _1 ( j _0 ) = 1 = \delta _1 ( j _1 )$, and $\underline{\#} J _k > \underline{\#} J _{k ^{\prime}}$;
\end{itemize}
\item $\epsilon _{j _0}$ for each $j _0 \in J _k$ such that $\delta ( j _0 ) = 1$.
\end{enumerate}
\end{lemma}

\begin{proof}
In this proof, we assume all integer index (which are {\it a priori} not necessarily positive) to be positive. By the choice of the sequence $\{ \gamma _k \} _k$, we have
$$| \left< s _{\sigma},  \epsilon _j + \epsilon _{j ^{\prime}} \right> | \ge e ^{- 2n \gamma} \min \{ e ^{\gamma _k + \gamma _{k ^{\prime}}} ; k, k ^{\prime} \in [ 1, n ] \}  > e ^{\gamma}$$
for each $j, j ^{\prime}$. It follows that weights of the form $\pm ( \epsilon _j + \epsilon _{j ^{\prime}} )$ does not belong to $\Psi ( \mathbb V ^{a _{\sigma}} )$. We examine the assertion $\pm \epsilon _j, \epsilon _j - \epsilon _{j ^{\prime}} \in \Psi ( \mathbb V ^{a _{\sigma}} \cap {} ^{w _{\sigma}} \mathbb V ^+ )$ by the case-by-case analysis. We have three cases:
\item (Case $j, j ^{\prime} \in J _k$) We have $\epsilon _j - \epsilon _j ^{\prime} \in \Psi ( \mathbb V ^{a _{\sigma}} )$ if and only if
$$\left< \epsilon _j - \epsilon _{j ^{\prime}}, s _{\sigma} \right> = e ^{( j ^{\prime} - j ) \gamma} = e ^{\gamma}.$$
This forces $j ^{\prime} - j = 1$. Hence, we put $i := j, j ^{\prime} = i + 1$. We have $\left< \epsilon _i - \epsilon _{i + 1}, D _{\sigma} \right> = 1 \le 1$, which verifies the first part of the assertion.
\item (Case $i \in J _k \neq J _{m} \ni j$) By (\ref{cnfg}), we deduce that
$$\left< \epsilon _j - \epsilon _{j ^{\prime}}, s _{\sigma} \right> \in e ^{\Gamma} = q _2 ^{\mathbb Z}$$
if and only if $j, j ^{\prime} \in J _k$ or $\delta _1 ( J _k ) = \{ 0, 1\} = \delta _1 ( J _m )$ holds. We choose $j _0 \in J _k$ and $j _1 \in J _m$ such that $\delta _1 ( j _0 ) = 1 = \delta _1 ( j _1 )$. We write $j = i + j _0$ and $j ^{\prime} := i ^{\prime} + j _1$. Then, we need
$$\left< \epsilon _{i + j _0} - \epsilon _{i ^{\prime} + j _1}, s _{\sigma} \right> = e ^{( i ^{\prime} - i ) \gamma } = e ^{\gamma}.$$
This happens if and only if $i ^{\prime} = i + 1$. By the definition of $D _{\sigma}$, we have
$$\left< \epsilon _{i + j _0} - \epsilon _{i + j _1 + 1}, D _{\sigma} \right> \le 1$$
if and only if $\underline{\#} J _k \ge \underline{\#} J _{k ^{\prime}}$. Since we assume $( \mathbf J, \vec{\delta} )$ to be a strict marked partition, it follows that $\underline{\#} J _k \neq \underline{\#} J _{k ^{\prime}}$ by \ref{snf} 4). This verifies the second part of the assertion.
\item (Case $j \in J _k$) If $\epsilon _j \in \Psi ( \mathbb V ^{a _{\sigma}} )$ or $- \epsilon _j \in \Psi ( \mathbb V ^{a _{\sigma}} )$, then we have $\left< s, \epsilon _j \right> = e ^{\gamma _0}$ or $e ^{- \gamma _0}$, respectively. By (\ref{cnfg}), this forces $\epsilon _j \in \Psi ( \mathbb V ^{a _{\sigma}} )$ and $\delta _1 ( J _k ) = \{ 0, 1 \}$. Let $j _0 \in J _k$ be such that $\delta _1 ( j _0 ) = 1$. Put $j = i + j _0$ for some $i \in \mathbb Z$. Then, we have $\left< \epsilon _{i + j _0}, s _{\sigma} \right> = e ^{i \gamma + \gamma _0} = e ^{\gamma _0}$ if and only if $i = 0$. Moreover, we have $\left< \epsilon _{j _0}, D _{\sigma} \right> = 1 \le 1$ and $j _0 > 0$, which verifies the final part of the assertion.
\end{proof}

\begin{lemma}\label{dePs}
The group $P _{\sigma}$ satisfies the following conditions:
\begin{enumerate}
\item $P _{\sigma} = T U _{\sigma} \subset B$, where $U _{\sigma}$ is an unipotent subgroup of $G$;
\item The Lie algebra ${\mathfrak u} _{\sigma}$ of $U _{\sigma}$ contains $\mathfrak g [ \alpha ] \subset \mathfrak g$ if and only if $\alpha = \epsilon _j - \epsilon _{j ^{\prime}}$, where $j, j ^{\prime}$ are as follows:
\begin{itemize}
\item $j \in J _k \neq J _{k ^{\prime}} \ni j ^{\prime}$ for some $k , k ^{\prime}$;
\item There exists $j _0 \in J _k$ and $j _1 \in J _{k ^{\prime}}$ such that $\delta _1 ( j _0 ) = 1 = \delta _1 ( j _1 )$;
\item $j - j _0 = j ^{\prime} - j _1$ and $j < j ^{\prime}$.
\end{itemize}
\end{enumerate}
\end{lemma}

\begin{proof}
Let $L _{\sigma}$ be the reductive Levi component of $P _{\sigma}$ which contains $T$. Let $A$ be the Zariski closure of the group generated by $D _{\sigma}$ and $s _{\sigma}$. Then, $G ( A )$ is connected by Theorem \ref{Steinberg}. Thus, we have $L _{\sigma} = T$ if we have $\alpha ( A ) \neq \{ 1 \}$ for each $\alpha \in R$. This is equivalent to $\alpha ( D _{\sigma} ) \neq 1$ or $\alpha ( s _{\sigma} ) \neq 1$ holds for each $\alpha \in R$ since $R$ is a finite set. By (\ref{cnfg}), we have $\left< \epsilon _i - \epsilon _j, s _{\sigma} \right> = 1$ only if
\begin{eqnarray}
\exists \kappa \in \{ \pm 1 \} \text{ s.t. }\kappa i \in J _k \neq J _{k ^{\prime}} \ni \kappa j \text{ and } \delta _1 ( J _k ) = \{ 0, 1 \} = \delta _1 ( J _{k ^{\prime}} ),\label{cdM}
\end{eqnarray}
where $J _k, J _{k ^{\prime}} \in \mathbf J$. By \ref{snf} 4), we deduce that
$$\left< \epsilon _i, D _{\sigma} \right> \neq \left< \epsilon _j, D _{\sigma} \right> \text{ for each }i \in J _k, j \in J _{k ^{\prime}}.$$
Therefore, we conclude $L _{\sigma} = T$.

Since $T$ normalizes the unipotent part of $P _{\sigma}$, we describe all one-parameter unipotent subgroup of $G$ belonging to $P _{\sigma}$ in order to prove the assertion. This is equivalent to count the set of weight spaces $\mathfrak g [ \alpha ] \subset \mathfrak g$ which is fixed by $s _{\sigma}$ and has eigenvalue $\le 1$ with respect to $D _{\sigma}$. We examine the case $\alpha = \epsilon _i - \epsilon _j$ with the assumption (\ref{cdM}) for $\kappa = + 1$. (This last part of the assumption is achieved by swapping the roles of $i$ and $j$ if necessary.) Fix $j _0 \in J _k$ and $j _1 \in J _{k ^{\prime}}$ such that $\delta  _1 ( j _0 ) = 1 = \delta _1 ( j _1 )$. Then, the definition of $s _{\sigma}$ further asserts $i - j _0 = j - j _1$. In order that $D _{\sigma}$ has eigenvalue $\le 1$, we need to have
$$\left< \epsilon _i, D _{\sigma} \right> \le \left< \epsilon _j, D _{\sigma} \right>,$$
which is equivalent to $\underline{\#} J _k \ge \underline{\#} J _{k ^{\prime}}$. This implies $\# J _k >\# J _{k ^{\prime}}$ by \ref{snf} 4). It follows that $i < j$, which verifies the second condition. Since $\alpha = \epsilon _i - \epsilon _j \in R ^+$ in this case, we also deduce the first condition. 
\end{proof}

\begin{proof}[Proof of Proposition \ref{prehom}]
Since $P _{\sigma} \subset G ( s _{\sigma} )$, we have $P _{\sigma} \mathbf v _{\sigma} \subset \mathbb V ^{a _{\sigma}}$. Since the reductive part of $P _{\sigma}$ is equal to $T$, we deduce $P _{\sigma} \mathbf v _{\sigma} \subset {} ^{w _{\sigma}} \mathbb V ^+$. Therefore, it suffices to prove the following equality at the level of tangent space
\begin{eqnarray}
T _{\mathbf v _{\sigma}} ( P _{\sigma} \mathbf v _{\sigma} ) \cong \mathfrak p _{\sigma} \mathbf v _{\sigma} = \mathbb V ^{a _{\sigma}} \cap {} ^{w _{\sigma}} \mathbb V ^+ \label{taneq}
\end{eqnarray}
in order to deduce the assertion. Consider a $T$-weight decomposition $\mathbf v _{\sigma} ^J = \sum _{\beta \in \Xi _J} v _{\beta}$, where $J \in \mathbf J$ and $0 \neq v _{\beta} \in \mathbb V [ \beta ]$. Each $\Xi _J$ consists of linearly independent weights of $X ^* ( T _J )$. Moreover, we have $X ^* ( T _J ) \cap X ^* ( T _{J ^{\prime}} ) = \{ 0 \}$ in $X ^* ( T )$ (by using the natural embeddings). Hence, we deduce
$$\mathfrak t \mathbf v _{\sigma} = \sum _{k \ge 1} \mathfrak t \mathbf v _{\sigma} ^{J _k} = \sum _{k \ge 1} \sum _{\beta \in \Xi _{J _k}} \mathbb C v _{\beta}.$$
It is easy to see that $\bigcup _{k \ge 1} \Xi _{J _k}$ is precisely the set of $T$-weights described in Lemma \ref{prehomwt} 1) and 3).

In the below, we apply the action of $\mathfrak u _{\sigma}$ (c.f. Lemma \ref{dePs}) to fill out each $\mathbb V [ \beta ]$ for each $T$-weight $\beta$ described in Lemma \ref{prehomwt} 2). Such a $\beta$ is written as $\epsilon _i - \epsilon _j$, where $i \in J _k, j \in J _{k ^{\prime}}$ are as in Lemma \ref{prehomwt} 2). By explicit calculation, we have a non-zero element of $\mathfrak g$ of weight $\epsilon _{m + j _0} - \epsilon _{m + j _1}$ which satisfies
$$\xi _m \mathbf v _{\sigma} = \begin{cases}
y _{m - 1 + j _0, m + j _1} - y _{m + j _0, m + j _1 + 1} & (m + j _0 + 1 \in J _{k ^{\prime}})\\
y _{m - 1 + j _0, m + j _1} & (m + j _0 + 1 \not\in J _{k ^{\prime}})
\end{cases}.$$
for each $m + j _0 \in J _{k ^{\prime}}$. (Here we implicitly used $\overline{\#} J _k > \overline{\#} J _{k ^{\prime}}$, which is deduced from $\underline{\#} J _k > \underline{\#} J _{k ^{\prime}}$ by \ref{snf} 4).) We know $\xi _m \in \mathfrak p _{\sigma}$ by Lemma \ref{dePs} 2). We have
$$\left( \sum _{m \in \mathbb Z; m + j _1 \in J _{k ^{\prime}}} \mathbb C \xi _m \right) \mathbf v _{\sigma} = \sum _{m \in \mathbb Z; m + j _1 \in J _{k ^{\prime}}} \mathbb V [ \epsilon _{m - 1 + j _0} - \epsilon _{m + j _1} ].$$
By summing up for all possible pairs $(J _k, J _{k ^{\prime}} ) \in \mathbf J$, the set of $T$-weights appearing in the RHS exhausts the $T$-weights described in Lemma \ref{prehomwt} 2).
\end{proof}

\begin{corollary}\label{prehomc}
Keep the setting of Proposition \ref{prehom}. Let $\nu = ( a, \mathbf v _{\sigma} ) = ( s, \vec{q}, \mathbf v _{\sigma} )$ be an admissible parameter. Then, the natural embedding
$$P _{\sigma} ( s ) \mathbf v _{\sigma} \subset \mathbb V ^a \cap \mathbb V ^{a _{\sigma}} \cap {} ^{w _{\sigma}} \mathbb V ^+$$
is dense open.
\end{corollary}

\begin{proof}
The assertion follows merely by taking $a$-fixed part of (\ref{taneq}) in the proof of Proposition \ref{prehom}.
\end{proof}

\section{A vanishing theorem}\label{avt}
We retain the setting of the previous section.

\begin{definition}[Exotic Springer fibers]
For each $X \in \mathfrak N _2$, we define $\mathcal E _X$ as the image of the projection of
$${\bm \mu} ^{-1} ( X ) \subset {\bm F} = G \times ^B \mathbb V ^+$$
to $G / B$. For a pre-admissible parameter $\nu = ( a, X )$, we have a subvariety ${\bm \mu} ^{-1} ( X ) ^a \subset {\bm \mu} ^{-1} ( X )$. We denote the image of ${\bm \mu} ^{-1} ( X ) ^a$ under the projection to $G/B$ by $\mathcal E _X ^a$. By construction, we have $\mathcal E _X \cong {\bm \mu} ^{-1} ( X )$ and $\mathcal E _X ^a \cong {\bm \mu} ^{-1} ( X ) ^a \subset {\bm F} ^a$. We call $\mathcal E _X$ and $\mathcal E _X ^a$ exotic Springer fibers.
\end{definition}

\begin{theorem}[Homology vanishing theorem]\label{ov}
Let $\nu = ( a, X )$ be a pre-admissible parameter which is admissible or $a = a _0$. Then, we have
$$H _{2i + 1} ( \mathcal E _X ^a ) = 0 \text{ for every } i = 0, 1, \ldots.$$
Moreover, we have an isomorphism
$$\mathsf{ch} : \mathbb C \otimes _{\mathbb Z} K ( \mathcal E _X ^a ) \stackrel{\cong}{\longrightarrow} H _{\bullet} ( \mathcal E _X ^a ).$$
\end{theorem}

\begin{remark}{\bf 1)} The map $\mathsf{ch}$ in Theorem \ref{ov} is the homology Chern character map. (See e.g. \cite{CG} \S 5.8.) It sends the class of the (embedded) structure sheaf $\mathcal O _C$ for a closed subvariety $C \subset \mathcal E _X ^a$ to
$$\mathsf{ch} [ \mathcal O _C ] = [ C ] + \text{lower degree terms} \in H _{2 \dim C} ( \mathcal E _X ^a ) \oplus \cdots \oplus H _0 ( \mathcal E _X ^a ).$$
{\bf 2)} The first part of Theorem \ref{ov} is valid even for integral coefficient case when $G ( s ) \subset \mathop{GL} ( n, \mathbb C )$ (c.f. \cite{AH})\footnote{Previous versions of this paper also contain a similar result (since math.RT/0601155v3, April/2006). The author decided to drop it since it is unnecessary to prove our main theorems and \cite{AH} contains a better proof.}. Here we present a proof along the line of earlier versions of this paper, with an enhancement (the proof of Theorem \ref{ov} modulo Proposition \ref{wver} given in \S 6.2) informed to the author by Eric Vasserot.
\end{remark}

\subsection{Review of general theory on homology vanishing}
In this subsection, we recall several definitions and results of \cite{BH} and \cite{DLP} which we need in the course of our proof of Theorem \ref{ov}.

\begin{definition}[$\alpha$-partitions]
A partition of a variety $\mathcal X$ over $\mathbb C$ is said to be an $\alpha$-partition if it is indexed as $\mathcal X _1, \mathcal X _2, \ldots \mathcal X _k$ in such a way that $\mathcal X _1 \cup \ldots \cup \mathcal X _i$ is closed for every $i = 1, \ldots, k$.
\end{definition}

\begin{theorem}[\cite{DLP} 1.7--1.10]\label{zeroH}
Let $\mathcal X$ be a variety with $\alpha$-partition $\mathcal X _1, \mathcal X _2, \ldots, \mathcal X _k$. If we have
$$H _{2 i + 1} ( \mathcal X _m ) = 0 \text{ for every }i = 0, 1, \ldots$$
for each $m = 1, \ldots, k$, then we have
$$H _{2 i + 1} ( \mathcal X ) = 0 \text{ for every }i = 0, 1, \ldots.$$
Moreover, we have
$$\sum _{i \ge 0} \dim H _{2 i} ( \mathcal X ) = \sum _{m \ge 1} \sum _{i \ge 0} \dim H _{2 i} ( \mathcal X _m ).$$
\end{theorem}



\begin{theorem}[\cite{BH} 9.1]\label{subB}
Let $\mathcal Z$ be a smooth variety with $\mathbb G _m$-action. Assume that for some $t \in \mathbb G _m$, we have
\begin{itemize}
\item $\mathcal Z ^{\mathbb G _m} = \mathcal Z ^t \text{ and } \lim _{N \to \infty} t ^N z \in \mathcal Z ^t \hskip 3mm \forall z \in \mathcal Z$;
\item For each $z _0 \in \mathcal Z ^t$, the set $\{ z \in \mathcal Z ; \lim _{N \to \infty} t ^N z = z _0\}$ defines an affine closed subscheme of $\mathcal Z$.
\end{itemize}
Then, $\mathcal Z$ is a vector bundle over $\mathcal Z ^t$. In particular, the two conditions
$$H _{2 i + 1} ( \mathcal Z ) = 0  \text{ for every }i = 0, 1, \ldots, \text{ and } H _{2 i + 1} ( \mathcal Z ^t ) = 0 \text{ for every }i = 0, 1, \ldots$$
are equivalent. Moreover, we have
$$\sum _{i \ge 0} \dim H _{2i} ( \mathcal Z ) = \sum _{i \ge 0} \dim H _{2i} ( \mathcal Z ^t )$$
if one of the above equivalent conditions hold.
\end{theorem}

\subsection{Proof of vanishing theorem}

This subsection is devoted to the proof of Theorem \ref{ov}.

By taking $G$-conjugation if necessary, we assume $a \in {\bm T}$. We have
$$( a, X ) = ( s, \vec{q}, X ) = \prod _{\mathbf c \in \mathcal C _a} ( s _{\mathbf c}, \vec{q}, X _{\mathbf c} ).$$
By the same argument as in the proof of Corollary \ref{cp}, each connected component of ${\bm \mu} ^{-1} ( X ) ^a$ is a product of connected components of
$$\mathcal E _{X _{\mathbf c}} ^{(s _{\mathbf c}, \vec{q})} \subset \mathop{Sp} ( 2n ^{\mathbf c} ) / ( B \cap \mathop{Sp} ( 2n ^{\mathbf c} ) ) \text{ for all } \mathbf c \in \mathcal C _a.$$
Therefore, by the K\"unneth formula, it suffices to prove the assertion when $\mathcal C _a$ consists of a unique clan $[1, n]$. By Proposition \ref{rEJNF}, we further assume that $s \in T$, and $X = \mathbf v _{\sigma}$ for a strict marked partition $\sigma = ( \mathbf J, \vec{\delta} )$ by taking $G$-conjugate if necessary.

\begin{proposition}[Weak version of Theorem \ref{ov}]\label{wver}
Let $\nu = ( a, X ) = ( s, \vec{q}, \mathbf v _{\sigma} )$ be a pre-admissible parameter which is admissible or $a = a _0$. Assume that $s \in T$ and $\sigma$ is a strict marked partition. For $s _{\sigma} \in T$ defined in the above of Lemma \ref{asv}:
\begin{itemize}
\item We have
$$H _{2i + 1} ( ( \mathcal E _{X} ^a ) ^{s _{\sigma}} ) = 0 \text{ for every } i = 0, 1, \ldots;$$
\item Each connected component of $( \mathcal E _{X} ^a ) ^{s _{\sigma}}$ is smooth projective;
\item We have an isomorphism
$$\mathsf{ch} : \mathbb C \otimes _{\mathbb Z} K ( ( \mathcal E _X ^a ) ^{s _{\sigma}} ) \stackrel{\cong}{\longrightarrow} H _{\bullet} ( ( \mathcal E _X ^a ) ^{s _{\sigma}} ).$$
\end{itemize}
\end{proposition}

Before giving the proof of Proposition \ref{wver}, we complete our proof of Theorem \ref{ov} for $( a, \mathbf v _{\sigma} )$ assuming Proposition \ref{wver} for $( a, \mathbf v _{\sigma} )$.

\begin{proof}[Proof of Theorem \ref{ov} for $( a, X ) = ( a, \mathbf v _{\sigma} )$]
Let $\mathcal E _1, \mathcal E _2, \ldots$ be a sequence of all connected components of $( \mathcal E _{X} ^a ) ^{s _{\sigma}}$. For each $\mathcal E _k$, we set
$$\mathcal B _k := \{ g B \in \mathcal E _X ^a ; \lim _{N \to \infty} s _{\sigma} ^{-N} g B \in \mathcal E _k \}.$$
Let $\pi _k : \mathcal B _k \to \mathcal E _k$ be the $s _{\sigma} ^{-1}$-attracting map. Let
$$P := \{ g \in G ; \lim _{N \to \infty} \mathrm{Ad} ( s _{\sigma} ^{-N} ) g \in G \}$$
be a parabolic subgroup of $G$. It is straight-forward to see
$$\lim _{N \to \infty} \mathrm{Ad} ( s _{\sigma} ^{-N} ) g \in G ( s _{\sigma} )$$
for $g \in P$. It follows that each $\mathcal B _k$ intersects with a unique $P$-orbit in $G / B$. In particular, we can assume that the sequence $\mathcal B _1, \mathcal B _2, \ldots$ forms an $\alpha$-partition of $\bigcup _{k \ge 1} \mathcal B _k \subset \mathcal E _X ^a$ by rearranging the sequence if necessary.

We choose $\{ \gamma _i \} _{i}$ (in the definition of $s _{\sigma}$) so that we have
\begin{align}\nonumber
& \min \{\left< s _{\sigma}, \epsilon _i \right> ; i \in J _k \} < \min \{\left< s _{\sigma}, \epsilon _i \right> ; i \in J _{k ^{\prime}} \} \text{ , and}\\
& \gamma \max \{\left< s _{\sigma}, \epsilon _i \right> ; i \in J _k \} > \max \{\left< s _{\sigma}, \epsilon _i \right> ; i \in J _{k ^{\prime}} \}\label{arrG}
\end{align}
for each $k < k ^{\prime}$. (This choice is possible by $\# J _k \ge \# J _{k ^{\prime}}$ and Definition \ref{snf}.)

\begin{claim}\label{affS}
Each fiber of $\pi _k$ is an irreducible affine scheme.
\end{claim}

\begin{proof}
Let $P = G ( s _\sigma ) U$ be the Levi decomposition. Let $g B \in \mathcal E _k$. We set $F ( gB ) := \{ u \in U ( s ) ; X \in u g \mathbb V ^+ \}$. This is a closed subset of $U ( s )$. Set $U ^{\flat} := ( U ( s ) \cap g U g^{-1} )$. We have a free right $U ^{\flat}$-action on $F ( g B )$. We have $\pi _k ^{-1} ( g B ) \cong F ( gB ) / U ^{\flat}$. Since $\pi _k ^{-1} ( g B )$ is a closed subspace of an affine space $U ( s ) / U ^{\flat}$, it suffices to prove that $F ( gB )$ is an affine space. We have a product decomposition $U ( s ) = U _2 U _1$, where $U _1$ is the product of $U _{\epsilon _i - \epsilon _j} \subset U ( s )$ $(i, j > 0)$ and $U _2$ is the product of $U _{\epsilon _i + \epsilon _j} \subset U ( s )$ ($i,j>0$). By (\ref{arrG}), the space $( U _1 X - X )$ is a linear subspace of $\mathbb V ^+$. Hence, $( ( U _1 X - X ) \cap g \mathbb V ^+ )$ is an affine space. Here $\mathsf{Stab} _{U} X$ is a unipotent group, which is automatically an affine space.

Since $U _2$ acts $V _1 ^+$ trivially, we have $u \in F ( gB )$ only if $u \in U _2 u _1$ with $u _1 \in F ( gB )$. The closed subset $( U _2 u _1 \cap F ( gB ) ) \subset U _2 u _1$ define linear conditions on $U _2$ since $U _2$ is commutative. Let $A \subset T$ denote the Zariski closure of the group generated by $s$ and $s _\sigma$. The group $G ( A )$ normalizes $U _1$ and $U _2$, and $U _1$ normalizes $U _2$. Hence, the condition along different points of $( U _1 \mathcal E _k \cap \mathcal B _k )$ are isomorphic via conjugation of $U _1 G ( A )$. It follows that $F ( gB )$ is a successive fibration of affine spaces by affine spaces, which is itself an affine space.
\end{proof}

\begin{claim}\label{affBE}
The variety $\mathcal B _k$ is a smooth affine bundle over $\mathcal E _k$.
\end{claim}

\begin{proof}
We keep the setting of the proof of Claim \ref{affS}. Let $\mathfrak f ( gB ) := \{ \xi \in \mathrm{Lie} U ( a ) ; \xi X \in g \mathbb V ^+ \}$. Since $g B$ is the unique $s _{\sigma}$-fixed point of $F ( gB )$, we deduce
\begin{eqnarray}
\dim F ( g B ) \le \dim \mathfrak f ( g B ). \label{dimestC}
\end{eqnarray}
Notice that $U$ is invariant under the $G ( s _{\sigma} )$-action. Thus, $\dim \mathfrak f ( gB )$ is invariant along $\mathcal E _k$. In view of Claim \ref{affS}, the assertion follows if the equality of (\ref{dimestC}) holds for each $g B \in \mathcal E _k$. If $\xi \in \mathfrak f ( gB )$ is a $s _{\sigma}$-eigenvector, then we have $\exp ( - \xi ) = 1 - \xi \in F ( gB )$. In particular, we have $\dim F ( gB ) \ge \dim \mathfrak f ( gB )$ since we have enough number of linearly independent tangent lines.
\end{proof}

We return to the proof of Theorem \ref{ov}

By Claim \ref{affBE} and Theorem \ref{subB}, the first assertion reduces to
$$H _{2 i + 1} ( \mathcal E _k ) = 0  \text{ for every }i = 0, 1, \ldots$$
for each $k$. Hence, Proposition \ref{wver} for $(a, X)$ yields the first assertion.

Similarly, Proposition \ref{wver} for $(a,X)$ and the Thom isomorphism give
$$\mathsf{ch} : \mathbb C \otimes _{\mathbb Z} K ( \mathcal B _k ) \stackrel{\cong}{\longrightarrow} H _{\bullet} ( \mathcal B _k )$$
for each $k$. The Chern character map commutes with localization sequences. (c.f. \cite{CG} \S 5.8.) Therefore, a successive application of localization sequences implies the second assertion.
\end{proof}

The rest of this section is devoted to the proof of Proposition \ref{wver}.

We set $( s _\sigma, \vec{q} _\sigma ) := a _{\sigma} \in \bm T$ defined in \S \ref{sss}. We have $a _{\sigma} X = X$. Hence, $a _{\sigma}$ acts on $\mu ^{-1} ( X ) ^a$. Its projection gives the $s _\sigma$-action on $\mathcal E _X ^a$. Let $A$ be the Zariski closure of the subgroup of $\bm T$ generated by $a$ and $a _{\sigma}$. We put $W _A := \{ w \in W ; B ( A ) \subset {} ^{w} B\}$. We put $F ^A ( w ) := G ( A ) \times ^{B ( A )} ( \mathbb V ^{A} \cap {} ^w \mathbb V ^+ )$ for each $w \in W_A$. We have $\bigcup _{w \in W _A} F ^A ( w ) = ( G \times ^B \mathbb V ^+ ) ^A$. Consider the map
$${} ^w \mu ^A : F ^A ( w ) = G ( A ) \times ^{B ( A )} ( \mathbb V ^{A} \cap {} ^{w} \mathbb V ^+ ) \longrightarrow \mathbb V ^{A},$$
for each $w \in W_A$.

\begin{lemma}[Part of Proposition \ref{wver}]\label{spp}
Each connected component of $( \mathcal E _X ^a ) ^{a _{\sigma}}$ is smooth projective.
\end{lemma}

\begin{proof}
Projectivity follows from that of $\mathcal E _X$, which itself follows by Theorem \ref{fgeom} 3). By Lemma \ref{dePs} 1) and Corollary \ref{prehomc}, we deduce that
$$\overline{B ( A ) \mathbf v _{\sigma}} \subset \mathbb V ^{A}$$
is a linear subspace. It follows that $( {} ^w \mu ^A ) ^{-1} ( \overline{B ( A ) \mathbf v _{\sigma}} )$ is a smooth subvariety of $G ( A ) \times ^{B ( A )} ( \mathbb V ^{A} \cap {} ^w \mathbb V ^+ )$. Hence, $( {} ^w \mu ^A ) ^{-1} ( B ( A ) \mathbf v _{\sigma} )$ is a smooth subvariety of $F ^A ( w )$. Since changing $\mathbf v _{\sigma}$ by $B ( A )$-action gives isomorphic fibers, we deduce that $( {} ^w \mu ^A ) ^{-1} ( \mathbf v _{\sigma} )$ is a smooth subvariety of $F ^A ( w )$ as required.
\end{proof}

\begin{corollary}[of the Proof of Lamma \ref{spp}]
The variety $( {} ^w \mu ^A ) ^{-1} ( \overline{B ( A ) \mathbf v _{\sigma}} )$ is smooth. \hfill $\Box$
\end{corollary}

We return to the proof of Proposition \ref{wver}.

We prove the rest of assertions by the induction on the cardinality $n ( \sigma )$ of the set
$$\mathtt N ( \sigma) := \{J \in \mathbf J; \delta _1 ( J ) = \{ 0, 1 \} \}.$$
In other words, we assume Theorem \ref{ov} for every admissible parameter of the form $( a, \mathbf{v} _{\sigma ^{\prime}} )$ such that $n ( \sigma ^{\prime} ) < n ( \sigma )$.
If $n ( \sigma ) = 0$, then Lemma \ref{dePs} 2) asserts that $G ( s _{\sigma} ) = T$. This implies that $( \mathcal E _X ^a ) ^{s _{\sigma}}$ is a union of points. Thus, we obtain the assertion for $n ( \sigma ) = 0$.

We prove the assertion for $n ( \sigma ) = k$ by assuming that the assertion holds for all $n ( \sigma ) < k$. Let $J \in \mathtt N ( \sigma )$ be the member such that $\# J \ge \# J ^{\prime}$ for every $J ^{\prime} \in \mathtt N ( \sigma )$. Let $\sigma ^{\prime}$ be a strict marked partition obtained from $\sigma$ by replacing $\delta _1$ by $\delta _1 ^{\prime}$ defined as:
$$\delta _1 ^{\prime} ( J ) = \{ 0 \}, \text{ and } \delta _1 ^{\prime} ( j ) = \delta _1 ( j ) \text{ for all } j \in [1, n] \backslash J.$$
Let $j _0 \in J$ be the unique element such that $\delta _1 ( j _0 ) = 1$. By Lemma \ref{stabTJ}, there exists $t \in T _J$ such that
$$\lim _{N \to \infty} t ^N \mathbf v _{\sigma} = \mathbf v _{\sigma} - \mathbf v _{1, \sigma} ^J = \mathbf v _{\sigma ^{\prime}}.$$

By Lemma \ref{dePs} 2), every $T$-weight of $P _{\sigma}$ containing $i \in J$ is of the form $\epsilon _i - \epsilon _j$ for some $j \in J ^{\prime}$. Moreover, we have $P _{\sigma ^{\prime}} \subset P _{\sigma}$. It follows that the action of $t \in T$ contracts $P _{\sigma}$ to $P _{\sigma ^{\prime}}$. By Corollary \ref{prehomc}, the $t$-action also contracts $\mathbb V ^A \cap {} ^{w _{\sigma}} \mathbb V ^+$ to
$$\mathbb V ^A \cap \mathbb V ^{a _{\sigma ^{\prime}}} \cap {} ^{w _{\sigma ^{\prime}}} \mathbb V ^+ = \mathbb V ^{A ^{\prime}} \cap {} ^{w _{\sigma ^{\prime}}} \mathbb V ^+,$$
where $A ^{\prime}$ is the Zariski closure of $\left< a, a _{\sigma ^{\prime}} \right> \subset {\bm T}$.

Therefore, the $t$-action contracts $( {} ^w \mu ^A ) ^{-1} ( \overline{B ( A ) \mathbf v _{\sigma}} )$ to $( {} ^w \mu ^{A ^{\prime}} ) ^{-1} ( \overline{B ( A ^{\prime} ) \mathbf v _{\sigma ^{\prime}}} )$. By taking the quotient of
$$\mathcal S := B ( A ) \mathbf v _{\sigma} \cup B ( A ) \mathbf v _{\sigma ^{\prime}}$$
by $\mathsf{Stab} _{B ( A )} \mathbf v _{1, \sigma} ^J$, we obtain an affine plane $\mathbb A ^1$ with contracting $t$-action to the origin. Therefore, we obtain a smooth family of smooth projective varieties over $\mathbb A ^1$ whose fiber over $0 \in \mathbb A ^1$ is $\mathcal E _{\mathbf v _{\sigma ^{\prime}}} ^{A}$ and whose general fiber $\mathcal E _{\mathbf v _{\sigma}} ^{A}$ contracting to $\mathcal E _{\mathbf v _{\sigma ^{\prime}}} ^{A ^{\prime}}$. Moving smooth projective varieties is the same as moving all cycles by rational equivalence. Therefore, it suffices to prove
\begin{eqnarray}
\mathbb C \otimes _{\mathbb Z} K ( \mathcal E _{\mathbf v _{\sigma ^{\prime}}} ^{A} ) \stackrel{\cong}{\longrightarrow} H _{\bullet} ( \mathcal E _{\mathbf v _{\sigma ^{\prime}}} ^{A} ).\label{AA}
\end{eqnarray}
Since $\mathcal E _{\mathbf v _{\sigma ^{\prime}}} ^{A}$ is smooth, the Bia\l ynicki-Birula theorem asserts that $\mathcal E _{\mathbf v _{\sigma ^{\prime}}} ^{A}$ is a union of vector bundles over connected components of $\mathcal E _{\mathbf v _{\sigma ^{\prime}}} ^{A ^{\prime}}$. Since the Chern character map commutes with pullbacks and localization sequences, we deduce (\ref{AA}) from Theorem \ref{subB} and Proposition \ref{wver} for $( a, \mathbf v _{\sigma ^{\prime}} )$. Therefore, we have Proposition \ref{wver} for every admissible parameter of the form $( a, \mathbf v _{\sigma} )$ with $n ( \sigma ) = k$. Hence, the induction proceeds and we have proved Proposition \ref{wver} (and hence Theorem \ref{ov}).

\section{Standard modules and an induction theorem}\label{smi}
We retain the setting of the previous section.

\begin{definition}[Standard modules]\label{stdmod}
Let $\nu = ( a, X )$ be a pre-admissible parameter. We define
$$M _{\nu} := H _{\bullet} ( \mathcal E _X ^a ) \text{ and } M ^{\nu} := H ^{\bullet} ( \mathcal E _X ^a ).$$
By the Ginzburg theory \cite{CG} 8.6, each of $M _{\nu}$ or $M ^{\nu}$ is a $\mathbb H$-module.
\end{definition}

By the symmetry of the construction of varieties involved in $M _{( a, X )}$ and $\mathbb H _a$, we deduce $M _{(a, X)} \cong M _{( \mathrm{Ad} ( g ) a, g X )}$ as $\mathbb H _a = \mathbb H _{\mathrm{Ad} ( g ) a}$-modules for each $g \in {\bm G}$.

Let $s _Q \in T ( \mathbb R )$ be an element such that
\begin{eqnarray}
0 < \left< \alpha, s _Q \right> \le 1 \text{ for all } \alpha \in R ^+.\label{dsQ}
\end{eqnarray}
Let $Q := G ( s _Q )$ and ${\bm Q} := Q \times ( \mathbb C ^{\times} ) ^3$. These are subgroups of $G$ and ${\bm G}$, respectively. We put $\mathbb V _Q := \mathbb V ^{s _Q} _2$ and $\mathfrak N _Q = \mathfrak N _2 ^{s _Q} \subset \mathbb V _Q$. We put $F _Q := Q \times ^{( Q \cap B )} ( \mathbb V _Q \cap \mathbb V ^+ _2 )$. We have a map
$$\mu _Q : F _Q = Q \times ^{( Q \cap B )} ( \mathbb V _Q \cap \mathbb V ^+ _2 ) \longrightarrow \mathfrak N _Q.$$
We define $Z _Q := F _Q \times _{\mathfrak N _Q} F _Q$.

The natural inclusion map
$$F _Q = Q \times ^{( Q \cap B )} ( \mathbb V _Q \cap \mathbb V ^+ _2 )  \hookrightarrow \bigcup _{w \in W} Q \times ^{( {} ^w B ) ( s _Q )} ( \mathbb V _Q \cap {} ^w \mathbb V ^+ _2 ) = {\bm F} ^{s _Q}$$
gives an identification of $F _Q$ with a connected component of ${\bm F} ^{s _Q}$. This equips an action of ${\bm Q}$ on $\mathfrak N _Q$, $F _Q$, and $Z _Q$ by restricting the ${\bm G}$-actions on their ambient spaces.

We put
$$\mathbb H _{Q} := \mathbb C \otimes _{\mathbb Z} K ^{\bm Q} ( Z _Q ),$$
where the convolution algebra structure on $K ^{\bm Q} ( Z _Q )$ are equipped by the restrictions of the maps $p _{1}$ and $p _{2}$ from ${\bm Z} \to {\bm F}$ to $Z _Q \to F _Q$.

\begin{lemma}\label{pHecke}
Keep the above setting. Form an increasing sequence of integers
$$1 \le n _1 \le n _2 \le \cdots $$
by requiring that
$$\alpha _i ( s _Q ) < 1 \text{ if and only if } i = n _k \text{ for some }k.$$
Then, we have
\begin{enumerate}
\item $\mathbb H _Q$ is a subalgebra of $\mathbb H$ generated by $\mathcal A [ T ]$ and the set
$$\{ T _i ; i \neq n _k \text{ for some } k \};$$
\item For a pre-admissible parameter $\nu = ( s, \vec{q}, X )$ such that $s \in T$ and $X \in \mathfrak N _Q$, the vector space
$$M _{\nu} ^Q := H _{\bullet} ( \mu _Q ^{-1} ( X ) ^{( s, \vec{q} )})$$
is a $\mathbb H _Q$-module.
\end{enumerate}
\end{lemma}

\begin{proof}
By the condition $(\ref{dsQ})$, we have $\left< \alpha + \beta, s _Q \right> = 1$ for $\alpha, \beta \in R ^+$ if and only if $\left< \alpha, s _Q \right> = 1$ and $\left< \beta, s _Q \right> = 1$. This implies that $Q$ is generated by $T$ and the one-parameter unipotent subgroups corresponding to simple roots $\alpha _i$ (and $- \alpha _i$) such that $\alpha _i ( s _Q ) = 1$.

The variety $F _Q$ decomposes into a product of vector bundles over the flag varieties of simple components of $Q$. By explicit computation, we deduce that the vector bundles we concern are either {\bf a)} the cotangent bundle of the flag variety when the simple component is type $A$, or {\bf b)} the variety ${\bm F}$ for a (possibly smaller) symplectic group which arose as a simple component of $Q$. Moreover, the map $\mu _Q$ is the product of the moment maps of the cotangent bundles of flag varieties of type $A$ and our map ${\bm \mu}$ (for some symplectic group).

Hence, taking account into the argument in \S \ref{fHecke}, both statements are straight-forward modifications of \cite{CG} \S 7.6 and \S 8.6. Thus, we leave the details to the reader.
\end{proof}

\begin{corollary}
Under the assumption of Lemma \ref{pHecke}, $QB$ is a parabolic subgroup of $G$.
\end{corollary}

\begin{proof}
See the first paragraph of the proof of Lemma \ref{pHecke}.
\end{proof}

Let $\mathbb V _U$ be the unique $T$-equivariant splitting of the map $\mathbb V ^+ _2 \longrightarrow \mathbb V ^+ _2 / \mathbb V ^+ _Q$. Let $U$ be the unipotent radical of $QB$. If $X \in \mathbb V$ satisfies $s _Q X = X$, then $s _Q$ has eigenvalue $< 1$ on $\mathfrak u X$. Hence, we have necessarily $\mathfrak u X \subset \mathbb V _U$.

For an admissible parameter $( s, \vec{q}, X )$, we can regard $X = ( X _0 + X _1 \oplus X _2 ) \in \mathbb V$ as an element $X _0 \oplus  X _1 \oplus X _2 \in \mathbb V _2$ so that $s X _i = q _i X _i \in V _1$ for $i = 0, 1$.  

\begin{theorem}[Induction theorem]\label{indt}
We put $P := QB$. Let $P = Q U$ be its Levi decomposition. Let $\nu = ( a,X ) = ( s, \vec{q}, X )$ be an admissible parameter regarded as an element of ${\bm G} \times \mathbb V _2$. Assume $s \in Q$ and $X \in \mathfrak N _Q$. If we have
\begin{eqnarray}
\mathbb V _U ^a \subset \mathfrak u X,\label{denom}
\end{eqnarray}
then we have an isomorphism
$$\mathsf{Ind} ^{\mathbb H} _{\mathbb H _{Q}} M _{\nu} ^{Q} \cong M _{\nu}$$
as $\mathbb H$-modules, where $M _{\nu} ^Q$ is as in Lemma \ref{pHecke}.
\end{theorem}


The rest of this section is devoted to the proof of Theorem \ref{indt}.

By taking $Q$-conjugation if necessary, we assume $X \in \mathbb V ^+ _2$.

Let $W _Q := N _Q ( T ) / T \subset W$. We define
$$W ^Q := \{ w \in W ; \ell ( w ) \le \ell ( v w ) \text{ for all } v \in W _Q\}.$$
Let $w \in W ^Q$. Let $\mathcal O _w$ be the $P$-orbit of $G / B$ which contains $\dot{w} B$. By counting the weights, we have $( \mathbb V ^+ _2 \cap \mathbb V _Q ) \subset ( \mathbb V ^+ _2 \cap {} ^{w} \mathbb V ^+ _2 )$. It follows that $X \in ( \mathbb V ^+ _2 \cap {} ^{w} \mathbb V ^+ _2 )$. Hence, the map
$$( \mathcal E _X \cap \mathcal O _1 ) = \mu _Q ^{-1} ( X ) \ni g B \mapsto g \dot{w} B \in \mathcal E _X \cap \mathcal O _w$$
gives rise to an isomorphism $( \mathcal E _X \cap \mathcal O _1 ) \cong \mathcal E _X \cap \mathcal O ^{s _Q} _w$. Let $B ^-$ be the opposite Borel subgroup of $B$ with respect to $T$. We put $U _w := U \cap {} ^w B ^-$. Since $s _Q$ attracts points of $\mathcal O _w$, we obtain a map
$$\psi _w : \mathcal E _X \cap \mathcal O _w \rightarrow ( \mathcal E _X \cap \mathcal O _w ) ^{s _Q} \cong ( \mathcal E _X \cap \mathcal O _1 )$$
which sends a point $\mathsf{p}$ to $\lim _{N \to \infty} s _Q ^N \mathsf{p}$. We have an expression of a point $g u \dot{w} B \in \mathcal E _X \cap \mathcal O _w$ as $g \in Q$, $g B \in \mu _Q ^{-1} ( X )$, and $u \in U _w$. Let $w _Q$ be the longest element of $W ^Q$.

\begin{lemma}\label{fdesc}
The fiber of the map $\psi _w$ at $g B \in QB$ is given as
$$\psi ^{-1} _w ( gB ) = \{ u \in g U _w g ^{-1} ; u X - X \in g {} ^w \mathbb V ^+ _2 \cap \mathbb V _U \}.$$
In particular, $\psi ^{-1} _{w _Q} ( g B )$ is isomorphic to $\mathsf{Stab} _U ( X )$.
\end{lemma}

\begin{proof}
The variety $\mathcal O _w$ is a $U _w$-fibration over $\mathcal O _1$. The condition $X \in g u \dot{w} \mathbb V ^+ _2$ is equivalent to $( g u ^{-1} g ^{-1} ) X - X \in g \dot{w} \mathbb V ^+ _2$. Moreover, $U$ is $Q$-stable and $( g u ^{-1} g ^{-1} ) X - X \in \mathbb V _U$, which implies the first result. Since $U _{w _Q} = U$ and $g {} ^{w _Q} \mathbb V ^+ _2 \cap \mathbb V _U = \{ 0 \}$, we conclude the second assertion.
\end{proof}

\begin{lemma}\label{dimeqK}
We have
$$\dim H _{\bullet} ( \mathcal E _X ^a \cap \mathcal O _w ) = \dim H _{\bullet} ( \mathcal E _X ^a \cap \mathcal O _1 ).$$
\end{lemma}

\begin{proof}
By the proof of Theorem \ref{ov} and \cite{DLP} 3.9, we have an $\alpha$-partition $\mathcal X _1, \mathcal X _2, \ldots$ of $\mathcal E _X ^a \cap \mathcal O _1$ such that each $\mathcal X _m$ ($m = 1, 2, \ldots$) is a smooth variety without odd-term homology. By condition (\ref{denom}), we see that
$$\dim ( g \dot{w} \mathbb V ^+ \cap ( X + \mathbb V _U ^a ) ) = \dim ( \dot{w} \mathbb V ^+ \cap ( X + \mathbb V _U ^a ) )$$
for all $g = hu \in QU ( s )$ such that $hB \in \mathcal E _X ^a$. We denote its (common) dimension by $d$. Here $( g \dot{w} \mathbb V ^+ \cap ( X + \mathbb V _U ^a ) )$ is an affine space contained in $G X$. The fiber of the map
$$\varphi _{w} : ( \mathcal E _X ^a \cap \mathcal O _w ) \to ( \mathcal E _X ^a \cap \mathcal O _w ^{s _Q} ) \cong ( \mathcal E _X ^a \cap \mathcal O _1 )$$
is isomorphic to a fiber of the following map at $X$:
$$U (s) \times ^{( U (s) \cap \mathrm{Ad} ( g \dot{w} ) B )} ( g \dot{w} \mathbb V ^+ \cap ( X + \mathbb V _U ^a ) ) \to X + \mathbb V _U ^a = U (s) X.$$
In particular, $\varphi _w$ is a smooth affine fibration of relative dimension $\dim U _w (s) + d - \dim \mathbb V _U ^a$. Therefore, Theorem \ref{subB} implies that $\varphi _w ^{-1} ( \mathcal X _m )$ is a vector bundle over $\mathcal X _m$ for each $m$. Hence, we deduce that
$$\sum _{i \ge 0} \dim H _{2i} ( \mathcal X _m ) = \sum _{i \ge 0} \dim H _{2i} ( \varphi _w ^{-1} ( \mathcal X _m) )$$
and
$$H _{2i + 1} ( \mathcal X _m ) = H _{2i + 1} ( \varphi _w ^{-1} ( \mathcal X _m) ) = 0 \text{ for } i = 1, 2, \ldots$$
for each $m$. Since $\varphi _w ^{-1} ( \mathcal X _1 ), \varphi _w ^{-1} ( \mathcal X _2 ), \ldots$ forms an $\alpha$-partition of $\mathcal E _X ^a \cap \mathcal O _w$, we obtain the result.
\end{proof}

We return to the proof of Theorem \ref{indt}.

It is easy to see that
$$\mathcal E _X ^a = \bigsqcup _{w \in W ^Q} ( \mathcal E _X ^a \cap \mathcal O _w )$$
forms an $\alpha$-partition. Together with Theorem \ref{ov} and Lemma \ref{dimeqK}, this implies
\begin{eqnarray}
\dim M _{\nu} = ( \# W ^Q ) \dim M _{\nu} ^Q = ( \# W / \# W _Q ) \dim M _{\nu} ^Q. \label{dimc}
\end{eqnarray}
Moreover, the natural map
$$\imath : M _{\nu} ^Q = H _{\bullet} ( \mathcal E _X ^a \cap \mathcal O _1 ) \hookrightarrow H _{\bullet} ( \mathcal E _X ^a ) = M _{\nu}$$
is injective. Since we have
$$p _1 ( Z _{\le s _i} \cap p _2 ^{-1} ( \mathcal O _1 ) ) \subset \mathcal O _1 \text{ if } i \neq n _k \text{ for some } k = 1, 2, \ldots,$$
the map $\imath$ is an embedding of $\mathbb H _Q$-modules. (The sequence $\{ n _k \} _k$ is borrowed from Lemma \ref{pHecke}.) Hence, we have an induced map
$$\phi : \mathsf{Ind} ^{\mathbb H} _{\mathbb H _{Q}} M _{\nu} ^{Q} \longrightarrow M _{\nu}.$$
Thanks to (\ref{dimc}), we have:
\begin{lemma}
Theorem \ref{indt} follows if $\phi$ is surjective. \hfill $\Box$
\end{lemma}

We return to the proof of Theorem \ref{indt}.

For each $w \in W ^Q$, we define
$$R _{w} := [ \mathcal O _{Z _{\le w ^{-1}}} ] \in K ^{\bm G} ( {\bm Z} ).$$
By the construction of \S 2, we have
$$R _w H _{\bullet} ( \mathcal E _X ^a \cap \mathcal O _1 ) \subset H _{\bullet} ( \mathcal E _X ^a \cap \overline{\mathcal O _w} ) \subset H _{\bullet} ( \mathcal E _X ^a ).$$
Since $W ^Q$ has a partial order $\le _Q$ induced by the Bruhat order, we put
$$H _{\bullet} ( \mathcal E _X ^a ) _{\le w} := \sum _{v \le _Q w; v \in W ^Q} R _v H _{\bullet} ( \mathcal E _X ^a \cap \mathcal O _1 ).$$
Consider the composition map
$$\tau _w : H _{\bullet} ( \mathcal E _X ^a \cap \mathcal O _1 ) \stackrel{R _w \circ}{\longrightarrow} H _{\bullet} ( \mathcal E _X ^a \cap \overline{\mathcal O _w} ) \stackrel{\mathsf{res}}{\longrightarrow} H _{\bullet} ( \mathcal E _X ^a \cap \mathcal O _w ).$$

\begin{lemma}\label{twt}
Theorem \ref{indt} follows if each $\tau _w$ is surjective.
\end{lemma}
\begin{proof}
We have
$$\dim \mathrm{Im} \phi \ge \sum _{w \in W ^Q} \dim \mathsf{gr} H _{\bullet} ( \mathcal E _X ^a ) _{\le w} = \sum _{w \in W ^Q} \dim M _{\nu} ^Q = ( \# W ^Q ) \dim M _{\nu} ^Q,$$
where we used the assumption at the first inequality. Here $\mathsf{gr}$ stands for the graded quotient with respect to some completed order on $W ^Q$ which extends $\le _Q$. By (\ref{dimc}), we conclude that $\phi$ must be surjective.
\end{proof}

We return to the proof of Theorem \ref{indt}.

We have only to prove that each $\tau _w$ is surjective provided if (\ref{denom}) holds. We have an open enbedding
$$\left( p _2 ^{-1} ( \mathcal E _X \cap \mathcal O _1 ) \cap {\bm \pi} ^{-1} ( \mathsf O _{w ^{-1}} ) \right) \subset \left( p _2 ^{-1} ( \mathcal E _X \cap \mathcal O _1) \cap Z _{\le w ^{-1}} \right).$$

\begin{lemma}\label{intw}
For each subset $\mathcal E \subset ( \mathcal E _X \cap \mathcal O _1 )$, we have $\left( p _2 ^{-1} ( \mathcal E  ) \cap {\bm \pi} ^{-1} ( \mathsf O _{w ^{-1}} ) \right) \cong \psi _w ^{-1} ( \mathcal E )$.
\end{lemma}

\begin{proof}
By definition, the LHS is written as:
$$\{ ( g _1 B, g _2 B ) \in {\bm \mu} ^{-1} ( X ) \times ( \mathcal E _X \cap \mathcal O _1 ) ; g _1 ^{-1} g _2 \in B \dot{w} ^{-1} B \}.$$
Since $B \cap Q = {} ^w B \cap Q$, we have $B \dot{w} ^{-1} B = B \dot{w} ^{-1} U$. By taking the right $B$-translation if necessary, we can assume $g _1 \in g _2 U \dot{w}$. This forces $g _1 B$ to live in the fiber of the map $\psi _w$. This implies that $g _2 B$ is completely determined by the data of $\psi _w ^{-1} ( \mathcal E )$ and vice versa.
\end{proof}

Let $A$ be the Zariski closure of $\left< a, s _Q \right> \subset {\bm T}$. The set $\mathfrak u X \subset \mathbb V _U$ is an $A$-stable linear subspace. It follows that
$$S := \mathbb V _U / \mathfrak u X$$
has a $A$-stable splitting in $\mathbb V _U$. Using this splitting, we define
$$\mathcal E _X ^{\sim} := \{ ( g B, X + y ) \in {\bm F} ; g B \in ( \mathcal E _X \cap \mathcal O _1 ), y \in S \}.$$

Each element of $\mathbb V _U$ is contracted to $0$ by the $s _Q$-action. Hence, $S$ has a contraction to $0 \in S$. This gives a contraction
$$\theta : \mathcal E _X ^{\sim} \longrightarrow ( \mathcal E _X \cap \mathcal O _1 )$$
given by collecting $s _Q$-attracting points.

\begin{theorem}\label{tr}
For each $w \in W ^Q$, the intersection of ${\bm \pi} ^{-1} ( \mathsf O _{w ^{-1}})$ and $( {\bm F} \times \mathcal E _X ^{\sim} )$ is transversal inside ${\bm F} ^2$.
\end{theorem}

\begin{proof}
We prove the assertion by induction. The case $w = 1$ is clear. Assume that
\begin{itemize}
\item $w = w ^{\prime} s$ by $w ^{\prime} \in W ^Q$ and $s = s _i$ such that $\ell ( w ) = \ell ( w ^{\prime} ) + 1$;
\item The assertion holds for $w ^{\prime}$;
\end{itemize}
and prove the assertion for $w$. Let $\delta _{in}$ be Kronecker's delta which takes $1$ if $s = s _n$ and $0$ otherwise. For $v = w, w ^{\prime}$, we set
$$\mathcal E ^{\sim} ( v ) := \{ ( g \dot{v} B, X + y ) \in p _1 (( {\bm F} \times \mathcal E ^{\sim} _X ) \cap {\bm \pi} ^{-1} ( \mathsf O _{v ^{-1} }) ) \}.$$
We denote the fibers of the maps $\mathcal E ^{\sim} ( v ) \to ( \mathcal E _X \cap \mathcal O _{v} ^{s _Q} )$ over $g \dot{v} B$ as:
$$F _{v} ( g B ) := \{ ( u, X + y ) \in g U _{v} g ^{-1} \times S ; u ^{-1} X - X \in g {} ^{v} \mathbb V ^+ _2, y \in S \cap g u {} ^{v} \mathbb V ^+ _2 \}.$$
We have
$$\dim {\bm \pi} ^{-1} ( \mathsf O _{w ^{-1}}) = \dim {\bm \pi} ^{-1} ( \mathsf O _{( w ^{\prime} )^{-1}}) - \delta _{in}.$$
Taking account into Lemma \ref{intw}, the failure of the dimension condition of transversal intersection implies
\begin{align}\nonumber
&\dim F _w ( g B ) > \dim F _{w ^{\prime}} ( g B ) \text{ if } s \neq s _n \text{ or }\\
&\dim F _w ( g B ) \ge \dim F _{w ^{\prime}} ( g B ) \text{ if } s = s _n\label{fdim}
\end{align}
for some $g B \in ( \mathcal E _X \cap \mathcal O _1 )$. We assume (\ref{fdim}) (for some $gB$) to deduce contradiction. Being transversal intersection is an open condition. Since the intersection has a contraction to its $s _Q$-fixed point, it suffices to consider the situation near $s _Q$-fixed points. Hence, we replace $F _{v} ( g B )$ in (\ref{fdim}) by the following tangent space version
$$\mathfrak f _v ( g B ) := \{ ( \xi, y ) \in \mathrm{Ad} ( g ) \mathfrak u _{v} \times S ; \xi X \in g {} ^{v} \mathbb V ^+ _2, y \in S \cap g {} ^{v} \mathbb V ^+ _2 \}.$$
Let $\mathfrak u ^v := \mathrm{Lie} ( U \cap {} ^v B )$. It is clear that $( \mathrm{Ad} ( g ) \mathfrak u ^v ) X \subset g {} ^v \mathbb V ^+ _2$ since $X \in g {} ^v \mathbb V ^+ _2$. We have $\mathfrak u = \mathfrak u ^v \oplus \mathfrak u _v$. We put $\Delta ( v ) := \dim ( \mathfrak u ^v Y \cap \mathfrak u _v Y )$. It follows that
$$\dim \mathfrak f _v  ( g B ) = \dim \mathsf{Stab} _{\mathfrak u _v} Y + \dim ( \mathbb V _U \cap {} ^v \mathbb V ^+ _2 ) / \mathfrak u ^v Y + \Delta ( v )$$
where we put $Y := g ^{-1} X$. We have
$$\mathfrak u _w = \mathfrak u _{w ^{\prime}} \oplus \mathfrak g [ \alpha ], \text{ and }\mathfrak u ^w \oplus \mathfrak g [ \alpha ] = \mathfrak u ^{w ^{\prime}}$$
for $\alpha = w ^{\prime} \alpha _i$.
According to the behavior of $\mathfrak g [\alpha] Y$, we have four cases:
\item (Case $\mathfrak g [\alpha] Y \subset \mathfrak u ^w Y \cap \mathfrak u _{w ^{\prime}} Y$) We have
\begin{align*}
& 1 + \delta _{in} + \dim ( \mathbb V _U \cap {} ^w \mathbb V ^+ _2 ) / \mathfrak u ^w Y = \dim ( \mathbb V _U \cap {} ^{w ^{\prime}} \mathbb V ^+ _2 ) / \mathfrak u ^{w ^{\prime}} Y,\\
& \dim \mathsf{Stab} _{\mathfrak u _w} Y - 1 = \dim \mathsf{Stab} _{\mathfrak u _{w ^{\prime}}} Y, \text{ and } \Delta ( w ) = \Delta ( w ^{\prime} ).
\end{align*}
\item (Case $\mathfrak u ^w Y \supset \mathfrak g [\alpha] Y \not\subset \mathfrak u _{w ^{\prime}} Y$) We have
\begin{align*}
& 1 + \delta _{in} + \dim ( \mathbb V _U \cap {} ^w \mathbb V ^+ _2 ) / \mathfrak u ^w Y = \dim ( \mathbb V _U \cap {} ^{w ^{\prime}} \mathbb V ^+ _2 ) / \mathfrak u ^{w ^{\prime}} Y,\\
& \dim \mathsf{Stab} _{\mathfrak u _w} Y = \dim \mathsf{Stab} _{\mathfrak u _{w ^{\prime}}} Y, \text{ and } \Delta ( w ) - 1 = \Delta ( w ^{\prime} ).
\end{align*}
\item (Case $\mathfrak u ^w Y \not\supset \mathfrak g [\alpha] Y \subset \mathfrak u _{w ^{\prime}} Y$)
\begin{align*}
& \delta _{in} + \dim ( \mathbb V _U \cap {} ^w \mathbb V ^+ _2 ) / \mathfrak u ^w Y = \dim ( \mathbb V _U \cap {} ^{w ^{\prime}} \mathbb V ^+ _2 ) / \mathfrak u ^{w ^{\prime}} Y,\\
& \dim \mathsf{Stab} _{\mathfrak u _w} Y - 1 = \dim \mathsf{Stab} _{\mathfrak u _{w ^{\prime}}} Y, \text{ and } \Delta ( w ) + 1 = \Delta ( w ^{\prime} ).
\end{align*}
\item (Case $\mathfrak u ^w Y \not\supset \mathfrak g [\alpha] Y \not\subset \mathfrak u _{w ^{\prime}} Y$)
\begin{align*}
& \delta _{in} + \dim ( \mathbb V _U \cap {} ^w \mathbb V ^+ _2 ) / \mathfrak u ^w Y = \dim ( \mathbb V _U \cap {} ^{w ^{\prime}} \mathbb V ^+ _2 ) / \mathfrak u ^{w ^{\prime}} Y,\\
& \dim \mathsf{Stab} _{\mathfrak u _w} Y = \dim \mathsf{Stab} _{\mathfrak u _{w ^{\prime}}} Y, \text{ and } \Delta ( w ) = \Delta ( w ^{\prime} ).
\end{align*}
\item This case-by-case analysis claims that we cannot achieve the infinitesimal version of (\ref{fdim}). Hence, we have contradiction. It follows that the intersection of ${\bm \pi} ^{-1} ( \mathsf O _{w ^{-1}})$ and $( {\bm F} \times \mathcal E _X ^{\sim} )$ must have proper dimension inside ${\bm F} ^2$ under the induction hypothesis.

Now the linear independence of the normal vectors follows as an immediate consequence of the fact that they are concentrated on the first factor and $S$ on the second factor (of ${\bm F} \times \mathcal E _X ^{\sim}$), or the diagonal part (of ${\bm \pi} ^{-1} ( \mathsf O _{w ^{-1}})$), respectively.

Therefore, the induction proceeds and we obtain the result.
\end{proof}

\begin{lemma}\label{ti}
The map $\tau _{w}$ is an isomorphism.
\end{lemma}

\begin{proof}
By \cite{CG} 2.7.26 and Theorem \ref{tr}, we deduce that the map $\tau _{w}$ induces an isomorphism
\begin{eqnarray}
H _{\bullet} ( \mathcal E _X ^{\sim} ) \stackrel{\cong}{\longrightarrow} H _{\bullet} ( p _1 ( ( {\bm F} \times \mathcal E _X ^{\sim} ) \cap {\bm \pi} ^{-1} ( \mathsf O _{w ^{-1}}) ) ).\label{isow}
\end{eqnarray}
The spaces appearing in the homologies are given as fibrations over $( \mathcal E _X \cap \mathcal O _1 )$ and $( \mathcal E _X \cap \mathcal O _w )$ with its fiber linear subspaces of $S$. Here $\mathcal E _X ^{\sim}$ has larger fiber. We switch to algebraic $K$-theory by Theorem \ref{ov} and Lemma \ref{intw}. We have
$$R ( A ) _a \otimes _{R ( A )} K ^A ( \mathcal E _X ^{\sim} ) \cong R ( A ) _a \otimes _{R ( A )} K ^A ( \mathcal E _X ^a \cap \mathcal O _1 ) = R ( A ) _a \otimes _{\mathbb Z} K ( \mathcal E _X ^a \cap \mathcal O _1 )$$
by the Thomason localization theorem and the fact that $A$ fixes $( \mathcal E _X ^a \cap \mathcal O _1 )$. Hence, the map $\tau _{w}$ itself is surjective if $[ Y ], [ \theta ^{-1} ( Y ) ] \in R ( A ) _a \otimes _{R ( A )} K ^A ( \mathcal E _X ^{\sim} )$ define the same cycle up to an invertible factor for each $A$-stable closed subvariety $Y \subset ( \mathcal E _X \cap \mathcal O _1 )$. This is true if the alternating sum of the Koszul complex of $S$ is invertible in $R ( A ) _a$. This is equivalent to $S ^a = 0$, which is further re-phrased as
$$\mathbb V _U ^a \subset \mathfrak u X.$$
This is (\ref{denom}).
\end{proof}

We return to the proof of Theorem \ref{indt}.

Thanks to Lemma \ref{ti}, we have finished the proof of Theorem \ref{indt} by Lemma \ref{twt}.

\section{Exotic Springer correspondence}\label{esc}
We keep the setting of the previous section.

As we see in \S \ref{gen conv}, we know that the action of $\mathbb H$ on $M _{\nu}$ factors through the isomorphism
$$\mathbb C \otimes _{\mathbb Z} K ( {\bm Z} ^a ) \stackrel{{\rm RR}}{\longrightarrow} H _{\bullet} ( {\bm Z ^a} ) \cong \mathbb H _a.$$

Let $\left< \mathbb C [ \mathfrak t ] ^W _+ \right>$ be the ideal of $\mathbb C [ \mathfrak t ]$ generated by the set of $W$-invariant polynomials without constant terms. Let $\mathbb C [ W ] \# \left( \mathbb C [ \mathfrak t ] / \left< \mathbb C [ \mathfrak t ] ^W _+ \right> \right)$ be the smash-product, which means that its product is given as
$$( w _1, f _1 ) ( w _2, f _2 ) := ( w _1 w _2, f _1 w _1 ( f _2 ) ) \text{ for } w _1, w _2 \in W, f _1, f _2 \in \mathbb C [ \mathfrak t ] / \left< \mathbb C [ \mathfrak t ] ^W _+ \right>.$$

It is clear that ${\bm F} ^{a _0} \cong F$, ${\bm Z} ^{a _0} \cong Z$, and the restriction of the natural projections ${\bm Z} \to {\bm F}$ restrict to natural projections $Z \to F$.

\begin{proposition}\label{eW}
We have an isomorphism
$$\mathbb C [ W ] \# \left( \mathbb C [ \mathfrak t ] / \left< \mathbb C [ \mathfrak t ] ^W _+ \right> \right) \cong H _{\bullet} ( {\bm Z} ^{a _0} )$$
as algebras. 
\end{proposition}

\begin{proof}
We have
$$H _{\bullet} ( {\bm Z} ^{a _0} ) \cong \mathbb C _{a _0} \otimes _{R ( {\bm G} )} K ^{\bm G} ( Z ).$$
Here the RHS is written as
$$\mathbb C \otimes _{R ( G )} \mathbb H / ( \mathbf q _0 = - \mathbf q _1 = \mathbf q _2 = 1 ).$$
Thus, we have
$$H _{\bullet} ( {\bm Z} ^{a _0} ) \cong \mathbb C \otimes _{R ( G )} \mathbb C [ \widetilde{W} ],$$
where $\widetilde{W} := W \ltimes X ^* ( T )$ is the affine Weyl group of type $C _n ^{(1)}$. (Here $\mathbb C$ is the $R ( G )$-module given by the evaluation at $1 \in G$. The algebra $R ( G )$ acts on $\mathbb C [ \widetilde{W} ]$ by $R ( G ) \cong \mathbb Z [ X ^* ( T ) ] ^W$.) Thus, it suffices to show
$$\mathbb C [ X ^* ( T ) ] / \mathbb C [ X ^* ( T ) ] \mathfrak m _1 ^{\sim} \cong \mathbb C [ \mathfrak t ] / \left< \mathbb C [ \mathfrak t ] ^W _+ \right>,$$
where $\mathfrak m _1 ^{\sim} \subset \mathbb C [ X ^* ( T ) ] ^W = \mathbb C [ T ] ^W$ is the defining ideal of the image of $1 \in T$ in $\mathrm{Spec} \mathbb C [ T ] ^W$. This follows from the fact that the neighborhoods of $1 \in T$ and $0 \in \mathfrak t$ are $W$-equivariantly diffeomorphic through the exponential map.
\end{proof}

\begin{corollary}
Keep the setting of \ref{eW}. We have a surjection
$$H _{\bullet} ( {\bm Z} ^{a _0} ) \longrightarrow \!\!\!\!\! \rightarrow \mathbb C [ W ].$$
\end{corollary}

\begin{proof}
Keep the notation of the proof of Theorem \ref{eW}. We have
$$\left< \mathbb C [ \mathfrak t ] ^W _+ \right> \subset \mathfrak m _1 \subset \mathbb C [ \mathfrak t ],$$
where $\mathfrak m _1$ is the defining ideal of $0 \in \mathfrak t$. Since $0$ is a $W$-fixed point of $\mathfrak t$, we deduce that $\mathfrak m_1$ is a $W$-invariant maximal ideal. It follows that
$$
H _{\bullet} ( {\bm Z} ^{a _0} ) \cong \mathbb C [ W ] \# \left( \mathbb C [ \mathfrak t ] / \left< \mathbb C [ \mathfrak t ] ^W _+ \right> \right) \longrightarrow \!\!\!\!\! \rightarrow \mathbb C [ W ] \# \left( \mathbb C [ \mathfrak t ] / \mathfrak m _1 \right) \cong \mathbb C [ W ]
$$
as desired.
\end{proof}

\begin{theorem}[Exotic Springer correspondence]\label{eS}
There exist one-to-one correspondences between the sets of the following three kinds of objects:
\begin{itemize}
\item a strict marked partition $\sigma$;
\item the $G$-orbit of $\mathfrak N$ given as $G \mathbf v _{\sigma}$;
\item an irreducible $W$-module.
\end{itemize}
\end{theorem}

\begin{remark}
Our proof of Theorem \ref{eS} does not tell which representation is obtained from a given orbit. Such information can be found in \cite{K08}, which employs totally different argument.
\end{remark}

\begin{proof}[Proof of Theorem \ref{eS}]
Let $\mathcal P$ be the set of isomorphism classes of $G$-equivariant irreducible perverse sheaves on $\mathfrak N$. Each $I \in \mathcal P$ is isomorphic to the minimal extension from a $G$-orbit of $\mathfrak N$. By Proposition \ref{2con}, the (perverse) sheaf $I$ must be the extension of a constant sheaf on a $G$-orbit. This implies $\# \mathcal P \le \# ( G \backslash \mathfrak N )$. Let $S$ be the set of strict normal forms. By Proposition \ref{wN1} 1), we have $\# ( G \backslash \mathfrak N ) \le \# S$. Hence, we have
\begin{eqnarray}
\# \mathsf{Irrep} W \le \# \mathcal P \le \# ( G \backslash \mathfrak N ) \le \# S \le \# \mathsf{Irrep} W,\label{ineqeq}
\end{eqnarray}
where the first inequality comes from Theorem \ref{Gd} and the last inequality is Proposition \ref{wN1} 2). This forces all the inequalities in (\ref{ineqeq}) to be equalities as required.
\end{proof}

The following is a summary of the consequences of \S \ref{gen conv} Theorem \ref{Gd}:

\begin{theorem}[Ginzburg, \cite{CG} \S 8.5]\label{sds}
Let $a$ be a finite pre-admissible element. Let $L$ be an irreducible $\mathbb H _a$-module. Then, there exists a unique ${\bm G} ( a )$-orbit $\mathcal O \subset \mathfrak N ^a$ with the following properties:
\begin{enumerate}
\item There exists a surjective $\mathbb H _a$-module homomorphism $M _{(a, X)} \rightarrow L$ for every $X \in \mathcal O$;
\item If $L$ appears in the composition factor of $M _{(a, Y)}$ as $\mathbb H _a$-modules for some $Y \in \mathfrak N ^a$, then we have $Y \in \overline{\mathcal O}$. \hfill $\Box$
\end{enumerate}
\end{theorem}

Theorems \ref{eS} and \ref{sds} claim that each strict marked partition $\sigma$ gives a unique simple quotient of $M _{(a _0, \mathbf v _{\sigma})}$. We denote this $W$-module by $L _{\sigma}$ or $L _X$ for $X \in G \mathbf v_{\sigma}$, depending on the situation.

\begin{corollary}\label{BM}
Keep the setting of Theorem \ref{eS}. A $\mathbb C [ W ]$-module $M _{(a _0, X)}$ contains $L _{\sigma}$ only if $X \in \overline{G \mathbf v _{\sigma}}$ holds. \hfill $\Box$
\end{corollary}

\section{A deformation argument on parameters}\label{dap}

We retain the setting of the previous section.

\begin{theorem}\label{reg}
Let $a = ( s, \vec{q} ) \in {\bm T}$ be an admissible element such that $\mathcal C _a = \{ [ 1, n ] \}$. Then, there exists an admissible element $a ^{\prime} := ( s ^{\prime}, \vec{q} ^{\prime} )$ such that
\begin{itemize}
\item The $s ^{\prime}$-action on $V _1$ has only positive real eigenvalues;
\item We have $q _0 ^{\prime}, q _1 ^{\prime}, q _2 ^{\prime} \in \mathbb R ^{\times} _{> 0}$;
\item We have equalities $\mathfrak N _+ ^a = \mathfrak N _+ ^{a ^{\prime}}$ and $G ( s ) = G ( s ^{\prime} )$.
\end{itemize}
Moreover, we have an isomorphism $\mathbb H _a \cong \mathbb H _{a ^{\prime}}$ as algebras.
\end{theorem}

\begin{proof}
Let $N$ be the largest positive integer such that $1, q _{2}, \ldots, q _2 ^{N}$ are distinct. (If $q _2$ is not a root of unity, then we regard $N = \infty$.) For each $i = 1, \ldots, n$, we set  $\chi _i : = \epsilon _i ( s )$. By rearranging $s$ by the $W$-action if necessary, we assume $| \chi _i  |  \ge 1$ (if $N = \infty$) or $\chi _i = q _2 ^{j}$ for some $j \in \frac{1}{2} [0, N]$. We set $\mathtt E := \{ \chi _i ; 1 \le i \le n \}$. We choose a representative $j _{0} \in [ 1, n ]$ which satisfies the following condition:
\begin{itemize}
\item If $\pm 1 \in \mathtt E$, then we require $\chi _{j _{0}} = \pm 1$;
\item If $\pm q _2 ^{1 / 2} \in \mathtt E$, then we require $\chi _{j _{0}} = \pm q _2 ^{1 / 2}$;
\item If $\pm q _2 ^{1 / 2} \not\in \mathtt E$ and $\pm q _2 ^{- 1 / 2} \in \mathtt E$, then we require $\chi _{j _{0}} = \pm q _2 ^{- 1 / 2}$;
\end{itemize}
For each pair $i, j \in [1,n]$, we have
\begin{eqnarray}
\chi _i ^{\kappa _{i, j}} = \chi _j ^{\kappa _{i, j} ^{\prime}} q _2 ^{m _{i, j}} \text{ for some } \kappa _{i, j}, \kappa _{i, j} ^{\prime} \in \{ \pm 1\}, m _{i,j} \in [0, n]. \label{relc}\end{eqnarray}
Since $q _2$ is not a root of unity of order $\le 2n$, it follows that the choice of $m _{i, j}$ is at most one if $( \kappa _{i, j}, \kappa _{i, j} ^{\prime} )$ is fixed. For each pair $( i, j )$ in $[1,n]$, we set $\mathtt I _{( i, j )}$ to be the set of triples $( \kappa _{i, j}, \kappa _{i,j} ^{\prime}, m _{i,j} )$ which satisfies (\ref{relc}). Choose two real numbers $q \gg q _2 ^{\prime} \gg 1$ such that $q$ and $q _2 ^{\prime}$ have no algebraic relation. Then, we set
$$( \chi _i  ^{\prime} ) ^{\kappa _{i, j _0}} := \begin{cases}
( q _2 ^{\prime} ) ^{m _{i, j _{0}}} & ( \chi _{j _{0}} = \pm 1 ) \\
( q _2 ^{\prime} ) ^{m _{i, j _{0}} + \kappa _{i, j _{0}} ^{\prime} / 2} & ( \chi _{j _{0}} = \pm q _2 ^{1/2})\\
( q _2 ^{\prime} ) ^{m _{i, j _{0}} - \kappa _{i, j _{0}} ^{\prime} / 2} & ( \chi _{j _{0}} = \pm q _2 ^{- 1/2})\\
q ( q _2 ^{\prime} ) ^{m _{i, j _{0}}} & ( \chi _{j _{0}} \neq \pm 1, \pm q _2 ^{\pm 1/2} ) \\
\end{cases}.$$
Since the relation (\ref{relc}) for $(i, j)$ is determined by that of $(i, j _{0})$ and $(j, j _{0})$ for each pair $i, j$ in $[1,n]$, it follows that
\begin{eqnarray}
( \chi _i ^{\prime} ) ^{\kappa _{i, j}} = ( \chi _j ^{\prime} ) ^{\kappa _{i, j} ^{\prime}} ( q _2 ^{\prime} ) ^{m _{i, j}} \text{ for some } \kappa _{i, j}, \kappa _{i, j} ^{\prime} \in \{ \pm 1\}, m _{i,j} \in [0, n]\label{relp}
\end{eqnarray}
for all $( \kappa _{i,j}, \kappa _{i,j} ^{\prime}, m _{i, j} ) \in \mathtt I _{(i, j)}$. Conversely, we have $( \kappa _{i,j}, \kappa _{i,j} ^{\prime}, m _{i, j} ) \in \mathtt I _{(i, j)}$ if the relation (\ref{relp}) holds. It is clear that $\chi _i ^2 = 1$ if and only if $( \chi _i ^{\prime} ) ^2 = 1$. We put $s ^{\prime} \in T$ so that $\epsilon _i ( s ^{\prime} ) = \chi _i ^{\prime}$ for each $i = 1, 2, \ldots, n$. By the above consideration, it follows that $\mathfrak g ( s ^{\prime} ) = \mathfrak g ( s )$. Since both $G ( s ^{\prime} )$ and $G ( s )$ are connected by Theorem \ref{Steinberg}, we deduce $G ( s ) = G ( s ^{\prime} )$.

Since the relation of (\ref{relc}) is preserved, we have $V _2 ^{( s, q _2 ) } = V _2 ^{( s ^{\prime}, q _2 ^{\prime} )}$. If we have $\chi _i ^{\kappa _i} = q _k$ for some $i \in [1, n]$, $\kappa _i \in \{ \pm 1\}$, and $k = 0, 1$, then we set $q _k ^{\prime} := ( \chi _i ^{\prime} ) ^{\kappa _i}$. Otherwise, we set $q _k ^{\prime}$ $(i = 0, 1)$ to be an arbitrary real number which is not an eigenvalue of $s ^{\prime}$ on $V _1$. (I.e. not equal to any of $( \chi _i ^{\prime} ) ^{\pm 1}$.) Since we have infinitely many possibilities, we can assume $q _0 ^{\prime} \neq q _1 ^{\prime}$ and $q _k ^{\prime} \gg 1$ in this case. This gives $\mathbb V ^{a} = \mathbb V ^{a ^{\prime}}$ by setting $a ^{\prime} := ( s ^{\prime}, \vec{q} ^{\prime} )$. We have $q _0 ^{\prime} \neq q _1 ^{\prime}$ in all cases since $q _0 \neq q _1$. Hence, the isomorphism $\mathbb V ^a _2 \cong \mathbb V ^a$ implies $\mathbb V ^{a ^{\prime}} _2 \cong \mathbb V ^{a ^{\prime}}$.

Therefore, as subvarieties of ${\bm F}$ and $\mathfrak N _2$, we have equalities
$${\bm F} ^a = \bigcup _{w \in W} G ( s ) \times ^{{} ^w B ( s )} ( {} ^w \mathbb V ^+ \cap \mathbb V ^a) = \bigcup _{w \in W} G ( s ^{\prime} ) \times ^{{} ^w B ( s ^{\prime} )  } ( {} ^w \mathbb V ^+ \cap \mathbb V ^{a ^{\prime}} ) = {\bm F} ^{a ^{\prime}}$$
and $\mathfrak N _+ ^a = \mathfrak N _+ ^{a ^{\prime}}$.

The projection map ${\bm F} ^a \to \mathfrak N _+ ^a$ is induced by the projection ${\bm \mu}$. Hence, so is ${\bm F} ^{a ^{\prime}} \to \mathfrak N _+ ^{a ^{\prime}}$. Therefore, we have an equality of convolution algebras
$$\mathbb H _a \cong \mathbb C \otimes _{\mathbb Z} K  ( {\bm Z} ^a ) = \mathbb C \otimes _{\mathbb Z} K  ( {\bm Z} ^{a ^{\prime}} ) \cong \mathbb H _{a ^{\prime}},$$
which proves the last assertion.

Since $q, q _2 ^{\prime} \gg 1$, each of $q _i ^{\prime}$ ($i = 0, 1$) is positive real. This verifies the requirement about $\vec{q} ^{\prime}$ as desired.
\end{proof}

\begin{proposition}\label{inclW}
Let $a = ( s, q _0, q _1, q _2) \in {\bm T}$ be an admissible element such that:
\begin{itemize}
\item We have $\mathcal C _a = \{ [1, n] \}$;
\item The $s$-action on $V _1$ has only positive real eigenvalues;
\item We have $V _1 ^{( s, q _1 )} = \{ 0 \}$;
\item Each $q _i$ $(i = 0, 1, 2)$ is a positive real number;
\end{itemize}
Let $\underline{a} := ( s, q _0, q _2 )$ and let
$$\log \underline{a} := ( \log s, r _0, r _2 ), \text{ where } q _0 = e ^{r _0}, q _2 = e ^{r _2}.$$
Let $A$ be the Zariski closure of $\left< \underline{a} \right> \subset {\bm T}$. Then $H _{\bullet} ^A ( Z )$ is a $\mathbb C [\mathfrak a]$-algebra such that
\begin{enumerate}
\item The quotient of $H _{\bullet} ^A ( Z )$ by the ideal generated by functions of $\mathbb C [ \mathfrak a ]$ which is zero along $\log \underline{a}$ is isomorphic to $H _{\bullet} ( Z ^{\underline{a}} )$;
\item The images of the natural inclusions $\mathbb C [ W ] \subset H _{\bullet} ( Z ) \subset H _{\bullet} ^A ( Z )$ induces an injection
$$\mathbb C [ W ] \hookrightarrow H _{\bullet} ( Z ^{\underline{a}} ) = H _{\bullet} ( {\bm Z} ^{a} ).$$
\end{enumerate}
Moreover, we have
$$\mathbb C [ \mathfrak a ] \otimes H _{\bullet} ( \mathcal E _X ) \cong H _{\bullet} ^A ( \mathcal E _X ) \text{ for } X \in \mathfrak N ^{\underline{a}}$$
as a compatible $( \mathbb C [ W ], \mathbb C [ \mathfrak a ] )$-module, where $W$ acts on $\mathfrak a$ trivially. 
\end{proposition}

\begin{corollary}\label{Wmod}
Keep the setting of Proposition \ref{inclW}. We have
$$M _{(a _0, X)} = H _{\bullet} ( \mathcal E _X ) \cong H _{\bullet} ( \mathcal E _X ^a ) = M _{(a, X)}$$
as $\mathbb C [ W ]$-modules. \hfill $\Box$
\end{corollary}

The rest of this section is devoted to the proof of Proposition \ref{inclW}.

\begin{lemma}\label{Ac}
Keep the setting of Proposition \ref{inclW}. Then, $A$ is connected.
\end{lemma}
\begin{proof}
The group $A$ is defined to be the spectrum of the quotient of $\mathbb C [ T _1 ]$ by the ideal generated by monomials $m$ such that $m ( s, q _0, q _2 ) = 1$. Since all the values of $\epsilon _i ( s )$, $q _0, q _1$ are positive real number, the conditions $m ( s, q _0, q _2 ) = 1$ and $m ( s ^r, q _0 ^r, q _2 ^r ) = 1$ are the same for all $r \in \mathbb R _{>0}$, where the branch of powers are taken so that all of $\epsilon _i ( s ^r ), q _0 ^r, q _2 ^r$ ($i = 1, \ldots, n$) are positive real numbers. It follows that a monomial $m \in \mathbb C [ T _1 ]$ satisfies $m ( s, q _0, q _2 ) ^k = 1$ for some positive integer $k$ if and only if $m ( s, q _0, q _2 ) = 1$. Therefore, such monomials form a saturated $\mathbb Z$-sublattice of $X ^* ( T _1 )$. In particular, its quotient lattice is a free $\mathbb Z$-lattice, which implies that $A$ is connected.
\end{proof}

We return to the proof of Proposition \ref{inclW}.

For each $m \ge 0$, let $ET _m := ( \mathbb C ^m \backslash \{ 0 \} ) ^{\dim {\bm T}}$ be a variety such that $i$-th $\mathbb C ^{\times}$-factor of ${\bm T} = ( \mathbb C ^{\times} ) ^{\dim {\bm T}}$ acts as dilation of the $i$-th factor for each $1 \le i \le n + 3$. By the standard embedding $\mathbb C ^m \hookrightarrow \mathbb C ^{m + 1}$ sending $( x )$ to $( x , 0 )$, we form a sequence of ${\bm T}$-varieties
\begin{eqnarray}
\emptyset = ET _0 \hookrightarrow ET _1 \hookrightarrow ET _2 \hookrightarrow \cdots.\label{Borel}
\end{eqnarray}
We define $ET := \varinjlim _m ET _m$, which is an ind-quasiaffine scheme with free ${\bm T}$-action. When we consider the homology of $ET$, we refer to the homology of its underlying classical topological space $\bigcup _{m \ge 0} ET_m$. Since $ET$ is contractible manifold with respect to the classical topology, we regard $ET$ as the universal vector bundle of each subgroup of ${\bm T}$. (Hence we regard $BA := A \backslash ET$ in the below.)

\begin{corollary}[of Lemma \ref{Ac}]\label{oba}
Keep the above setting. We have $H ^{odd} ( BA ) = 0$. \hfill $\Box$
\end{corollary}

We return to the proof of Proposition \ref{inclW}.

For a $A$-variety $\mathcal X$, we set
$$\mathcal X _A := \triangle A \backslash \left( ET \times \mathcal X \right),$$
where $\triangle A$ represents the diagonal action of $A$. We have a forgetful map
$$f ^A _{\mathcal X} : \mathcal X _A \rightarrow BA = A \backslash ET.$$
Let $\mathbb D _{\mathcal X} ^A$ be the relative dualizing sheaf with respect to $f _{\mathcal X} ^A$ (c.f. Bernstein-Lunts \cite{BL} \S 1.6). We define
$$H _{i} ^A ( \mathcal X ) \cong H ^{-i} ( \mathcal X _A, \mathbb D _{\mathcal X} ^A ).$$
In the below, we understand that $H _{\bullet} ^A ( \mathcal X ) := \bigoplus _{m} H _{m} ^A ( \mathcal X )$. Notice that this homology group is the same as the one obtained by replacing $ET$ with an ind-object of the direct system $\{ ET _m \}$ and take the limit of the associated inverse system since $\mathcal X$ is homotopic to a finite dimensional CW-complex. The projection maps $p _i : Z _A \rightarrow F _A$ ($i = 1, 2$) equip $H _{\bullet} ^A ( Z )$ a structure of convolution algebra. It is straight-forward to see that the diagonal subsets $\triangle F \subset Z$ and $( \triangle F ) _A \subset Z _A$ represents $1 \in H _{\bullet} ( Z )$ and $1 \in H _{\bullet} ^A ( Z )$, respectively.

\begin{lemma}
The algebra $H _{\bullet} ^A ( Z )$ contains $H _{\bullet} ( Z )$ as its subalgebra. In particular, we have $\mathbb C [ W ] \subset H _{\bullet} ^A ( Z )$ as subalgebras. Moreover, the center of $H _{\bullet} ^A ( Z )$ contains $H ^{\bullet} ( BA ) [ ( \triangle F ) _{A} ] \subset H _{\bullet} ^A ( Z )$.
\end{lemma}

\begin{proof}
In the Leray spectral sequence
$$H ^i ( BA ) \otimes H _{j} ( Z ) \Rightarrow H _{- i + j} ^A ( Z ),$$
we have $H ^{odd} ( BA ) = 0$ and $H _{odd} ( Z ) = 0$ (since $Z$ is paved by affine spaces). It follows that this spectral sequence degenerates at the level of $E _2$-terms. Moreover, the image of the natural map $\imath : H _j ( Z ) \hookrightarrow H _{j} ^A ( Z )$ represents cycles which are locally constant fibration over the base $BA$. It follows that the map $\imath$ is an embedding of convolution algebras.\\
Multiplying $H ^{\bullet} ( BA )$ is an operation along the base $BA$, which commutes with the convolution operation (along the fibers of $f ^A _{Z}$). It follows that $H ^{\bullet} ( BA ) \rightarrow H ^{\bullet} ( BA ) [ ( \triangle F ) _A ] \subset H _{\bullet} ^A ( Z )$ is central subalgebra as desired.
\end{proof}

We return to the proof of Proposition \ref{inclW}.

By the Thomason localization theorem (see e.g. \cite{CG} \S 8.2), we have an isomorphism
$$R ( A ) _{\underline{a}} \otimes _{R ( A )} K ^A ( Z ^{\underline{a}} ) \cong R ( A ) _{\underline{a}} \otimes _{R ( A )} K ^A ( Z )$$
as algebras. For each of $\mathcal X = Z$, or $Z ^{\underline{a}}$, we have an embedding
$$\imath : K ^A ( \mathcal X ) \hookrightarrow \varprojlim _m K ^A ( ET _m \times \mathcal X )\cong \varprojlim _m K ( A \backslash ( ET _m \times \mathcal X ) )$$
obtained by pulling back an $A$-equivariant vector bundle on an irreducible component of $\mathcal X$ to each $ET_m \times \mathcal X$ by the second projection. The latter inverse limits are formed by the pullbacks via the closed embeddings coming from (\ref{Borel}). Here the last inverse limit has a natural topology whose open sets are formed by the formal sum of vector bundles which are trivial on $K ( A \backslash ( ET _m \times \mathcal X ) )$ for some fixed choice of $m$. By construction, the image of $\imath$ must be dense open with respect to the topology on the RHS.

We regard the RHS as a substitute of $K ( \mathcal X _A )$. Let $\mathbb C [[ \mathfrak a ]] _{a}$ and $\mathbb C [ \mathfrak a ] _a$ be the formal power series ring and the localized ring of $\mathbb C [ \mathfrak a ]$ along $\log \underline{a}$, respectively. The Chern character map relative to $BA$ gives an isomorphism
$$\mathbb C [[ \mathfrak a ]] _{a} \otimes _{\mathbb C [ \mathfrak a ]} H ^A _{\bullet} ( Z ^{\underline{a}} ) \cong \mathbb C [[ \mathfrak a ]] _{a} \otimes _{\mathbb C [ \mathfrak a ]} H ^A _{\bullet} ( Z ).$$
By restricting this to the sum of vectors of finitely many degrees, we obtain
\begin{eqnarray}
\mathbb C [ \mathfrak a ] _{a} \otimes _{\mathbb C [ \mathfrak a ]} H ^A _{\bullet} ( Z ^{\underline{a}} ) \cong \mathbb C [ \mathfrak a ] _{a} \otimes _{\mathbb C [ \mathfrak a ]} H ^A _{\bullet} ( Z ).\label{localH}
\end{eqnarray}
Since localization along $\mathbb C [ \mathfrak a ] _a$ commutes with the quotient by its unique maximal ideal, we deduce the first assertion.

The isomorphism (\ref{localH}) is an algebra isomorphism. This implies that $1 \in \mathbb C [ W ]$ goes to $1 \in H _{\bullet} ( Z ^{\underline{a}} )$. It follows that each of $s _i$ goes to a non-zero element of $H _{\bullet} ( Z ^{\underline{a}} )$ with its square equal to $1$. By construction, there exists $f _i \in \mathbb C ( \mathfrak a )$ ($i = 1, \ldots n$) such that $1, f _1 s _1, \ldots f _n s _n \in H _{\bullet} ^A ( Z ^{\underline{a}} )$ define linearly independent vectors in $H _{\bullet} ( Z ^{\underline{a}} )$. It follows that $f _i ^2 \in \mathbb C [ \mathfrak a ]$. This forces $f _i \in \mathbb C [ \mathfrak a ]$ (since a polynomial ring is integrally closed), which implies that the images of $1, s _1, \ldots, s _n \in H _{\bullet} ( Z ^{\underline{a}} )$ are linearly independent. This verifies the second assertion.

The vector space $H ^A _{\bullet} ( \mu ^{-1} ( X ) )$ admits an action of $H ^A _{\bullet} ( Z )$. By the Leray spectral sequence, we have
\begin{eqnarray}
H ^{\bullet} ( BA ) \otimes H _{\bullet} ( \mu ^{-1} ( X ) ) \Rightarrow H ^A _{\bullet} ( \mu ^{-1} ( X ) ).\label{stdsp}
\end{eqnarray}
By Theorem \ref{ov} and Corollary \ref{oba}, we know that $H _{odd} ( \mu ^{-1} ( X ) ) = 0 = H ^{odd} ( BA )$. It follows that (\ref{stdsp}) is $E _2$-degenerate, which proves the last part of Proposition \ref{inclW}.

\section{Main Theorems}\label{MAIN}
We retain the setting of \S \ref{fHecke}.

\begin{theorem}[Deligne-Langlands type classification]\label{DLmain}
Let $a \in {\bm G}$ be a finite pre-admissible element. Then, $\mathfrak R _a$ is in one-to-one correspondence with the set of isomorphism classes of simple $\mathbb H _{a}$-modules.
\end{theorem}

\begin{proof}
By definition, each element of $\mathfrak R _a$ corresponds to at least one isomorphism class of $\mathbb H _a$-modules. Since $a$ is finite, each irreducible direct summand of $( \mu ^{a} _+ ) _* \mathbb C _{F ^{a} _+}$ is the minimal extension of a local system (up to degree shift) from a ${\bm G} ( a )$-orbit $\mathbb O$. By Theorem \ref{pi1 desc}, a ${\bm G} ( a )$-equivariant local system on $\mathbb O$ is a constant sheaf. As a result, every element of $\mathfrak R _a$ corresponds to at most one irreducible module as desired.
\end{proof}

\begin{theorem}[Effective Deligne-Langlands type classification]\label{EDL}
Let $a \in {\bm G}$ be an admissible element. Then, the set $\Lambda _{a}$ is in one-to-one correspondence with the set of isomorphism classes of simple $\mathbb H _{a}$-modules.
\end{theorem}

\begin{proof}
The proof is given at the end of this section.
\end{proof}

As in Remark \ref{remHecke}, the quotient $\mathbb H / ( \mathbf q _0 + \mathbf q _1 )$ is isomorphic to an extended Hecke algebra $\mathbb H _B$ of type $B _n ^{( 1 )}$ with two parameters. Hence, we have

\begin{corollary}[Effective Deligne-Langlands type classification for type $B$]\label{BDL}
Let $a = ( s, q _0, - q _0, q _2 ) \in {\bm G}$ be a pre-admissible element such that $- q _0 ^2 \neq q _2 ^{\pm m}$ holds for every $0 \le m < n$. Then, the set $\Lambda _{a}$ is in one-to-one correspondence with the set of isomorphism classes of simple $\mathbb H _{a}$-modules. \hfill $\Box$
\end{corollary}

\begin{remark}\label{inv}
The Dynkin diagram of type $C _n ^{( 1 )}$ is written as:\\
\xymatrix@R=2pt{
& & *{0} & *{1} & *{2} &  & *{n - 2} & *{n - 1} & *{n}\\
& & *{\circ} \ar@{=}[r]|{>} & *{\circ} \ar@{-}[r] & *{\circ} \ar@{-}[rr]|{\cdots\cdots} &  & *{\circ} \ar@{-}[r] & *{\circ} \ar@{=}[r]|{<} & *{\circ}
}\\

This Dynkin diagram has a unique non-trivial involution $\varphi$. We define $t _0, t _1, t _n$ to be
$$t _1 ^2 = {\mathbf q} _2, t _n ^2 = - {\mathbf q} _0 {\mathbf q} _1, t _n ( t _0 - t _0 ^{- 1} ) = {\mathbf q} _0 + {\mathbf q} _1 \quad \text{(c.f. Remark \ref{remHecke} 1))}.$$
Let $T _0, \ldots, T _n$ be the Iwahori-Matsumoto generators of $\mathbb H$ (c.f. \cite{M, L4}). Their Hecke relations read
$$( T _0 + 1 ) ( T _0 - t _0 ^2 ) = ( T _i + 1 ) ( T _i - t _1 ^2 ) = ( T _n + 1 ) ( T _n - t _n ^2 ) = 0,$$
where $1 \le i < n$. The natural map $\varphi ( T _i ) = T _{n - i}$ ($0 \le i \le n$) extends to an algebra map $\varphi : \mathbb H \rightarrow \mathbb H ^{\prime}$, where $\mathbb H ^{\prime}$ is the Hecke algebra of type $C _n ^{( 1 )}$ with parameters $t _n, t _1, t _0$. We have $t _n = \pm \sqrt{- {\mathbf q} _0 {\mathbf q} _1}$ and $t _0 = \pm \sqrt{- {\mathbf q} _0 / {\mathbf q} _1}$ or $\pm \sqrt{- {\mathbf q} _1 / {\mathbf q} _0}$. In particular, $\varphi$ changes the parameters as $( {\mathbf q} _0, {\mathbf q} _1, {\mathbf q} _2 ) \mapsto ( {\mathbf q} _0, {\mathbf q} _1 ^{- 1}, {\mathbf q}  _2 )$ or $( {\mathbf q}  _0 ^{- 1}, {\mathbf q}  _1, {\mathbf q}  _2 )$. Therefore, the representation theory of $\mathbb H _{a}$ ($a = ( s, \vec{q} )$) is unchanged if we replace $q _0$ with $q _0 ^{-1}$, or $q _1$ with $q _1 ^{-1}$.
\end{remark}

The rest of this section is devoted to the proof of Theorem \ref{BDL}. In the course of the proof, we use:

\begin{proposition}\label{exclorbits}
Let $a = ( s, \vec{q} ) \in {\bm T}$ be an admissible element. Let $\mathcal O \subset \mathfrak N$ be a $G$-orbit. For any two distinct $G ( s )$-orbits $\mathcal O _1, \mathcal O _2 \subset \mathcal O \cap \mathfrak N ^a _+$, we have
$$\overline{\mathcal O _1} \cap \mathcal O _2 = \emptyset.$$
\end{proposition}

\begin{proof}
By Proposition \ref{rEJNF} and Lemma \ref{stabTJ}, we deduce that the scalar multiplication of a normal form of $\mathfrak N$ is achieved by the action of $T$. It follows that each $G ( s )$-orbit of $\mathfrak N ^a$ is a $Z _{\bm G} ( a )$-orbit. Let $X \in \mathcal O _1$. Let $G _X$ be the stabilizer of $X$ in ${\bm G}$. Assume that $\mathcal O _2 \cap \overline{\mathcal O _1} \neq \emptyset$ to deduce contradiction. Since $\mathcal O _2$ is a $Z _{\bm G} ( a )$-orbit, we have $\mathcal O _2 \subset \overline{\mathcal O _1}$. Fix $X ^{\prime} \in \mathcal O _2$. Consider an open neighborhood $\mathcal U$ of $1$ in $G$ (as complex analytic manifolds). Then, $\mathcal U X ^{\prime} \in \mathcal O$ is an open neighborhood of $X ^{\prime}$. It follows that $\mathcal U X ^{\prime} \cap \mathcal O _1 \neq \emptyset$. We put $\mathfrak g _{a, X  ^{\prime}} := \mathrm{Lie} G _{X ^{\prime}} + \mathrm{Lie} Z _{\bm G} ( a )$. We have
$$N _{\mathcal O _2 / \mathcal O, X ^{\prime}} = \mathfrak g / \mathfrak g _{a, X ^{\prime}}.$$
Every non-zero vectors of $N _{\mathcal O _2 / \mathcal O, X ^{\prime}}$ is expressed as a linear combination of eigenvectors with respect to the $a$-action. These $a$-eigenvectors can be taken to have non-zero weights and does not contained in $G _{X ^{\prime}}$. It follows that
$$\mathcal U X ^{\prime} \cap \mathcal O _1 \not\subset \mathbb V ^a,$$
which is contradiction (for an arbitrary sufficiently small $\mathcal U$). Hence, we have necessarily $\mathcal O _2 \cap \overline{\mathcal O _1} = \emptyset$ as desired.
\end{proof}

\begin{proof}[Proof of Theorem \ref{EDL}]
By taking $G$-conjugate if necessary, we assume $a \in {\bm T}$. By Corollary \ref{cp}, it suffices to prove Theorem \ref{EDL} when $\mathcal C _a$ consists of a unique clan $[ 1, n ]$. By Corollary \ref{V10}, we can further assume $V _1 ^{(s, q _1)} = \{ 0 \}$ by swapping the roles of $q _0$ and $q _1$ if necessary. By Theorem \ref{sds} (c.f. Theorem \ref{Gd}), an admissible parameter $( a, X )$ is regular if there exists a simple $\mathbb H _a$-constituent of $M _{( a, X )}$ which does not appear in any $M _{(a, X ^{\prime})}$ such that $\overline{G ( s ) X} \subsetneq \overline{G ( s ) X ^{\prime}}$.

We apply Proposition \ref{reg} (if necessary) to modify $a$ so that the assumption of Proposition \ref{inclW} is fillfulled. By Proposition \ref{inclW}, each $M _{(a, X)}$ has a $W$-module structure given by the restriction of the $\mathbb H _a$-module structure. Moreover, the simple $W$-module $L _X$ corresponding to the $G$-orbit $GX \subset \mathfrak N$ (by the exotic Springer correspondence) appears in $M _{(a, X)}$. By Proposition \ref{exclorbits}, we have $G X \neq G X ^{\prime}$ for every $X ^{\prime} \in \mathfrak N ^a$ such that $\overline{G ( s ) X} \subsetneq \overline{G ( s )X ^{\prime}}$. By Corollary \ref{Wmod} and Corollary \ref{BM}, $M _{( a, X ^{\prime} )}$ does not contain $L _{X}$ as $W$-modules. Hence, the simple $\mathbb H _a$-constituent of $M _{(a, X)}$ which contains $L _X$ as $W$-type does not occur in any $M _{(a, X ^{\prime})}$ such that $\overline{G ( s ) X} \subsetneq \overline{G ( s ) X ^{\prime}}$ as required.
\end{proof}

\section{Consequences}\label{conq}
In this section, we present some of the consequences of our results. We retain the setting of the previous section.

\begin{definition}\label{stdmod}
Let $\nu = ( a, X )$ be an admissible parameter. Let $L _{\nu}$ be the simple module of $\mathbb H$ corresponding to $\nu$. Let $IC ( \nu )$ be the corresponding ${\bm G} ( a )$-equivariant simple perverse sheaf on $\mathfrak N _+ ^{a}$. (c.f. \S \ref{gen conv}) We denote by $P _{\nu}$ the projective cover of $L _{\nu}$ as $\mathbb H _{a}$-modules. (It exists since $\mathbb H _{a}$ is finite dimensional.)
\end{definition}

Let $K$ be a $\mathbb H$-module and let $L$ be a simple $\mathbb H$-module. We denote by $[ K : L ]$ the multiplicity of $L$ in $K$.

Applying \cite{CG} 8.6.23 to our situation, we obtain:

\begin{theorem}[The multiplicity formula of standard modules]\label{stdmlt}
Let $\nu = ( a, X )$, $\nu ^{\prime} = ( a, X ^{\prime} )$ be admissible parameters. We have:
$$[ M _{\nu} : L _{\nu ^{\prime}} ] = \sum _{k} \dim H ^k ( i _{X} ^! IC ( \nu ^{\prime} ) ) \text{ and } [ M ^{\nu} : L _{\nu ^{\prime}} ] = \sum _{k} \dim H ^k ( i _{X} ^* IC ( \nu ^{\prime} ) ),$$
where $i _X : \{ X \} \hookrightarrow \mathfrak N _+ ^{a}$ is an inclusion. \hfill $\Box$
\end{theorem}

The following result is a variant of the Lusztig-Ginzburg character formula of standard modules in our setting.

\begin{theorem}[The character formula of standard modules]
Let $\nu = ( a, X ) = ( s, \vec{q}, X )$ be an admissible parameter. Let $\mathfrak B _{\nu}$ be the set of connected components of $\mathcal E _X ^a$. For each $\mathcal B \in \mathfrak B _{\nu}$, we define a linear form $\left< \bullet , s \right> _{\mathcal B}$ as a composition map\\
\xymatrix@R=8.0pt{&\left< \bullet , s \right> _{\mathcal B} :\!\!\!\!\!\!\!\!\!\!\!\!\!\!\!\! & R ( T ) \ar[r] ^{\cong\quad } & R ( g B g ^{-1} ) \ar[r] ^{\quad \quad \mathrm{ev} _{s}} & \mathbb C\\
& & R ^+ \ar@{^{(}->}[u] \ar[r] & \{ \text{weights of } g B g ^{-1} \} \ar@{^{(}->}[u] & &
}\\
by some $g \in G$ such that $g B \in \mathcal B$. Then, $\left< \bullet, s \right> _{\mathcal B}$ is independent of the choice of $g$ and the restriction of $M _{\nu}$ to $R ( T )$ is given as
$$\mathrm{Tr} ( e ^{\lambda} ; M _{\nu} ) := \sum _{\mathcal B \in \mathfrak B _{\nu}} \left< \lambda, s \right> _{\mathcal B} \sum _{j \ge 0} \dim H _{2 j} ( \mathcal B, \mathbb C ).$$
\end{theorem}

\begin{proof}
Taking account into Theorem \ref{pi1 desc}, the proof is exactly the same as in \cite{CG} \S 8.2.
\end{proof}

\begin{definition}\label{decmat}
Let $a = ( s, \vec{q} ) \in {\bm T}$ be an admissible element. We form three $| \Lambda _{a} | \times | \Lambda _{a} |$-matrices
$$[ P : L ] ^{a} _{\nu, \nu ^{\prime}} := [ P _{\nu}, L _{\nu ^{\prime}}], D ^{a} _{\nu, \nu ^{\prime}} := \delta _{\nu, \nu ^{\prime}} \chi _c ( \nu ), \text{ and } IC ^{a} _{\nu, \nu ^{\prime}} := [ M ^{\nu}, L _{\nu ^{\prime}}],$$
where $\chi _c ( \nu ) := \sum _{i \ge 0} ( - 1 ) ^i \dim H ^i ( {\bm G} ( a ) X, \mathbb C)$ ($\nu = ( a, X )$).
\end{definition}

The following result is a special case of the Ginzburg theory \cite{CG} Theorem 8.7.5 applied to our particular setting:

\begin{theorem}[The multiplicity formula of projective modules]\label{BGG}
Keep the setting of Definition \ref{decmat}. We have
$$[ P : L ] ^{a} = IC ^{a} \cdot D ^a \cdot {} ^t IC ^{a},$$
where ${} ^t$ denotes the transposition of matrices. \hfill $\Box$
\end{theorem}

\begin{center}
{\bf Index of notation}\\
{\footnotesize (Sorted by the order of appearance)}
\end{center}
{\scriptsize
\begin{multicols}{3}
$G, B, T, G ( s ), U _{\alpha}, \ldots$ \hfill \S \ref{prem}\\

$R, R ^+, \mathbb E, \epsilon _i, \alpha _i$ \hfill \S \ref{prem}\\

$W, \dot{w} \in N _G ( T ), s _i, \ell$ \hfill \S \ref{prem}\\

${} ^w H := \dot{w} H \dot{w} ^{-1}$ \hfill \S \ref{prem}\\

$\mathsf{Stab} _H x$ ($x \in \mathcal X$) \hfill \S \ref{prem}\\

$\mathfrak g, \mathfrak t, \mathfrak g ( s ), \ldots$ \hfill \S \ref{prem}\\

$V [ \lambda ], V ^+, V ^-, \Psi ( V )$ \hfill  \S \ref{prem}\\

$H _{\bullet} ( \mathcal X ), H _{\bullet} ( \mathcal X, \mathbb Z)$ \hfill  \S \ref{prem}\\

$I, I ^*, \Gamma _0, \mathrm{exp}$ \hfill  \S \ref{prem}\\

$V _1 = \mathbb C ^n, V _2 = \wedge ^2 V _1$ \hfill \S \ref{exgeom}\\

$\mathbb V _{\ell}$: $\ell$-exotic rep. \hfill \S \ref{exgeom}\\

$F _{\ell}, \mu _{\ell}, \mathfrak N _{\ell}$ \hfill \S \ref{exgeom}\\

$F, \mu, \mathfrak N, \ldots$ \hfill \S \ref{exgeom}\\

$G _{\ell}, Z _{\ell}, p _i, \pi _{\ell}$ \hfill \S \ref{exgeom}\\

$\mathbb C _{a}$ \hfill \S \ref{exgeom}\\

$\mathsf p _w \in \mathsf O _w$ \hfill \S \ref{exgeom}\\

$\star, \circ$ \hfill \S \ref{exgeom}\\

$a _0 := ( 1, 1, -1, 1 )$ \hfill \S \ref{sJNF}\\

$\vec{q}, \log _i ( s )$ ($s \in T$) \hfill \S \ref{sJNF}\\

$\Lambda _a$ \hfill \S \ref{sJNF}\\

$x _i, y _{i,j} \in \mathbb V$ \hfill \S \ref{os}\\

$\mathbf J, T _J, \vec{\delta}$ \hfill \S \ref{os}\\

$\sigma = ( \mathbf J, \vec{\delta} )$ \hfill \S \ref{os}\\

$\mathbf v _{\sigma}, \mathbf v _{i, \sigma}, \mathbf v _{\sigma} ^J, \ldots$ \hfill \S \ref{os}\\

$\overline{\#} J, \underline{\#} J$ ($J \in \mathbf J$) \hfill \S \ref{os}\\

$T _{\ell}, F _{\ell} ^{a}, \nu ^{a} _{\ell}, \mathfrak N ^{a} _{\ell}, \ldots$ \hfill \S \ref{gen conv}\\

$\mathbf G = G _2, \mathbf T = T _2, \ldots$ \hfill \S \ref{fHecke}\\

$\mathcal A, \mathbb H$ \hfill \S \ref{fHecke}\\

$T _i, \mathbf q _i, e ^{\lambda} \in \mathbb H$ \hfill \S \ref{fHecke}\\

$Z _{\le w}, \mathbb O _i, \widetilde{T} _i, \ldots$ \hfill \S \ref{fHecke}\\

$\mathbb H _{a}, F ^{a} _+, \mu ^{a} _+, \mathfrak N ^{a} _+,\ldots$ \hskip -1mm \hfill \S \ref{fHecke}\\

$\mathfrak R _a$ \hfill \S \ref{fHecke}\\

$\mathbf c, n ^{\mathbf c}, \Gamma$ \hfill \S \ref{clan}\\

$\mathfrak g ( s ) _{\mathbf c}, G ( s ) _{\mathbf c}$ \hfill \S \ref{clan}\\

$\mathbb V ^{a}, \mathbb V ^{a} _{\mathbf c}, F _+ ^{a}, F _+ ^{a} ( w )$ \hfill \S \ref{clan}\\

${} ^w \mu ^{a} _{\mathbf c}$ \hfill \S \ref{clan}\\

$G _{\mathbf c}, \mathbb V ( \mathbf c ), X _{\mathbf c}, \ldots$ \hfill \S \ref{clan}\\

$\nu _{\mathbf c}$ \hfill \S \ref{clan}\\

$s _{\sigma}, D _{\sigma}, P _{\sigma}$ \hfill \S \ref{sss}\\

$\mathcal E _X, \mathcal E _X ^a, \mathsf{ch}$ \hfill \S \ref{avt}\\

$M _{\nu}, M ^{\nu}$ \hfill \S \ref{smi}\\

$s _Q, \mathbb V _Q, \mathbb H _Q, \ldots$ \hfill \S \ref{smi}\\

$L _{\sigma} = L _X$ ($X \in G \mathbf v _{\sigma}$) \hfill \S \ref{esc}\\

$ET, BA, H _{\bullet} ^A ( \mathcal X )$ \hfill \S \ref{dap}\\

$L _{\nu}, IC ( \nu )$ \hfill \S \ref{conq}\\



\end{multicols}}

{\footnotesize 
\begin{thebibliography}{NONAME}
\bibitem[AH08]{AH} Pramod N. Achar, and Anthony Henderson, Orbit closures in the enhanced nilpotent cone, Adv. Math. 219 (2008), 27--62
\bibitem[BH85]{BH} Hyman Bass and William J. Haboush, Linearizing certain reductive group actions, Trans. Amer. Math. Soc. {\bf 292} no.2 463--482 (1985).
\bibitem[BL94]{BL} Joseph Bernstein and Valery Lunts, Equivariant Sheaves and Functors, LNM 1578 Springer (1994).
\bibitem[BS49]{BS} Armand Borel, and Jean de Siebenthal, Les sous-groupes fermés de rang maximum des groupes de Lie clos. Comment. Math. Helv. 23, (1949). 200--221. 
\bibitem[Ca85]{C85} Roger W. Cartier, Finite Groups of Lie type: Conjugacy classes and complex characters, Wiley, (1985). ISBN 0-471-50683-4 
\bibitem[CG97]{CG} Neil Chriss, and Victor Ginzburg, Representation theory and complex geometry. Birkh\"auser Boston, Inc., Boston, MA, 1997. x+495 pp. ISBN 0-8176-3792-3
\bibitem[CK]{CK} Dan Ciubotaru, and Syu Kato, in praparation. 
\bibitem[DLP88]{DLP} Corrado De Concini, George Lusztig, and Claudio Procesi, Homology of the zero-set of a nilpotent vector field on a flag manifold, J. Amer. Math. Soc. {\bf 1}, (1988), no.1, 15--34.
\bibitem[DK85]{DK} Jiri Dadok, and Victor Kac, Polar representations. J. Algebra 92 (1985), no. 2, 504--524.
\bibitem[En06]{E} Naoya Enomoto, Classification of the Irreducible Representations of Affine Hecke Algebras of Type $B_2$ with unequal parameters, J. Math. Kyoto Univ. {\bf 46} no.2. (2006) 259--273.
\bibitem[En08]{E2} Naoya Enomoto, A quiver construction of symmetric crystals, 	arXiv:0806.3615, preprint.
\bibitem[FGT08]{FGT}Michael Finkelberg, Victor Ginzburg, and Roman Travkin, Mirabolic affine Grassmannian and character sheaves, arXiv:0802.1652.
\bibitem[Gi97]{Gi} Victor Ginzburg, Geometric methods in representation theory of Hecke algebras and quantum groups. Notes by Vladimir Baranovsky. NATO Adv. Sci. Inst. Ser. C Math. Phys. Sci., 514, Representation theories and algebraic geometry (Montreal, PQ, 1997), 127--183, Kluwer Acad. Publ., Dordrecht, 1998.
\bibitem[Ha77]{H} Robin Hartshorne, Algebraic Geometry, GTM 52. Springer-Verlag, 1977. xvi+496 pp. ISBN 0-387-90244-9.
\bibitem[Ig73]{Ig} Jun-ichi, Igusa, Geometry of absolutely admissible representations, in: Number theory, algebraic geometry and commutative algebra, in honor of Yasuo Akizuki, 373--452, Kinokuniya, Tokyo, (1973).
\bibitem[Jo98]{J1} Roy Joshua, Modules over convolution algebras from derived categories, I, J. Algebra, {\bf 203}, 385--446, (1998).
\bibitem[Ka06b]{K} Syu Kato, An exotic Springer correspondence for symplectic groups, math.RT/0607478, preprint (never to appear).
\bibitem[Ka08]{K08} Syu Kato, Deformations of nilpotent cones and Springer correspondences, arXiv:0801.3707, preprint.
\bibitem[KL87]{KL2} David Kazhdan, and George Lusztig, Proof of the Deligne-Langlands conjecture for Hecke algebras, Invent. Math. 87 (1987), no. 1, 153--215.
\bibitem[Lu88]{L1} George Lusztig, Cuspidal local systems and graded Hecke algebras. I, Inst. Hautes \'Etudes Sci. Publ. Math., {\bf 67}, (1988), 145--202.
\bibitem[Lu89]{L2} George Lusztig, Affine Hecke algebras and their graded version, J. Amer. Math. Soc. {\bf 2}, (1989), no.3, 599--635.
\bibitem[Lu95a]{L5} George Lusztig, Classification of unipotent representations of simple p-adic groups. Internat. Math. Res. Notices 1995, no. 11, 517--589.
\bibitem[Lu95b]{L3} George Lusztig, Cuspidal local systems and graded Hecke algebras. II, Representations of groups (Banff, AB, 1994), CMS Conf. Proc., {\bf 16}, 217--275. (1995)
\bibitem[Lu02]{L6} George Lusztig, Cuspidal local systems and graded Hecke algebras. III, Representations Theory {\bf 6}, 202--242. (2002)
\bibitem[Lu03]{L4} George Lusztig, Hecke algebras with unequal parameters, CRM Monograph Series, {\bf 18}, American Mathematical Society, Providence, RI, (2003), vi+136, ISBN: 0-8218-3356-1
\bibitem[Mc03]{M} Ian G. Macdonald, Affine Hecke algebras and orhogonal polynomials, Cambridge Tracts in Mathematics 157 (2003)
\bibitem[Oh86]{Oh} Takuya Ohta, The singularities of the closures of nilpotent orbits in certain symmetric pairs, T\^ohoku Math. J., 38 (1986), 441-468.
\bibitem[OS07]{OS} Eric Opdam, and Maarten Solleveld, Homological algebra for affine Hecke algebras, arXiv:0708.1372, to appear in Adv. in Math.
\bibitem[OS08]{OS2} Eric Opdam, and Maarten Solleveld, Discrete series characters for affine Hecke algebras and their formal degrees, arXiv:0804.0026.
\bibitem[Po04]{P} Vladimir L. Popov, The cone of Hilbert nullforms, Proc. Steklov Math. Inst. 241 (2003), 177-194.
\bibitem[Ra01]{R} Arun Ram, Representations of rank two affine Hecke algebras, in "Advances in Algebra and Geometry, University of Hyderabad conference 2001", C. Musili ed., Hindustan Book Agency, 2003, 57-91, available online at {\tt http://www.math.wisc.edu/{\textasciitilde}ram/bib.html}
\bibitem[Sa88]{S} Morihiko Saito, Modules de Hodge polarisables, Publ. Res. Inst. Math. Sci. {\bf 24} (1988), no. 6, 849--995 (1989).
\bibitem[Sc78]{Sc} Gerald W. Schwarz, Representations of simple groups with a free module of covariants, Invent. Math. 49, 1--12 (1978).
\bibitem[Se84]{Se} Jiro Sekiguchi, The nilpotent subvariety of the vector space associated to a symmetric pair, Publ. Res. Inst. Math. Sci. {\bf 20} (1984), 155--212.
\bibitem[So07]{So} Maarten Solleveld, Periodic cyclic homology of reductive p-adic groups, 	arXiv:0710.5815.
\bibitem[Sp07]{Sp} Tonny A. Springer, The exotic nilcone of a symplectic group, preprint.
\bibitem[Sp82]{Spal} Nicolas Spaltenstein, Classes unipotents et sous groupes de Borel, Lect. Notes. Math. 946 Springer-Verlag, 1982.
\bibitem[Th86]{Th} Robert W. Thomason, Lefschetz-Riemann-Roch theorem and coherent trace formula, Invent. Math. 85 515--543 (1986).
\bibitem[Tr08]{Tr} Roman Travkin, Mirabolic Robinson-Schensted-Knuth correspondence, preprint.
\end{thebibliography}}
\end{document}